\documentclass[12pt]{amsart}
\usepackage{amssymb}
\usepackage{amsbsy}
\usepackage{amscd}
\usepackage[mathscr]{eucal}
\usepackage{verbatim}

\include{diagrams}
\usepackage{pstricks}

\makeatletter
\def\cal{\mathcal}
\def\Bbb{\mathbb}
\def\frak{\mathfrak}

\newenvironment{pf*}[1]{\proof[#1]}{\endproof}
\makeatother
\newenvironment{aenume}{%
  \begin{enumerate}%
  }{\end{enumerate}}
\makeatletter
\@ifclasslater{amsart}{1999/11/24}{}{
\renewcommand*\subjclass[2][1991]{%
  \def\@subjclass{#2}%
  \@ifundefined{subjclassname@#1}{%
    \ClassWarning{\@classname}{Unknown edition (#1) of Mathematics
      Subject Classification; using '1991'.}%
  }{%
    \@xp\let\@xp\subjclassname\csname subjclassname@#1\endcsname
  }%
}
\renewcommand{\subjclassname}{%
  \textup{1991} Mathematics Subject Classification}
\@xp\let\csname subjclassname@1991\endcsname \subjclassname
\@namedef{subjclassname@2000}{%
  \textup{2000} Mathematics Subject Classification}
}
\makeatother
\newenvironment{NB}{
{\bf NB}. \footnotesize
}{}
\renewenvironment{NB}{%
 \comment
  }{\endcomment}
\newtheorem{Theorem}[equation]{Theorem}
\newtheorem{Corollary}[equation]{Corollary}
\newtheorem{Lemma}[equation]{Lemma}
\newtheorem{Proposition}[equation]{Proposition}

\theoremstyle{definition}
\newtheorem{Definition}[equation]{Definition}

\newtheorem{Notation}[equation]{Notation}

\theoremstyle{remark}
\newtheorem{Remark}[equation]{Remark}

\numberwithin{equation}{section}

\newcommand{\thmref}[1]{Theorem~\ref{#1}}
\newcommand{\secref}[1]{\S\ref{#1}}
\newcommand{\lemref}[1]{Lemma~\ref{#1}}
\newcommand{\propref}[1]{Proposition~\ref{#1}}
\newcommand{\corref}[1]{Corollary~\ref{#1}}
\newcommand{\subsecref}[1]{\S\ref{#1}}

\newcommand{\remref}[1]{Remark~\ref{#1}}
\newcommand{\C}{{\Bbb C}}
\newcommand{\Z}{{\Bbb Z}}
\newcommand{\Q}{{\Bbb Q}}
\newcommand{\R}{{\Bbb R}}
\newcommand{\proj}{{\Bbb P}}
\newcommand{\CP}{\proj}

\newcommand{\End}{\operatorname{End}}
\newcommand{\Hom}{\operatorname{Hom}}
\newcommand{\Ext}{\operatorname{Ext}}

\newcommand{\len}{\mathop{\text{\rm len}}\nolimits}
\newcommand{\rank}{\operatorname{rank}}

\newcommand{\tr}{\operatorname{tr}}

\newcommand{\id}{\operatorname{id}}
\newcommand{\ve}{\varepsilon}

\newcommand{\linf}{\ell_\infty}
\newcommand{\shfO}{\mathcal O}
\newcommand{\bp}{{\widehat\proj}^2}
\newcommand{\bM}{{\widehat M}}

\newcommand{\Pic}{\operatorname{Pic}}

\newcommand{\Supp}{\operatorname{Supp}}
\newcommand{\ch}{\operatorname{ch}}

\newcommand{\Todd}{\operatorname{Todd}}

\newcommand{\Li}{\operatorname{Li}}
\newcommand{\Zin}{Z^{\text{\rm inst}}}

\newcommand{\Fper}{F^{\text{\rm pert}}}
\newcommand{\Fin}{F^{\text{\rm inst}}}

\newcommand{\Finz}{\mathcal F^{\text{\rm inst}}}
\newcommand{\q}{\mathfrak q}

\newcommand{\hT}{\widetilde T}

\newcommand{\res}{\mathop{\text{\rm res}}}
\def\oo{{\cal O}}
\def\I{{\cal I}}
\def\AA{{\cal A}}
\def\A{{\mathbb A}}
\def\P{{\mathbb P}}
\def\<{\langle}
\def\>{\rangle}

\def\F{{\mathcal F}}

\def\E{{\mathcal E}}
\def\cho{\overline{\ch}}

\newcommand{\TT}{\Gamma}

\begin{document}
\title[Donaldson invariants via instanton counting]
{Instanton counting and Donaldson invariants}
\author{Lothar G\"ottsche}
\address{International Centre for Theoretical Physics, Strada Costiera 11, 
34014 Trieste, Italy}
\email{gottsche@ictp.trieste.it}
\author{Hiraku Nakajima}
\address{Department of Mathematics, Kyoto University, Kyoto 606-8502,
Japan}
\email{nakajima@math.kyoto-u.ac.jp}
\thanks{The second author is supported by the Grant-in-aid
for Scientific Research (No.15540023, 17340005), JSPS}

\author{K\={o}ta Yoshioka}
\address{Department of Mathematics, Faculty of Science, Kobe University,
Kobe 657-8501, Japan}
\email{yoshioka@math.kobe-u.ac.jp}
\subjclass[2000]{Primary 14D21; Secondary 57R57, 81T13, 81T60}


\begin{abstract}
For a smooth projective toric  surface we determine the Donaldson 
invariants and their wallcrossing
in terms of the Nekrasov partition function. 
Using the solution of the Nekrasov conjecture
\cite{NY1,NO,BE} and its refinement \cite{NY2}, we apply this result to give a generating function for the wallcrossing of Donaldson invariants 
of good walls of simply connected projective surfaces with $b_+=1$ in terms of modular forms. This formula was proved earlier in \cite{G} more generally
for simply connected $4$-manifolds with $b_+=1$, assuming the 
Kotschick-Morgan conjecture and it was also derived by physical arguments in \cite{MW}.
\end{abstract}

\maketitle
\section*{Introduction} 
Donaldson invariants have for a long time played an important r\^ole in
the study and classification of differentiable $4$-manifolds (see
\cite{DK}).
They are defined by moduli spaces of anti-self-dual connections on a
principal $SO(3)$-bundle. The anti-self-duality equation depends on
the choice of a Riemannian metric $g$. For generic $g$ there are no
reducible solutions to the equation and moduli spaces are smooth
manifolds.
In case $b_+ > 1$ two generic Riemannian metrics can be connected by 
a path. Then Donaldson invariants are independent of the choice of the
metric, and they are invariants of a $C^\infty$ compact oriented
$4$-manifold $X$.

On the other hand, in case $b_+ = 1$ nongeneric metrics form a real
codimension $1$ subset in the space of Riemannian metrics, i.e.\ a
collection of {\it walls\/}, and two generic metrics cannot be
connected by a path in general. As a consequence, Donaldson invariants
are only piecewise constants as functions of the Riemannian metric 
$g$ \cite{K,KM}. More precisely we have a chamber structure on the period
domain, which is a connected component $\mathcal C$ of the positive
cone in the second cohomology group $H^2(X,\R)$, and the
Donaldson
invariants stay constant only when the period $\omega(g)$, which is
the cohomology class of the self-dual harmonic $2$-form modulo
scalars, stays in a chamber.
The {\it wallcrossing terms\/} are  the differences of Donaldson invariants
when the metric moves to another chamber passing through a wall.
In \cite{G} the first author gave a formula for their generating
function in terms of modular forms, assuming the Kotschick-Morgan
conjecture\footnote{There are two preprints by Chen \cite{Chen} and by
Feehan-Leness \cite{FL4}, giving a proof and an announcement of a proof
of the conjecture respectively. 
Fr\o yshov also gave a talk on a proof.
Their approaches are differential geometric and quite different
from ours, and the authors believe they are correct, but  unfortunately do not have
the ability to check their papers in full detail.}, which
states that the wallcrossing term is a polynomial in the intersection
form and the multiplication by $\xi$, the cohomology class defining
the wall (see \subsecref{backDon} for more detail).
The method of the proof was indirect and did not give a clear reason why
modular forms appear.

A physical derivation of the wallcrossing formula was given by
Moore-Witten \cite{MW}. We shall review their derivation and the
physical background  only very briefly here (see
\cite[Introduction]{NY2} for a more detailed exposition for
mathematicians).
The work of Moore-Witten was based on Seiberg-Witten's ansatz \cite{SW}
of the $\mathcal N=2$ supersymmetric Yang-Mills theory on $\mathbb
R^4$, which is a physical theory underlying Donaldson invariants
\cite{Wit}. The theory is controlled by a family of elliptic curves
parametrized by a complex plane (called the $u$-plane). The modular
forms that appear in the wallcrossing formula are related to this
family.
They
expressed Donaldson invariants in terms of two contributions, the integral
over the $u$-plane and the contribution from the points $\pm 2$, where the
corresponding elliptic curves are singular.
The latter contribution corresponds to Seiberg-Witten invariants,
which conjecturally contain the same information as Donaldson
invariants \cite{Wit2}.
Moore-Witten further studied the $u$-plane integral and its
contribution to Donaldson invariants. They recovered the wallcrossing
formula, as well as Fintushel-Stern's blowup formula \cite{FS}, and
also obtained new results, such as Seiberg-Witten contributions and
calculation for $\mathbb P^2$ in terms of Hurwitz class numbers.

Seiberg-Witten and Moore-Witten's arguments clarified the reason why
modular forms appear in Donaldson invariants. But they were physical
and have no mathematically rigorous justification so far.
A more rigorous approach was proposed much later by Nekrasov \cite{Nek}.
He introduced the partition function 
\[
  \Zin(\ve_1,\ve_2,a;\Lambda) = \sum_{n\ge 0} \Lambda^{4n}\int_{M(n)} 1,
\]
where $M(n)$ is the Gieseker's partial compactification of the framed
moduli space of $SU(2)$-instantons on $\mathbb A^2$ and $\int_{M(n)}$
denotes the pushforward homomorphism to a point in the equivariant
homology groups, defined by a formal application of
\begin{NB} removed "the" 6.9. LG \end{NB}
Bott's fixed point
formula to the noncompact space $M(n)$. The variables $\ve_1$, $\ve_2$
are generators of the equivariant cohomology
$H^*_{\C^*\times\C^*}(\mathrm{pt})$ of a point with respect to the two
dimensional torus $\C^*\times\C^*$ acting on $\mathbb A^2$. The
remaining variable $a$ is also a generator of
$H^*_{\C^*}(\mathrm{pt})$, where $\C^*$ acts on $M(n)$ by the change
of the framing.
This definition can be viewed as the generating function of the
equivariant Donaldson invariants of $\mathbb R^4 = \mathbb
A^2$. Although Nekrasov was motivated by a physical argument, the partition
function is mathematically rigorously defined. 
He then conjectured
\begin{equation*}
   \ve_1\ve_2 \log \Zin(\ve_1,\ve_2,a;\Lambda) = \Finz_0(a,\Lambda)
   + \text{higher terms}
   \qquad \text{as $\ve_1,\ve_2\to 0$},
\end{equation*}
where $\Finz_0(a,\Lambda)$ is the instanton part of the Seiberg-Witten
prepotential defined via periods of the elliptic curves mentioned
above. The conjecture was proved by three groups, the second and third
named authors \cite{NY1}, Nekrasov-Okounkov \cite{NO}, and
Braverman-Etingof \cite{BE} by completely different methods.

In this paper we express the wallcrossing terms of Donaldson
invariants in terms of the Nekrasov partition function, under the
assumption that the wall is good (see \subsecref{wallmod} for the
definition). Thereby we give a partial mathematical justification of
Moore-Witten's argument, where Seiberg-Witten's ansatz is replaced by
the Nekrasov partition function. More precisely, we take a smooth
toric surface $X$ and consider equivariant Donaldson invariants. They
also depend on the choice of a Riemannian metric as ordinary Donaldson
invariants, and we have an equivariant wallcrossing term. The first
main result (\thmref{main1}) expresses it as the residue at $a=\infty$
(corresponding to $u=\infty$ of the $u$-plane) of a product over
contributions from fixed points in $X$, and the local contribution is
essentially the Nekrasov partition function. 
This result comes from the following: In the wallcrossing
the moduli space changes by replacing certain sheaves lying in extensions
of ideal sheaves of zero-dimensional schemes twisted by line bundles by 
extensions the other way round. Using this fact one can express the 
change of Donaldson invariants under wallcrossing in terms of intersection 
numbers  on the  Hilbert schemes $X_2^{[l]}$ of points on two copies of $X$.
For the wallcrossing of 
the Donaldson invariants without higher Chern characters this was already 
shown in \cite[Th.~6.13]{EG1} and \cite[Th.~5.4,\ Th.~5.5]{FQ}.
These intersection numbers can be computed via equivariant localization on 
$X_2^{[l]}$. Every $\Gamma$-invariant scheme in $X_2^{[l]}$ is a
union of $\Gamma$-invariant
schemes with support one of the fixed points of $X$,
and the contribution to the intersection number coming from invariant 
subschemes with support one of the 
fixed points of $X$ is given by the Nekrasov partition function.

Then the second main result (\thmref{main2}) is about the
nonequivariant limit $\ve_1,\ve_2\to 0$ and we recover the formula in
\cite{G} via the solution of Nekrasov's conjecture and its refinement
\cite{NY2}, i.e.\ determination of several higher terms of
$\varepsilon_1\varepsilon_2 \log
Z(\varepsilon_1,\varepsilon_2,a;\Lambda)$.
It is worthwhile remarking that the variable $a$ appears in the
wallcrossing term as an auxiliary variable, which is eventually
integrated out. By contrast it plays a fundamental role in the
Seiberg-Witten ansatz as a period of the Seiberg-Witten curve.

It is natural to expect that our equivariant wallcrossing formula is a
special case of that for the Donaldson invariants for {\it families\/}
whose definition was mentioned in \cite{Don}. Then we expect that
higher coefficients of the Nekrasov partition function, which are
higher genus Gromov-Witten invariants for a certain noncompact toric
Calabi-Yau 3-fold, also play a role in $4$-dimensional topology.

In \secref{cobor} we show that the wallcrossing term for a good wall of
an arbitrary projective surface $X$ can be given by a universal
polynomial depending on Chern classes $c_i(X)$, $\xi$ and the intersection
product on $H^*(X)$. The proof of this result does not yield an
explicit form of the universal polynomial directly. But combining with
the explicit form obtained for toric surfaces, we conclude that the
same explicit formula holds for an  arbitrary surface with $b_+=1$. In
particular, it does not depend on $c_1(X)$ and satisfies the
statements in the Kotschick-Morgan conjecture. (See \remref{rem:c_1}
for more explanation.) 
The `goodness' of the wall means that the moduli space is smooth along
sheaves replaced by the wallcrossing. 

Results of Mochizuki show that  the goodness assumption  can be removed: 
\cite[Thm 1.12]{Moc} gives 
\propref{wallcr}  for arbitrary walls if we replace vector
bundles $\mathcal A_{\xi,+}$, $\mathcal A_{\xi,-}$ by the
corresponding classes in $K$-theory. In the proof Mochizuki uses virtual fundamental classes and virtual localization. Therefore our main results
(\thmref{main2}, \corref{cor:main3}) are true for any wall on a
simply-connected projective surface.
\begin{NB} This section was changed 6.9 LG \end{NB}

In \secref{P2} we express the equivariant Donaldson invariants
themselves for $\proj^2$, instead of the wallcrossing terms, in terms
of the Nekrasov partition function. The result here is independent of
those in previous sections.
However we do not know how to deduce an explicit formula for ordinary
Donaldson invariants via nonequivariant limit $\ve_1,\ve_2\to 0$.
Note also that we cannot extend this result to other toric surfaces,
as fixed points are no longer isolated.

The Nekrasov partition function is defined for any rank. A higher rank
generalization of Donaldson invariants is given recently by Kronheimer
\cite{Kr}. Though they are defined for $b_+ > 1$, many of his results
are applicable to the $b_+ = 1$ case also. Therefore it is natural to
hope that our results can be generalized to the higher rank cases. One
of new difficulties appearing in higher rank cases is a recursive
structure of the wallcrossing. We hope to come back this problem in
future.

Finally let us mention that Nekrasov proposed that the equivariant
Donaldson invariants for toric surfaces can be expressed as products
of his partition functions over fixed points, integrated over $a$ in
any rank \cite{Nek2}. As equivariant Donaldson invariants vanish for a
certain chamber for toric surfaces, our wallcrossing formula gives
such an expression together with an explicit choice of contour for the
$a$-integral, which was not specified in [loc.\ cit.]. It is an
interesting problem to justify his argument more directly.

\subsection*{Acknowledgement}\begin{NB} corrected 6.9.LG\end{NB}
The project started in 2004 Jan.\ when the first author visited Kyoto
for a workshop organized by the second and third authors. They are
grateful to the Kyoto University for its hospitality.
A part of this paper was written when all authors visited MSRI at
Berkeley. They thank MSRI for its hospitality.
They also thank Nikita Nekrasov for discussions over years and Takuro
Mochizuki for explanations of his results and Jun Li for useful discussions.
\begin{NB} added 6.9 LG \end{NB}

%


\tableofcontents

\section{Background Material}
We will  work over $\C$. We usually consider homology and cohomology with
rational coefficients and for a variety $Y$ we will 
 write $H_i(Y)$, and $H^i(Y)$ for $H_i(Y,\Q)$ and
 $H^i(Y,\Q)$ respectively. If $Y$ 
is projective and $\alpha \in H^*(Y)$, we denote 
$\int_Y\alpha$ its evaluation on the fundamental cycle of $Y$.
If $Y$ carries an action of a torus $T$, $\alpha$ is a $T$-equivariant class, 
and $p:X\to pt$ is the projection to a point, we denote  $\int_Y\alpha:=p_*(\alpha)\in H^*(pt)$.

In this whole paper $X$ will be a nonsingular projective surface over $\C$. Later 
we will specialize $X$ to a smooth projective toric surface.
For a class $\alpha\in H^*(X)$, we denote $\<\alpha\>:=\int_X\alpha$.
If $X$ is a toric surface, we use the same notation for the equivariant
pushforward to a point.

\subsection{Donaldson invariants}\label{backDon}
Let $X$ 
be a smooth simply connected  compact oriented $4$-manifold with a 
Riemannian metric $g$. For $P\to X$ an $SO(3)$-bundle over $X$ let 
$M(P)$ be the 
moduli space of irreducible anti-self-dual 
connections on $P$. 
For generic $g$ this will be a manifold of dimension
$d:=-2p_1(P)-3(1+b^+(X))$. Let ${\cal P}\to X\times M(P)$ be the universal 
bundle.
Then the Donaldson invariant of $Y$ is a polynomial on $H_0(X)\oplus H_2(X)$,
defined by 
$$
D^g_{c_1,d}(\alpha^n p^b)=\int_{M(P)}\mu(\alpha)^n\mu(p)^b.
$$
Here $c_1$ is a lift of $w_2(P)$ to $H^2(X,\Z)$, 
$p\in H_0(X)$ is the class of a point and  $\alpha\in H_2(X)$, and for $\beta\in H_i(X)$ we define $\mu(\beta):=-\frac{1}{4} p_1({\cal P})/\beta$.
As $M(P)$ is not compact, this integral must be justified using the Uhlenbeck compactification of $M(P)$.
Note that the orientation of $M(P)$ depends on the lift $c_1$ and
 a choice of a connected
component ${\cal C}$ of the positive cone  in $H^2(X,\R)$ 
which for algebraic surfaces we 
always take to be the component containing the ample cone.
The generating function is
$$D^g_{c_1}(\exp(\alpha z+px)):=\sum_{d\ge 0}\sum_{n,m\ge 0} D_{c_1}^g 
\Big(\frac{\alpha^n}{n!}\frac{p^m}{m!}\Big)z^nx^m.$$
When $b_+(X)>1$, then $D^g_{c_1,d}$ is independent of $g$ as long as $g$ is
generic. 
If $b_+(X)=1$, then $D^g_{c_1,d}$ depends on the period
point $\omega(g)\in {\cal C}$. 

In fact the positive cone in $H^2(X,\R)$ has 
a chamber structure (see \cite{K},\cite{KM}):  
For a class $\xi\in H^2(X,\Z)\setminus \{0\}$, 
\begin{NB} corrected 6.9.LG, we do not want $0$ to define a wall.
\end{NB}
we put $W^\xi:=\big\{x\in {\cal C}\bigm|
\<x\cdot \xi\>=0\big\}$. Assume $W^\xi\ne \emptyset$. Then we call $\xi$ a class of 
type $(c_1,d)$ and  call $W^\xi$ 
a wall of type $(c_1,d)$, if 
the following conditions hold
\begin{enumerate}
\item
$\xi+c_1$ is divisible by $2$ in $H^2(X,\Z)$,   
\item $d+3+\xi^2\ge 0$. 
\end{enumerate}
We call $\xi$ a class of type $c_1$ and call $W^\xi$ a wall of type $c_1$, 
if $\xi+c_1$ is divisible by $2$ in $H^2(X,\Z)$.
The chambers of type $(c_1,d)$ are the connected components of 
the complement of the union of all walls of type $(c_1,d)$ in ${\cal C}$.
In \cite{KM} it is shown that $D^g_{c_1,d}$ depends only on the chamber of $\omega(g)$.
\begin{NB}
  The case $\xi = 0$ is not excluded in the definition of the wall.
  But it never appear in the wallcrossing formula. Therefore it does
  not give a wall in practice. May 31, HN
\end{NB}

Let $C_+$, $C_-$ be chambers of type $(c_1,d)$ in $\mathcal C$ and
$g_+$, $g_-$ be Riemannian metrics with $\omega(g_\pm)\in C_\pm$. Then
\begin{equation*}
   D^{g_+}_{c_1,d}(\alpha^n p^b) - D^{g_-}_{c_1,d}(\alpha^n p^b)
   = \sum_\xi \Delta_{\xi,d}^X(\alpha^n p^b),
\end{equation*}
where the summation runs over the set of all classes $\xi$ of type
$(c_1,d)$ with $\langle\xi\cdot C_+\rangle > 0 > \langle\xi\cdot
C_-\rangle$.  The term $\Delta_{\xi,d}^X$ is called the {\it wallcrossing
  term}. The Kotschick-Morgan \cite{KM} conjecture says that
$\Delta_{\xi,d}^X$ is a polynomial in the multiplication by $\xi$ and the
intersection form with coefficients depending only on $\xi^2$ and the
homotopy type of $X$. Wallcrossing terms with small $d$ had been
calculated by various authors \cite{K,KM,EG1,EG2,FQ,KL,MaW}. Then the
first named author \cite{G} gave a formula for the generating function
of $\Delta_\xi^X$ in terms of modular forms, assuming the
Kotschick-Morgan conjecture. See also \cite{GZ}.

Now we specialize to the case of a smooth projective surface $X$ with $p_g(X)=0$, in particular $b_+(X)=1$. Let $H$ be an ample divisor on $X$. Then the cohomology
class $H$ is a representative of the period point of the Fubini-Study metric of $X$ 
associated to $H$. We write $D_{c_1,d}^H$ for the corresponding Donaldson invariants.
By \cite{Li},\cite{Mo},
the $D_{c_1,d}^H$ can also be computed
using moduli spaces of sheaves on $X$.
We denote by $M_H^X(c_1,d)$  the moduli space of torsion-free $H$-semistable 
sheaves
(in the sense of Gieseker and Maruyama) of rank $2$ and with $c_1(E)=c_1$ and
$4c_2(E)-c_1(E)^2-3=d$. 
 Let $M^X_H(c_1,d)_s$ be the open 
subset of stable sheaves.
Assume that $M_H^X(c_1,d)=M_H^X(c_1,d)_s$ and that there exists a
universal sheaf ${\E}$ on $X\times M^X_H(c_1,d)$.
\begin{NB} corrected 7.9 LG\end{NB}
If there is no universal sheaf, we can replace it by a quasiuniversal sheaf.
When $p_g=0$ (the case of our primary interest), then $\operatorname{Pic}(X) \to H^2(X,{\Bbb Z})$ is surjective,
which means that $\chi(*,*)$ is unimodular on $K(X)$.
Hence there is a universal sheaf, if $M_H^X(c_1,d)=M_H^X(c_1,d)_s$.
\begin{NB} corrected 7.9. LG\end{NB}
For $\beta\in H_i(X,\Q)$, we put 
$\mu(\beta):=\big(c_2({\E})-\frac{1}{4}c_1({\E})^2\big)/\beta\in 
H^{4-i}(M_H^X(c_1,d),\Q),$
and define 
$$\Phi_{c_1,d}^H(\alpha^n p^m):=\int_{M_H^X(c_1,d)}\mu(\alpha)^n\mu(p)^m$$
and 
\begin{equation}\label{donn}\Phi_{c_1}^H(\exp(\alpha z+px)):=\sum_{d\ge 0}\Lambda^d \sum_{m,n}
\Phi^H_{c_1,d}\Big(\frac{\alpha^n}{n!}
\frac{p^m}{m!}\Big)z^n x^m=\sum_{d\ge 0} \Lambda^d
\int_{M_H^X(c_1,d)} \exp(\mu(\alpha z+px)).\end{equation}
Here if $Y$ is a compact variety  and $f=\sum_{i,j} a_{i,j} x^iz^j\in H^*(Y)[[x,z]]$,
we write 
$\int_Y f=\sum_{i,j} x^iz^j\int_Y a_{i,j}$.
Assume that $M_H^X(c_1,d)$ has the expected dimension $d$
or is empty,
and that $H$ does not lie on a wall of type $(c_1,d)$. Then 
by the results of \cite{Mo},\cite{Li} one has
\begin{equation}
\label{eq:Don=Phi}
\Phi_{c_1,d}^H(\alpha^n p^m)=(-1)^{(c_1^2+\<c_1\cdot K_X\>)/2}D_{c_1,d}^H(\alpha^n p^m).
\end{equation}

When $M^X_H(c_1,d)$ is not necessary of expected dimension, we define
the invariants as follows (cf.\ \cite[\S3.8]{FM}): we consider blowup
$P\colon \widehat X\to X$ at sufficiently many points $p_1,\dots,p_N$
disjoint from cycles representing $\alpha$, $p$.
Let $C_1, \dots, C_N$ denote the exceptional curves. We consider the
moduli space $M^{\widehat X}_{P^*H}(P^*c_1,d+4N)$, where the
polarization `$P^*H$' means $P^*H - \ve C_1 - \ve C_2 - \cdots \ve
C_N$ for sufficiently small $\ve > 0$. Then it has expected dimension
for sufficiently large $N$ by \secref{sec:app}.
We define
\begin{multline*}
   \int_{M^X_H(c_1,d)} \exp(\mu(\alpha z + p x))
\\
   := (-\frac12)^N \int_{M^{\widehat X}_{P^*H}(P^*c_1,d+4N)}
   \mu(C_1)^4\cdots \mu(C_N)^4 
   \exp(\mu(\alpha P^*z + p P^*x)).
\end{multline*}
By the blowup formula (see \cite[Th.~8.1]{FM}), this definition is
independent of $N$. From its definition, \eqref{eq:Don=Phi} remains to
hold. 

\subsection{Nekrasov partition function}\label{subsec:Nek}
We briefly review the Nekrasov partition function in the case of rank $2$.
For more details see \cite[sections 3.1, 4]{NY2}.
Let $\ell_\infty$ be the line at infinity in $\P^2$.
Let $M(n)$ be the moduli space of
pairs $(E,\Phi)$, where $E$ is a rank $2$ torsion-free sheaf on $\P^2$ with $c_2(E)=n$,
which is locally free in a neighbourhood of $\ell_\infty$ and 
$\Phi:E|_{\ell_\infty}\to \oo_{\ell_\infty}^{\oplus 2}$ is an isomorphism.
$M(n)$ is a nonsingular quasiprojective variety of dimension $4n$.

Let $\Gamma:=\C^*\times \C^*$ and $\widetilde T:=\Gamma\times \C^*$.
$\widetilde T$ acts on $M(n)$ as follows:
For $(t_1,t_2)\in \Gamma$, let $F_{t_1,t_2}$ be the automorphism of 
$\P^2$ defined by 
$F_{t_1,t_2}([z_0,z_1,z_2])\mapsto[z_0,t_1z_1,t_2z_2]$, and for $e\in \C^*$ let
$G_e$ be the automorphism of $\oo_{\ell_\infty}^{\oplus 2}$ given by 
$(s_1,s_2)\mapsto (e^{-1}s_1,es_2)$.
Then for $(E,\Phi)\in M(n)$ we put
$(t_1,t_2,e)\cdot(E,\Phi):=\big((F_{t_1,t_2}^{-1})^*E,\Phi'\big)$, where
$\Phi'$ is the composition 
\begin{diagram}[height=.8em,width=2em,scriptlabels]
(F^{-1}_{t_1,t_2})^*(E)|_{\ell_\infty}&\rTo^{(F^{-1}_{t_1,t_2})^*\Phi}& (F^{-1}_{t_1,t_2})^*\oo_{\ell_\infty}^{\oplus 2}&
\rTo{} &\oo_{\ell_\infty}^{\oplus 2}&\rTo^{G_e}& \oo_{\ell_\infty}^{\oplus 2}
\end{diagram}
where the middle arrow is the homomorphism given by the action.
%
Let $\ve_1,\ve_2,a$ be the coordinates on the Lie algebra of $\widetilde T$,
i.e. we can write $(t_1,t_2,e)=(e^{\ve_1},e^{\ve_2},e^a)$.

We briefly recall equivariant integration in the form we want to use it. 
Let $T$ be torus acting on a nonsingular variety $Y$ with finitely many fixed points
$q_1,\ldots,q_s$. Let $e_1,\ldots,e_n$ be the coordinates on the Lie algebra
of $T$. The equivariant cohomology of a point is 
$H^*_T(pt)=\Q[e_1,\ldots,e_n]$. 
If $\alpha\in H^*_T(Y)$ is an equivariant cohomology class, then we put
$$
\int_Y \alpha:=\sum_{i=1}^s \frac{\iota_{q_i}^*(\alpha)}{e_T(T_{q_i}Y)}\in 
\Q(e_1,\ldots,e_n).$$
Here $\iota_{q_i}^*$ is the equivariant pullback via the embedding $q_i\hookrightarrow Y$.
and $e_T(T_{q_i}Y)$ is the equivariant Euler class of the tangent space of $Y$ at $q_i$.
If $Y$ is also compact, then $\int_Y$ is the usual pushforward to a point 
in equivariant cohomology, in particular $\int_Y \alpha\in \Q[e_1,\ldots,e_n]$.

Let $x,y$ be the coordinates on $\A^2=\P^2\setminus \ell_\infty$.
The fixed point set $M(n)^{\widetilde T}$ is a set of 
$(\I_{Z_1},\Phi_1)\oplus (\I_{Z_2},\Phi_2)$, where the $\I_{Z_i}$ are
 ideal sheaves
of zero dimensional schemes $Z_1$, $Z_2$  with support in the origin of $\A^2$ 
with $\len(Z_1)+\len(Z_2)=n$ and $\Phi_\alpha$ ($\alpha=1,2$) are isomorphisms of $\I_{Z_\alpha}|_{\ell_\infty}$ with the $\alpha$-th factor of 
$\oo_{\ell_\infty}^{\oplus 2}$.
Write $I_\alpha$  for the ideal of $Z_\alpha$ in $\C[x,y]$. Then the above is a fixed point
if and only if $I_1$ and $I_2$ are generated by monomials in $x,y$.

A Young diagram is a set 
$$Y:=\big\{ (i,j)\in \Z_{> 0}\times \Z_{> 0}\bigm| j\le \lambda_i\big\},$$
where 
$(\lambda_i)_{i\in \Z_{> 0}}$ is a partition, i.e. 
$\lambda_i\in \Z_{\ge 0}$, $\lambda_i\ge \lambda_{i+1}$ for all $i$ and
only finitely many $\lambda_i$ are nonzero. Thus 
$\lambda_i$ is the length of the $i$-th column of $Y$.
Let $|Y|$ be number of elements of $Y$, so that 
$(\lambda_i)$ 
is a partition of $|Y|$.
We denote by 
$(\lambda_j')_j$ be the transpose of $\lambda$, thus 
$\lambda'_j$ is the length of the $j$-th row of $Y$.
For elements $s=(i,j)\in \Z_{>0}\times \Z_{>0}$ we put
$$a_Y(i,j)=\lambda_i-j,\quad l_Y(i,j)=\lambda_j'-i, \quad 
a'(i,j)=j-1, \quad l'(i,j)=i-1.$$

Let $I_Z\subset \C[x,y]$ be the ideal of a finite subscheme of $\A^2$ supported
in the origin which is generated by monomials in $x,y$.
To $Z$ we associate the Young diagram
$$Y=Y_Z:=\big\{(i,j)\in \Z_{>0}\times\Z_{>0} \bigm|
 x^{i-1}y^{j-1}\not\in I_Z\big\}.$$
 with $|Y|=\len(Z)$.
To a fixed point  
$(\I_{Z_1}\oplus \I_{Z_2},\phi)$  of the $\widetilde T$-action
on $M(n)$ we associate $\vec Y=(Y_1,Y_2)$ with $Y_i=Y_{Z_i}$. 
This gives  a bijection of the fixed point set $M(n)^{\widetilde T}$
with the set of pairs of Young diagrams $\vec Y=(Y_1,Y_2)$ with 
$|\vec Y|:=|Y_1|+|Y_2|=n$.

\begin{Notation}
We denote $e$ the one-dimensional 
$\widetilde T$-module given by $(t_1,t_2,e)\mapsto e$.
and similar we write $t_i$  ($i=1,2$) for the  $1$-dimensional
$\widetilde T$ modules given by $(t_1,t_2,e)\mapsto t_i$.
We also write $e_1:=e^{-1}$, $e_2:=e$.
We write $a_1:=-a$, $a_2:=a$.
\end{Notation}
Following \cite{NY1},\cite{NY2} let, for  $\alpha,\beta\in \{0,1\}$,  $N_{\alpha,\beta}^{\vec{Y}}(t_1,t_2,e)$ be the $\widetilde  T$-equivariant character of $\Ext^1(\I_{Z_{\alpha}},\I_{Z_{\beta}}(-\ell_{\infty}))$ and $n_{\alpha,\beta}^{\vec{Y}}(\ve_1,\ve_2,a)$
the equivariant Euler class.
\begin{NB}\begin{align*}
N^{\vec Y}_{\alpha,\beta}(t_1,t_2,e)&:=e_\beta e_\alpha^{-1}
\Big(\sum_{s\in Y_\alpha} t_1^{-l_{Y_\beta}(s)}t_2^{a_{Y_\alpha}(s)+1}+
\sum_{t\in Y_\beta} t_1^{l_{Y_\alpha}(t)+1}t_2^{-a_{Y_\beta}(t)}\Big).\\
n^{\vec Y}_{\alpha,\beta}(\ve_1,\ve_2,a)&:=
\prod_{s\in Y_{\alpha}}
\big(-l_{Y_\beta}(s)\ve_1+(a_{Y_\alpha}(s)+1)\ve_2+a_\beta-a_\alpha\big)
\\ &\qquad \qquad \qquad \qquad \cdot \prod_{s\in Y_{\beta}}
\big((l_{Y_\alpha}(s)+1)\ve_1-a_{Y_\beta}(s)\ve_2+a_\beta-a_\alpha\big).
\end{align*}
\end{NB}
Now the instanton part of the Nekrasov partition function is defined as
$$
\Zin(\ve_1,\ve_2,a,\Lambda):=\sum_{n \ge 0}\Lambda^{4n} \Big(\int_{M(n)} 1\Big)
=\sum_{\vec Y} \frac{\Lambda^{4|\vec Y|}}{\prod_{\alpha,\beta=1}^2 
n^{\vec Y}_{\alpha,\beta}(\ve_1,\ve_2,a)}.
$$
More generally we will consider the following:
For variables $\vec{\tau}:=(\tau_\rho)_{\rho\ge 1}$ let
\begin{equation}
\label{FY}
E^{\vec Y}(\ve_1,\ve_2,a,\vec{\tau}):=\exp\Bigg(\sum_{\rho=1}^\infty
\sum_{\alpha=1}^2\tau_\rho\Bigg[\frac{e^{a_\alpha}}{\ve_1\ve_2}
\Big(1-(1-e^{-\ve_1})(1-e^{-\ve_2})
\sum_{s\in Y_\alpha} e^{-l'(s)\ve_1-a'(s)\ve_2}\Big)\Bigg]_{\rho-1}\Bigg).
\end{equation}
(The sign in  \cite[(4.1)]{NY2} is not correct. See the first claim in the 
proof of \lemref{instan}.)
\begin{NB} Signs corrected.  6.9 LG
\end{NB}
Here $[\cdot]_{\rho-1}$ means the part of degree $\rho-1$, where $a,\ve_1,\ve_2$ have degree $1$.
Then the instanton part of the partition function is defined as
\begin{equation}
\label{Zinst}
\Zin(\ve_1,\ve_2,a;\Lambda,\vec{\tau})
:=\sum_{\vec Y}^\infty \frac{\Lambda^{4|\vec Y|}E^{\vec Y}(\ve_1,\ve_2,a,\vec{\tau})}
{\prod_{\alpha,\beta=1}^2 
n^{\vec Y}_{\alpha,\beta}(\ve_1,\ve_2,a)}\in \Q(\ve_1,\ve_2,a)[[\Lambda]].
\end{equation}
In particular $\Zin(\ve_1,\ve_2,a;\Lambda,\vec{0})=\Zin(\ve_1,\ve_2,a;\Lambda)$.
As a power series in $\Lambda$, \linebreak[2]$\Zin(\ve_1,\ve_2,\linebreak[1]a;\Lambda,\vec{\tau})$
starts with $1$. Thus 
$$\Fin(\ve_1,\ve_2,a;\Lambda,\vec{\tau}):=\log\Zin(\ve_1,\ve_2,a;\Lambda,\vec{\tau})
\in \Q(\ve_1,\ve_2,a)[[\Lambda]]$$
is well-defined and we put $\Fin(\ve_1,\ve_2,a;\Lambda):=\Fin(\ve_1,\ve_2,a;\Lambda,\vec{0}).$
Finally we define the perturbation part. 
We define $c_n$ ($n\in\Z_{\ge0}$) by
\begin{equation}\label{cm}
\frac{1}{(e^{\ve_1 t}-1)(e^{\ve_2t}-1)}=
\sum_{n \geq 0} \frac{c_n}{n!}t^{n-2},
\end{equation}
and define
\begin{equation}\label{perpart}
\gamma_{\ve_1,\ve_2}(x;\Lambda):=
\begin{aligned}[t]
\frac{1}{\ve_1\ve_2}\Big\{-\frac{1}{2}x^2\log\Big(\frac{x}{\Lambda}\Big)+\frac{3}{4}x^2\Big\}+\frac{\ve_1+\ve_2}{2\ve_1\ve_2}\Big\{-x\log\Big(\frac{x}{\Lambda}\Big)+x\Big\}\qquad\\
-\frac{\ve_1^2+\ve_2^2+3\ve_1\ve_2}{12\ve_1\ve_2}\log\Big(\frac{x}{\Lambda}\Big)+
\sum_{n=3}^\infty \frac{c_nx^{2-n}}{n(n-1)(n-2)}.
\end{aligned}
\end{equation}
We put 
$$\Fper(\ve_1,\ve_2,a;\Lambda):=-\gamma_{\ve_1,\ve_2}(2a;\Lambda)-\gamma_{\ve_1,\ve_2}(-2a;\Lambda).$$
Then $\Fper(\ve_1,\ve_2,a;\Lambda)$ is a Laurent series in $\ve_1,\ve_2$,
whose coefficients are multiple-valued meromorphic functions in $a,\Lambda$.
See \cite[Appendix E]{NY2} for the details.
Finally we define 
\begin{gather*}
F(\ve_1,\ve_2,a;\Lambda,\vec{\tau}):=\Fper(\ve_1,\ve_2,a;\Lambda)+\Fin(\ve_1,\ve_2,a;\Lambda,\vec{\tau}), 
\\
F(\ve_1,\ve_2,a;\Lambda):=F(\ve_1,\ve_2,a;\Lambda,\vec{0}).
\end{gather*}
Formally one defines $Z(\ve_1,\ve_2,a;\Lambda,\vec{\tau})):=\exp(\Fper(\ve_1,\ve_2,a;\Lambda))\Zin(\ve_1,\ve_2,a;\Lambda,\vec{\tau})$.

\section{Computation of the wallcrossing in terms of Hilbert schemes}

Let $X$ be a simply connected 
smooth projective  surface with $p_g=0$. 
In this section we will compute the wallcrossing of the Donaldson invariants of $X$
in terms of intersection numbers of Hilbert schemes of points on $X$. 
Our result will be more generally about a refinement of the Donaldson invariants, also
involving higher order $\mu$-classes.
In the next two sections 
we will specialize to the case that $X$ is a
smooth toric surface and relate this result to the Nekrasov partition
function.


\begin{Notation}
Let $t$ be a variable. If $Y$ is a variety and 
$b\in H^*(Y)[t]$, we denote
by
$[b]_d$ its part of degree $d$, where elements in 
$H^{2n}(Y)$ have degree $n$ 
and $t$ has degree $1$.

If $R$ is a ring, $t$ a variable and $b\in R((t))$, we will denote for 
$i\in \Z$ by
$[b]_{t^i}$ the coefficient of $t^i$ of $b$.

If $E$ is a torsion free sheaf of rank $r$ on $Y$, then we put 
$\cho(E):=\ch(E)e^{-\frac{c_1(E)}{r}}$.
We write $\cho_i(E):=[\cho(E)]_i$.
We can view this as Chern character of $E$ normalized  by  a twist  with a 
rational line bundle, so that its
first Chern class is zero.
Note that in case $r=2$, we have $-\cho_2(E)=c_2(E)-c_1^2(E)/4$.

If $E$ is a vector bundle of rank $r$ on $Y$
we write $c^{t}(E):=\sum_i c_i(E) t^{r-i}=t^{r}c_{1/t}(E)$, with
$c_t(E):=\sum_i c_i(E) t^i$.
\end{Notation}

Now we define a generalization of the $\mu$-map and of the 
Donaldson invariants.
\begin{Definition}
Fix an ample divisor $H$ on $X$ and fix $c_1$ and $d$.
Assume that there is a universal sheaf $\E$ over $X\times M_H^X(c_1,d)$,
and that $M_H^X(c_1,d)$ is of expected dimension $d$ or empty for all $d\ge 0$.
For a class $a\in H_i(X)$,  and an integer $\rho\ge 1$, we put
$\mu_\rho(a):=(-1)^\rho\cho_{\rho+1}(\E)/a\in H^{2\rho+2-i}(M_H^X(c_1,d))$.
Note that the universal sheaf is well-defined up to a twist by the pullback of 
a line bundle from $M_H^X(c_1,d)$, thus 
$\mu_\rho$ is independent of the choice of 
the universal sheaf.

Let $b_1,\ldots,b_s$ be a homogeneous basis of $H_*(X)$.
For all $\rho\ge 1$ let $\tau_1^\rho, \ldots \tau_s^\rho$ indeterminates,
and put $\alpha_\rho:=\sum_{k=1}^s a_k^\rho b_k \tau_k^\rho$, with 
$a_k^\rho\in\Q$. This means that $(\tau_k^\rho)_{k=1\ldots s}^{\rho\ge 1}$ 
is a coordinate system on the ``large phase space" $\bigoplus_{\rho\ge 1}H_*(X)[\rho]$. 
We define
\begin{equation}
\label{donna}\Phi_{c_1}^H\Bigg(\exp\Big(\sum_{\rho\ge 1} \alpha_\rho
\Big)
\Bigg)=\sum_{d\ge 0}\Lambda^d\int_{M_H^X(c_1,d)}\exp\Big(\sum_{\rho\ge 1}
\mu_\rho(\alpha_\rho)\Big).\end{equation}
This is an element of $\Q[[\Lambda,(\tau_k^\rho)]]$.
As by definition $\mu_1=\mu$, 
our previous definition of
$\Phi_{c_1}^H(\exp(\alpha z +p x))$ is obtained  by specializing $\alpha _\rho:=0$ for all $\rho>1$.
\end{Definition}

We believe that we can define the invariants without the assumption
that the moduli spaces are of expected dimensions as in the case of
ordinary Donaldson invariants. This can be done once we generalize the
blowup formula. This is a little delicate as higher Chern classes do
not descend to Uhlenbeck compactifications.

\subsection{The wallcrossing term}\label{wallmod}
Let $\xi\in  H^2(X,\Z)\setminus\{0\}$ \begin{NB} corrected 6.9 LG\end{NB}
be
a class of type $c_1$. 
We say that $\xi$ is {\it good\/} and $W^\xi$ is a {\it good wall\/} if
\begin{enumerate}
\item there is an ample divisor in $W^\xi$
\item  $D+K_X$ is not effective for any 
divisor $D$ with $W^{c_1(D)}=W^\xi$.
\end{enumerate}
\begin{NB}
The equation $W^{c_1(D)}=W^\xi$ is equivalent to $a c_1(D) = b\xi$ for
some $a$, $b\in\Z\setminus\{ 0\}$. We cannot replace this to
$c_1(D) = \xi$.
The point is that when one crosses the wall defined by $\xi$, at
the same time something happens for all $D$ with $a\xi=b c_1(D)$.
Apr. 30, H.N
\end{NB}
A sufficient condition for $\xi$ to be good  is that $W^\xi$ contains an 
ample divisor $H$ with $H\cdot K_X<0$. Let $\xi$ be a good class of type $c_1$.

Let $X^{[n]}$ be the Hilbert scheme of subschemes of length $n$ on $X$.
Let $Z_n(X)\subset X\times X^{[n]}$ be the universal subscheme. 
We write $X_2:=X\sqcup X$ and $X_2^{[l]}:=\coprod_{n+m=l} X^{[n]}\times X^{[m]}$.
Fix $l\in \Z_{\ge 0}$.
Let  $\I_1$ (resp.~$\I_2$) be the sheaf on $X\times X_2^{[l]}$ whose restriction
 to $X\times X^{[n]}\times X^{[m]}$ is $p_{1,2}^*(\I_{Z_n(X)})$ (resp.~
$p_{1,3}^*(\I_{Z_m(X)})$), where
$p_{1,2}\colon X\times X^{[n]}\times X^{[m]}\to X\times X^{[n]}$
(resp.\ $p_{1,3}\colon X\times X^{[n]}\times X^{[m]}\to X\times X^{[m]}$).
Let $p:X\times X_2^{[l]}\to X_2^{[l]}$ be the projection. 
On  $X_2^{[l]}$ we define 
\begin{align*}
\AA_{\xi,-}&:=\Ext^1_{p}(\I_{2},\I_{1}(\xi)),\qquad
\AA_{\xi,+}:=\Ext^1_{p}(\I_{1},\I_{2}(-\xi)).
\end{align*}
As $\xi$ is good, $\AA_{\xi,-}$, $\AA_{\xi,+}$ are locally free on 
$X_2^{[l]}$. If $\xi$ is understood, we also just write $\AA_-$ and $\AA_+$ instead
of  $\AA_{\xi,-}$, $\AA_{\xi,+}$. 
Let $\P_-:=\P(\AA_-^\vee)$ and $\P_+:=\P(\AA_+^\vee)$
 (we use the Grothendieck notation, i.e. this is the bundle of $1$-dimensional quotients). Let $\pi_{\pm}:\P_{\pm}
\to X_2^{[l]}$ be the projection. Then $\P_{\pm}=\coprod_{n+m=l} \P_{\pm}^{n,m}$
with  $\P_{\pm}^{n,m}=\pi_{\pm}^{-1}(X^{[n]}\times X^{[m]})$.

Now we define the wallcrossing term. We use the notations of 
the last section.
For a coherent sheaf $E$ of rank $r$ on a variety $Y$,
we view $\frac{1}{c^t(E)}$ as an element of $H^*(Y)[t^{-1}]$
via the formula
\begin{equation}\label{cct}\frac{1}{c^t(E)}=\frac{1}{t^r}\frac{1}{\sum_{i=0}^r c_i(E)\frac{1}{t^i}}=t^{-r}
\sum_i s_i(E) t^{-i},\end{equation}
where $r$ is the rank of $E$ and the $s_i(E)$ are the Segre classes of $E$. If $Y$ carries a $\Gamma$-action  and $E$ is equivariant, then  $\frac{1}{c^t(E)}\in H^*_\Gamma(Y)[[t^{-1}]]$.

\begin{Definition}
Let $\xi\in H^2(X,\Z)$ be a good class of type $c_1$.
For all $\rho\ge 1$ let $\alpha_\rho$ be as in (\ref{donna}).
The {\em wallcrossing terms} are
\begin{equation}\label{wallct}
\begin{split}&\delta^X_{\xi,t}\Big(\exp\Big(\sum_{\rho\ge 1} \alpha_\rho)
\Big)\Big)
\\
&\qquad\qquad:=
\begin{aligned}[t]
\sum_{l\ge0}
\Lambda^{4l-\xi^2-3}
\int_{X_2^{[l]}}\frac{\exp\Big(\sum_{\rho\ge 1}(-1)^{\rho}\big[\ch(
\I_{1})e^{\frac{\xi-t}{2}}+ \ch(\I_{2})e^{\frac{t-\xi}{2}}]_{\rho+1}
/\alpha_{\rho} \Big)}
{c^t(\AA_{\xi,+})c^{-t}(\AA_{\xi,-})},
\end{aligned}
\\
& \delta^X_\xi\Big(\exp\Big(\sum_{\rho\ge 1} \alpha_\rho
\Big)\Big):=\Big[\delta^X_{\xi,t}\Big(\exp\Big(\sum_{\rho\ge 1} \alpha_\rho
\Big)\Big)\Big]_{t^{-1}}.
\end{split}
\end{equation}
$\delta^X_{\xi,t}\big(\exp(\alpha z+p x)\big)$ and $\delta^X_{\xi}\big(\exp(\alpha z+p x)\big)$ are defined by replacing $\alpha_1$ by $\alpha z+p x$ and $\alpha_\rho$ by $0$ for $\rho\ge 2$ in $\delta^X_{\xi,t}\big(\exp\big(\sum_{\rho} \alpha_\rho)
\big)\big)$ and $\delta^X_{\xi}\big(\exp\big(\sum_{\rho} \alpha_\rho)
\big)\big)$.
By (\ref{cct}) we see that
$\delta^X_{\xi,t}\Big(\exp\Big(\sum_{\rho\ge 1} \alpha_\rho
\Big)\Big)\in\Lambda^{-\xi^2-3}\Q[t,t^{-1}][[\Lambda,(\tau_k^\rho)]]$ and
$\delta^X_{\xi}\Big(\exp\Big(\sum_{\rho\ge 1} \alpha_\rho
\Big)\Big)\in\Q[[\Lambda,(\tau_k^\rho)]]$ (see Remark \ref{weight}(1)).
\end{Definition}
\begin{Remark}\label{weight}
(1) Fix $l\ge 0$. Write $d:=4l-\xi^2-3$ and 
let $E(t):=\frac{1}{c^t(\AA_{\xi,+})c^{-t}(\AA_{\xi,-})}$ on $X^{[l]}_2$. 
Note that
$\rank(\AA_{\xi,+}\oplus\AA_{\xi,-})=d+1-2l$ (this follows from \cite[Lemma 4.3]
{EG1}). If $d<0$, then $d+1-2l\le 0$, thus $\AA_{\xi,+}=0=\AA_{\xi,-}$
and $E(t)=1$, and thus the coefficient of $\Lambda^d$ in 
$\delta_{\xi,t}^X\big(\exp\big(\sum \alpha_\rho
\big)\big)$ is a polynomial in $t$. 
Let again $d$ be arbitrary. We can write 
$E(t)=\sum_{i=0}^{2l} b_i t^{-(i+d+1-2l)}$, with 
$b_{i}\in H^{2i}(X^{[l]}_2)$. Thus if we give elements of $H^{2i}(X^{[l]}_2)$
the degree $i$ and $t$ the degree $1$, then $E(t)$ is homogeneous of degree $2l-d-1$.



(2) Note that the factor  $\Lambda^d$ both in the definition of $\Phi_{c_1}^H$
and of $\delta_\xi^X$ is redundant. The coefficient of 
$\Lambda^d$ in $\Phi^H_{c_1}(\exp(\sum \alpha_\rho))$ and 
 $\delta_\xi^X(\exp(\sum \alpha_\rho))$ is  
a polynomial of weight $d$ in the $\tau_k^{\rho}$. Here the weight of $\tau^\rho_k$ is
$\rho+1-i$ if $b_k\in H_{2i}(X)$. For $\Phi_{c_1}^X$, this is clear because
$d$ is the complex dimension of $M_X^H(c_1,d)$. 
For $\delta_\xi^X$, this follows easily 
from the last sentence of (1) and the fact that $X^{[l]}_2$ has dimension $2l$.
\end{Remark}
The aim of this section is to prove that the wallcrossing for the 
Donaldson invariants can be expressed as a sum over $\delta^X_\xi$.

\begin{Proposition}\label{wallcr}
Let $H_-$, $H_+$ be ample divisors on $X$, which do not lie on a wall 
of type $(c_1,d)$ for any $d\ge 0$. Let $B_+$ be the set of all 
classes $\xi$ of type $c_1$
with $\<\xi\cdot H_+\>>0 >\<\xi\cdot H_-\>$. Assume that all classes in $B_+$ 
are good. 
Then  
\begin{align*}
\Phi_{c_1}^{H_+}\Big(\exp\Big(\sum_{\rho} \alpha_\rho\Big)\Big)-\Phi_{c_1}^{H_-}
\Big(\exp\Big(\sum_{\rho} \alpha_\rho\Big)
\Big)
&=\sum_{\xi\in B_+}\delta^X_{\xi}
\Big(\exp\Big(\sum_{\rho} \alpha_\rho\Big)\Big).
\end{align*}
\end{Proposition}

\begin{Remark}
{}From our final expression in \corref{cor:main3},
$\delta^X_\xi(\exp(\alpha z+px))$ is compatible with Fintushel-Stern's
blowup formula \cite{FS}. (See \cite[\S4.2]{GZ} and \cite[\S6]{NY2}.)
Therefore it is enough to prove the proposition after we blowup $X$ at
sufficiently many times, as we did for the definition of
$\Phi^H_{c_1}$. In particular, we may assume $M^X_{H_\pm}(c_1,d)$ is
of expected dimension without loss of generality.
However the blowup does not make walls good in general, so we one needs 
different methods to prove the proposition for
general wall. Let $p:X\times X_2^{[l]}\to X_2^{[l]}$ and 
$q:X\times X_2^{[l]}\to X$, be the projections.  In \cite[Thm 1.12]{Moc} the proposition is proved for general walls
with $\AA_{\xi,+}$, $\AA_{\xi,-}$ replaced by  $-p_!(\I_2^\vee\otimes \I_1\otimes q^!\xi)$, $-p_!(\I_1^\vee\otimes \I_2\otimes q^!\xi^\vee)$). The proof uses virtual fundamental classes and virtual localization.
\begin{NB} changed 6.9. LG\end{NB}
\end{Remark}

In the rest of this section we will show Prop.~\ref{wallcr}. 
Let $d\ge 0$ be arbitrary.
It is enough to show that 
the coefficients of $\Lambda^d$ on both sides are equal. 
It is known that $M^X_H(c_1,d)$ and $\Phi^H_{c_1}$ is constant 
as long as $H$ stays in the same  chamber of type $(c_1,d)$ 
and only changes when
$H$ crosses a wall of type $(c_1,d)$. 
Following \cite{EG1} and \cite{FQ} we get the following description of the change
of moduli spaces. 
Let $B_d$ be the set of all $\xi\in B_+$ which define a wall of type $(c_1,d)$.
For the moment assume for simplicity that $B_d$ consists of a single element
$\xi$. 
Let $l:=(d+3+\xi^2)/4\in \Z_\ge 0$.
Write $M_{0,l}:=M^X_{H_-}(c_1,d)$.  Then successively for all $n=0,\ldots,l$ write 
$m:=l-n$. Then one has the following: 
$M_{n,m}$ contains a closed subscheme $E_-^{n,m}$ isomorphic to $\P_{-}^{n,m}$
and $M_{n,m}$ is nonsingular in a neighbourhood of $E_-^{n,m}$.  
Let $\widehat M_{n,m}$ be the blow up of $M_{n,m}$  along $E_-^{n,m}$.
The exceptional divisor is isomorphic to the fibre product  
$D^{n,m}:=\P_-^{n,m}\times_{X^{[n]}\times X^{[m]}} 
\P_{+}^{n,m}$. We can blow down $\widehat M_{n,m}$  in $D^{n,m}$ 
in the other fibre direction to obtain a new variety  $M_{n+1,m-1}$. 
The image of $D^{n,m}$ is a closed subset $E_+^{n,m}$ isomorphic to 
$\P_+^{n,m}$ and $M_{n+1,m-1}$ is smooth in a neighbourhood of $E_+^{n,m}$.

The transformation from $M_{n,m}$ to $M_{n+1,m-1}$  does not have to be  birational. 
It is possible that $E_+^{n,m}=\emptyset$, i.e. $\AA_+=0$. 
As we know that 
$\rank(\AA_-)+\rank(\AA_+)+2l=d+1$, this happens if and only if $E_-^{n,m}$ has dimension $d$
and thus by the smoothness of  $M_{n,m}$ near $E_-^{n,m}$, we get that $E_-^{n,m}$
is a connected component of  $M_{n,m}$. Then blowing up along $E_-^{n,m}$ just means
deleting $E_-^{n,m}$. Thus in this case $M_{n+1,m-1}=M_{n,m}\setminus E_-^{n,m}$.
Similarly  we have $E_-^{n,m}=\emptyset$, i.e. $\AA_-=0$, if and only if 
$E_+^{n,m}$ is a 
connected component of $M_{n+1,m-1}$ and  $M_{n+1,m-1}=M_{n,m}\sqcup E_+^{n,m}$.
Below, if the transformation from $M_{n,m}$ to $M_{n+1,m-1}$ is birational, we say 
we are in case (1), otherwise in case (2).

Finally we have $M_{l+1,-1}=M_{H_+}^X(c_1,d)$. 
If $B_d$ consists of more than one element, one obtains $M_{H_+}(c_1,d)$
from $M_{H_-}(c_1,d)$  by iterating  this procedure in a suitable order over  all
$\xi\in B_+$.

Fix $\xi$ in $B_d$. 
Fix $n,m\in \Z_{\ge 0}$ with $n+m=l:=(d+3+\xi^2)/4$. 
We write $M_-:=M_{n,m}$, $M_+:=M_{n+1,m-1}$. Let $\overline \E_\pm$ be 
universal sheaves on $X\times M_\pm$ respectively.
Let $E_-:=E_-^{n,m}$, $E_+=E_+^{n,m}$.
Let $\widetilde M$ be the blowup of $M_-$ along $E_-$, and denote
by $D$ the exceptional divisor (which is also the exceptional divisor 
of the blowup of $M_+$ along $E_+$). Write $D':=X\times D$ and let $j:D\to \widetilde M$,
$j':X\times D\to X\times \widetilde M$ be the embeddings.
Let $\E_-$, $\E_+$ be the pullbacks of $\overline \E_-$, $\overline 
\E_+$ to $X\times \widetilde M$.

\begin{Notation} 
%
Let $H\in \Q[(x_n)_{n>0}]$ be a polynomial. Let $a:=(a_n)_{n>0}$ with
$a_n\in H_*(X)$.
For any variety $Y$ and any class $A\in H^*(X\times Y)[[t]]$ we put
$H(A/a):=H(([A]_n/a_n)_{n> 0})\in H^*(Y)[[t]].$
On $X\times X^{[n]}\times X^{[m]}$, denote 
$C(t):=\ch(\I_1)e^{\frac{\xi-t}{2}}+\ch(\I_2)e^{\frac{-\xi+t}{2}}$, 
and $C_i(t):=[C(t)]_i$.

We denote by $\tau_-$ (resp. $\tau_+$) the universal quotient line bundle on 
$\P_-=\P(\AA_-^\vee)$ (resp. $\P_+=\P(\AA_+^\vee)$).
For a sheaf $\F$ and a divisor $B$, we write $\F(B)$ instead of $\F\otimes \oo(B)$.

For a class $a\in H^*(X)$ we also denote by $a$ its pullback to $X\times Y$
for a variety $Y$.
We write $\I_1,\I_2$ also for the pullback of $\I_1$, $\I_2$ to $D'$ and we write
$\tau_+,\tau_-$ also for their pullbacks to $D$ and $D'$.

%
\end{Notation}

We will show
\begin{equation}
\label{restr}
\int_{M_+}H(\cho(\overline \E_+)/a)-\int_{M_-}H(\cho(\overline\E_-)/a)=
\int_{X^{[n]}\times X^{[m]}}
\Bigg[\frac{H(C(t)/a)}{c^t(\AA_{\xi,+})c^{-t}(\AA_{\xi,-})}\Bigg]_{t^{-1}}.
\end{equation}
Formula (\ref{restr}) implies Proposition \ref{wallcr} by summing over all 
$\xi\in B_+$, all $d\ge 0$ and over all $n,m$ with $n+m=(d+\xi^2+3)/4$.

For the next three Lemmas assume that we are in case (1). 
Then by the projection formula 
$\int_{M_{\pm}}H(\cho(\overline \E_{\pm}/a)))=
\int_{\widetilde M} H(\cho(\E_\pm/a)),$ thus it is enough to prove (\ref{restr})
with the left-hand side replaced by 
$\int_{\widetilde M}\big(H(\cho(\E_+)/a)-H(\cho(\E_-)/a)\big)$.




\begin{Lemma}
$C(-\tau_-)=\cho((j')^*(\E_-))$, 
$C(\tau_+)=\cho((j')^*(\E_+))$ and
$$\cho(\E_+)-\cho(\E_-)=-j'_*\Big(\frac{C(t)-C(-s)}{s+t}|_{{s=\tau_-}
\atop{t=\tau_+}}
\Big).$$
\end{Lemma}
\begin{proof}
Write $\F_1:=\I_1(\frac{c_1+\xi}{2})$, $\F_2:=\I_2(\frac{c_1-\xi}{2})$.
By   \cite[section 5]{EG1} we have the following facts:
\begin{enumerate}
\item There exist a line bundle $\lambda$ on $D$ and an exact sequence
$0\to \F_1(\lambda)\to (j')^*(\E_-)\to \F_2(-\tau_-+\lambda)
\to 0$,
\item  $\E_+$ can be defined by the
exact sequence
$0\to \E_+\to \E_-\to j'_*(\F_2(-\tau_-+\lambda))\to 0.$
\item We have the exact sequence
$0\to \F_2(\tau_++\lambda)\to (j')^*(\E_+)\to 
\F_1(\lambda)
\to 0.$
\end{enumerate}
In particular $\cho((j')^*\E_-)=C(-\tau_-)$, 
$\cho((j')^*\E_+)=C(\tau_+)$.

Write $c_+:=c_1(\E_+),c_-:=c_1(\E_-)$.
As $c_1(j'_*(\F_2(-\tau_-+\lambda)))=D'$,
we see that $c_+=c_--D'$. We also have
$(j')^*(c_+)={c_1+\tau_+}+2\lambda$.
Thus we get
\begin{multline*}
\cho(\E_+)=\big(\ch(\E_-)-\ch(j'_*\F_2(-\tau_-+\lambda))\big)e^{-c_+/2}\\=
\cho(\E_-)e^{D'/2}-\ch(j'_*\F_2(-\tau_-+\lambda))e^{-c_+/2}.
\end{multline*}
Thus $\cho(\E_+)-\cho(\E_-)=(e^{D'/2}-1)\cho(\E_-)
-\ch(j'_*\F_2(-\tau_-+\lambda))e^{-c_+/2}$.
As $(j')^*D'=-\tau_+-\tau_-$
by \cite[Cor.~4.7]{EG1},
we get by 
the Grothendieck-Riemann-Roch Theorem and the 
projection formula
\begin{align*}
\ch\big(j'_*\F_2(-\tau_-+\lambda)\big)e^{-c_+/2}
&=j'_*\Big(\frac{1-e^t}{-t}|_{t=\tau_++\tau_-}\ch(\F_2(-\tau_-+\lambda))\Big)
e^{-c_+/2}\\
&=j'_*\Big(\frac{1-e^t}{-t}|_{t=\tau_++\tau_-}\ch(\I_2)
e^{\frac{-\xi-2\tau_--\tau_+}{2}}\Big)\\
&=-j'_*\Big(\Big(\frac{1}{s+t}\left(\ch(\I_2)e^{\frac{-\xi-2s-t}{2}}-\ch(\I_2)
e^{\frac{-\xi+t}{2}}\right)\Big)|_{
{s=\tau_-}\atop{t=\tau_+}}
\Big)
\end{align*}
On the other hand, as $e^{D'/2}-1$ is divisible by $D'$, we get
\begin{align*}
(e^{D'/2}&-1)\cho(\E_-)
=j'_*\Big(\frac{1-e^{-t/2}}{t}|_{t=\tau_++\tau_-}\cho((j')^*\E_-)\Big)\\
&=
j'_*\Big(\frac{1}{s+t}\big(\ch(\I_1)e^{\frac{\xi+s}2}+\ch(\I_2)e^{\frac{-\xi-s}2}
-\ch(\I_1)e^{\frac{\xi-t}2}
-\ch(\I_2)e^{\frac{-\xi-2s-t}2}\big)|_{{s=\tau_-}\atop{t=\tau_+}}\Big),
\end{align*}
and the result follows.
\end{proof}


\begin{Lemma}\label{restr1}
$$H(\cho(\E_+)/a)-H(\cho(\E_-)/a)=
-j_*\Big(\frac{H(C(t)/a)-H(C(-s)/a)}
{s+t}|_{{s=\tau_-}\atop{t=\tau_+}}\Big),$$
In particular 
$$
\int_{\widetilde M}\big(H(\cho(\E_+)/a)-H(\cho(\E_-)/a)\big)=
-\int_D\frac{H(C(t)/a)-H(C(-s)/a)}{s+t}|_{{s=\tau_-}\atop{t=\tau_+}}.$$
\end{Lemma}

\begin{proof}
We can assume that $H$ is homogeneous of degree $k$.
We make induction over $k$, the case $k=0$ being trivial. 
In case $k=1$, we have by the previous Lemma
\begin{multline*}
\cho_{i}(\E_+)/a_i-\cho_{i}(\E_-)/a_i\\ =
-j'_{*} \Big(\frac{C_{i}(t)-C_{i}(-s)}{s+t}|_{{s=\tau_-}\atop{t=\tau_+}}\Big)/a_i
=-j_*\Big(\frac{(C_{i}(t)-C_{i}(-s))/a_i}{s+t}|_{{s=\tau_-}\atop{t=\tau_+}}
\Big).
\end{multline*}
Now let $k$ be general. As the claim is linear in $H$, we can assume that 
$H=x_{i}H'$, with $\deg(H')=k-1$. 
Thus we get by induction
\begin{align*}
H(\cho(\E_+)/a)&-H(\cho(\E_-)/a)=
\big(\cho_{i}(\E_+)/a_{i}-\cho_{i}(\E_-)/a_{i}\big)
H'(\cho(\E_+)/a)\\&\qquad\qquad\qquad\qquad\qquad\qquad+
\cho_{i}(\E_-)/a_{i}\cdot\big(H'(\cho(\E_+)/a)-H'(\cho(\E_-)/a)\big)\\
&=-j_*\Big(\Big(\frac{(C_{i}(t)-C_{i}(-s))/a_{i}}{s+t}
H'\big(\cho((j')^*\E_+)/a)
\\&\qquad\qquad\qquad\qquad+\big(\cho_{i}((j')^*\E_-)/a_{i}\big)
\frac{H'(C(t)/a)-H'(C(-s)/a)}{s+t}\Big)|_{{s=\tau_-}\atop{t=\tau_+}}\Big)
\\
&=-j_*\Big(\frac{1}{s+t}\big((C_{i}(t)-C_{i}(-s))/a_{i}\cdot
H'(C(t)/a)\\
&\qquad \qquad\qquad \qquad \qquad+C_{i}(-s)/a_{i}\cdot(H'(C(t)/a)
-H'(C(-s)/a))\big)|_{{s=\tau_-}\atop{t=\tau_+}}\Big)\\
&=-j_*\Big(\frac{H((C(t)/a)-H(C(-s))/a)}{s+t}|_{{s=\tau_-}\atop{t=\tau_+}}\Big).
\end{align*}
This shows the first statement, the second follows immediately 
by the projection formula.
\end{proof}
Recall that  $D=\P(\AA_-^\vee)\times_{X^{[n]}\times X^{[m]}}\P(\AA_+^\vee)$.
Let $\pi:D\to X^{[n]}\times X^{[m]}$ and
$p_{\pm}:\P(\AA_\pm^\vee)\to X^{[n]}\times X^{[m]}$  be the projections.
We have reduced the computation of 
$\int_{\widetilde M}\big(H(\cho(\E_+)/a)-H(\cho(\E_+)/a)\big)$ 
to an integral over $D$, which we now 
push down to $X^{[n]}\times X^{[m]}$.

\begin{Lemma}\label{pushd}
$$\pi_*\Big(-\frac{H(C(t)/a)-H(C(-s)/a)}{s+t}|_{{s=\tau_-}\atop{t=\tau_+}}\Big)
=\Big[
\frac{H(C(t)/a)}{c^t(\AA_+)c^{-t}(\AA_-)}\Big]_{t^{-1}}.$$
\end{Lemma}
\begin{proof}
For a vector bundle $E$ of rank $e$ on a variety $Y$,
let $\tau$ be the tautological quotient line bundle on
$p: \P(E^{\vee})\to Y$. Then
$$
\sum_{n\ge 0} p_*(\tau^n) t^{-n-1}=
t^{-e} \sum_n p_*(\tau^{n+e-1}) t^{-n}
=t^{-e}\frac{1}{\sum_n c_n(E) t^{-n}}
=\frac{1}{c^t(E)},
$$
and similarly 
$\sum_n p_*((-\tau)^n) t^{-n-1}=-\frac{1}{c^{-t}(E)}$.
Thus we get 
\begin{align*}
\pi_*\Big(\frac{\tau_+^k-(-\tau_-)^k}{\tau_++\tau_-}\Big)&
=\sum_{i+j=k-1}\pi_{+*}(\tau_+^i)\pi_{-*}((-\tau_-)^j)\\
&=\Big[\Big(\sum_n \pi_{+*}(\tau_+^n) t^{-n-1}\Big)
\Big(\sum_n \pi_{-*}((-\tau_-)^n) t^{-n-1}\Big)\Big]_{t^{-k-1}}\\
&=-\Big[\frac{t^{k}}{c^t(\AA_+)c^{-t}(\AA_-)}\Big]_{t^{-1}}.
\end{align*}
We write $H(C(t)/a)=\sum_{k\ge 0}t^k Q_k$ with 
$Q_k\in H^*(X^{[n]}\times X^{[m]})$.
Then 
\begin{align*}
\pi_*\Big(-&\frac{H(C(t)/a)-H(C(-s)/a)}{t+s}
|_{{s=\tau_-}\atop{t=\tau_+}}\Big)=
-\sum_{k}\pi_*\Big(\frac{\tau_+^k-(-\tau_-)^k}{\tau_++\tau_-}\Big)
Q_k\\
&=\sum_k \Big[\frac{t^k Q_k}{c^t(\AA_+)c^{-t}(\AA_-)}\Big]_{t^{-1}}
=\Big[\frac{H(C(t)/a)}{c^t(\AA_+)c^{-t}(\AA_-)}\Big]_{t^{-1}}.
\end{align*}
\end{proof}
The projection formula and Lemmas \ref{restr1}, \ref{pushd} imply
formula (\ref{restr}). Thus we have  shown (\ref{restr}) in case (1). 

In case (2), we can assume by symmetry that $\P_+=\emptyset$, thus $\AA_{+}=0$ 
and $\AA_-$ has rank $d+1-2l$. 
Then we have 
\begin{multline*}
\int_{M_+}H(\cho(\overline \E_+)/a)-\int_{M_-}H(\cho(\overline \E_-)/a)
\\ =
-\int_{\P_-}H(\cho(\overline \E_-)/a)=
-\int_{X^{[n]}\times X^{[m]}}\pi_{-*}(H(\cho(\overline \E_-)/a)).
\end{multline*}
Denote by $j:\P_-\to M_-$ and $j':X\times \P_-\to X\times M_-$ the embeddings.
As before write $\F_1:=\I_1(\frac{c_1+\xi}{2})$, $\F_2:=\I_2(\frac{c_1-\xi}{2})$.
By \cite[Lemma 4.3]{EG1} and the universal property of $M_-$ there is line bundle
$\lambda$ on $\P_-$ and an exact sequence
$0\to \F_1(\lambda)\to (j')^*(\overline \E_-)\to \F_2(-\tau_-+\lambda)
\to 0$.
In particular, as before, $\cho((j')^*\overline \E_-)=C(-\tau_-)$.
The arguments of Lemma \ref{pushd} show that
$-\pi_{-*}(-\tau_-)^k=\Big[\frac{t^{k}}{c^{-t}(\AA_-)}\Big]_{t^{-1}}$, and in the
same way as in the proof of Lemma \ref{pushd} it follows that
$-\pi_{-*}(H(C(-\tau_-)/a))=\Big[\frac{H(C(t)/a)}{c^{-t}(\AA_-)}\Big]_{t^{-1}}$.
As $c^t(\AA_+)=1$, this shows (\ref{restr}) also  in  case (2) and thus finishes the 
proof of Proposition \ref{wallcr}.

\section{Comparison with the partition function}

For the next two sections let $X$ be a smooth projective  toric surface over $\C$,
in particular $X$ is simply connected and $p_g(X)=0$. 
$X$ carries an action of $\Gamma:=\C^*\times \C^*$ with finitely many 
fixed points, which we will denote by $p_1,\ldots,p_\chi$, where $\chi$ is the 
Euler number of $X$.  Let $w(x_i)$, $w(y_i)$ be the weights of the $\Gamma$-action on $T_{p_i}X$. Then there are local coordinates $x_i,y_i$ at $p_i$, so that $(t_1,t_2) x_i=e^{-w(x_i)} x_i$, $(t_1,t_2) y_i=e^{-w(y_i)} y_i$.
\begin{NB} Corrected 6.9 LG\end{NB}
 By definition $w(x_i)$ and $w(y_i)$
are linear forms in $\ve_1$ and $\ve_2$. For $\beta\in H^*_\Gamma(X)$ or $\beta\in H_*^\Gamma(X)$, we denote by
$\iota_{p_i}^*\beta$ its  pullback to the fixed point $p_i$.
More generally, if $\Gamma$ acts on a nonsingular variety $Y$ and $W\subset Y$ is 
invariant under the $\Gamma$-action, we denote by 
$\iota_W:H^*_\Gamma(Y)\to H^*_\Gamma(W)$ the pullback homomorphism.

Note that $T_X$ and the canonical bundle  are canonically equivariant.
Thus any polynomial in the Chern classes $c_i(X)$ and $K_X$ is canonically an element
of $H^*_\Gamma(X)$.  

\subsection{Equivariant Donaldson invariants and equivariant wallcrossing}
We start by defining an equivariant version of the Donaldson invariants and the wallcrossing terms.
For $t\in \Gamma$ denote by $F_t$ the automorphism $X\to X; x\mapsto t\cdot x$.
Then $\Gamma$ acts on $X^{[l]}_2$ by $t\cdot(\I_{Y_1},\I_{Y_2})=((F_{t}^{-1})^*\I_{Y_1},(F_{t}^{-1})^*\I_{Y_2})$ and on 
$X\times X^{[l]}$ by $t\cdot(x,\I_{Y_1},\I_{Y_2})=(F_t(x),(F_{t}^{-1})^*\I_{Y_1},
(F_{t}^{-1})^*\I_{Y_2})$ and the sheaves $\I_{1}$,
$\I_{2}$ are $\Gamma$-equivariant. Similarly $\Gamma$ acts on 
$X\times M^H_X(c_1,d)$ by $t\cdot(x,E)=(F_t(x),(F_{t}^{-1})^*E)$. Let $\E$ be a universal sheaf over $X\times M^H_X(c_1,d)$,
then one can show that $\E$ has a lifting to a $\Gamma$-equivariant sheaf, unique up
to twist by a character. Thus an equivariant universal sheaf is unique up to twist by an 
equivariant line bundle.

\begin{Definition}\label{eqwall}
We define the equivariant Donaldson invariants
$\widetilde\Phi_{c_1}^{H}(\exp(\alpha z +p x))$ by the right-hand side of
(\ref{donn}), where now  $\alpha\in H_2^\Gamma(X)$ and $p\in H_0^\Gamma(X)$ is a lift of the class of a point, $\mu$ is defined using the equivariant Chern classes of $\E$, and $\int_{M^X_H(c_1,d)}$ means pushforward to a point
in equivariant cohomology. We assume that the  moduli spaces $M^X_H(c_1,d)$
 have dimension equal to the expected dimension $d$.\begin{NB} added 7.9\end{NB}
If $\E$ is an equivariant torsion-free sheaf of rank $r$ we define
$\cho(\E):=\ch(\E)e^{-\frac{c_1(\E)}{r}}$, where we now use equivariant
Chern character and first Chern class and define $\mu_\rho(\beta):=(-1)^\rho \cho_{\rho+1}(\E)/\beta$.
Let $b_1,\ldots,b_s$ be a homogeneous basis of $H_\Gamma^*(X)$ as a free 
$\Q[\ve_1,\ve_2]$-module.
For all $\rho\ge 1$ let $\tau_1^\rho, \ldots \tau_s^\rho$ be indeterminates,
and put $\alpha_\rho:=\sum_{k=1}^s a_k^\rho b_k \tau_k^\rho$, with 
$a_k^\rho\in\Q[\ve_1,\ve_2]$. 
Using this we define
$\widetilde\Phi_{c_1}^{H}\big(\exp\big(\sum_{\rho}\alpha_\rho\big)\big)
\in \Q[\ve_1,\ve_2][[\Lambda,( \tau_k^\rho)]]$
by the right-hand side of (\ref{donna}).
As the equivariant universal sheaf is unique up to twist by an equivariant line bundle,
$\widetilde\Phi_{c_1}^{H}(\exp(\alpha z +p x))$ and 
$\widetilde\Phi_{c_1}^{H}\big(\exp\big(\sum_{\rho}\alpha_\rho\big)\big)$ 
are independent of the choice of equivariant universal bundle.

\begin{NB}
Let $T$ be the torus in the toric surface $X$.

Let $Q$ be an open subset of a suitable quot-scheme 
such that the moduli space is $M:=Q/PGL(V)$. We may assume that
$V={\cal O}_X(-nH)^{\oplus N}$ and 
${\cal O}_X(-nH)$ is $T$-linearized.
Let $q:V \otimes {\cal O}_Q \to {\cal Q}$ be the universal quotient.
For $t \in T$,
we have a family of quotient
$$
q_t:V \otimes {\cal O}_Q \cong (t^{-1})^*(V) \otimes {\cal O}_Q
\to (t^{-1})^*({\cal Q}).
$$
Since $T$ is connected, 
$(t^{-1})^*({\cal Q})$ is a deformation of
${\cal Q}$. Hence 
we have a (unique) morphism 
$f_t:Q \to Q$ with a commutative diagram:
\begin{equation}
\begin{CD}
V \otimes {\cal O}_Q  @>{q_t}>> (t^{-1})^*({\cal Q})\\
@| @VVV\\
V \otimes {\cal O}_Q  @>{f_t^*(q)}>> f_t^*({\cal Q})
\end{CD}
\end{equation}
$f_t$ gives an action of $T$ on $Q$ and 
${\cal Q}$ is equivariant with respect to the diagonal
action $t \mapsto (t,f_t) \in Aut(X \times Q)$.
This action is commutative with the action of $GL(V)$.
Hence it induces an action on the moduli space
$M$.

In order to have a universal family over $M$,
we have to find a $T \times GL(V)$-equivariant line bundle
$L$ on $Q$ such that $\lambda id_V \in GL(V)$ acts as
a multiplication by $\lambda$.
Let $E$ be an element of $K(X)$ corresponding to $M$.
If $\chi(*,E):K^T(X) \to K(X) \to {\Bbb Z}$ is surjective,
then we have such a line bundle $L$ as usual.   
Thus I think that there is an equivariant universal family
if and only if there is a universal family.
Am I right?
If this is true, then there is a quasi-universal family
in any case:
since ${\cal Q} \otimes (V(nH))^{\vee}$ is $PGL(V)$-equivariant,
it gives a quasi-universal family.
\end{NB}

We cannot hope to extend this naive definition without the assumption
that the moduli spaces are of expected dimensions. This is because 
we can blowup only at the fixed points of the torus action and cannot
avoid the support of the cycles representing $\alpha_\rho$. Here we
probably need to use virtual fundamental classes as in \cite{Moc}.
Then to prove that its specialization coincides with the ordinary
invariants, we need to prove the blowup formula in the context of
virtual fundamental classes.
\begin{NB} moved 7.9 LG \end{NB}

Let $\xi\in H^2(X,\Z)$ be an equivariant lifting of a good class of type $c_1$. 
Then $\I_1$, $\I_2$, $\AA_{\xi,+}$ and $\AA_{\xi,-}$ are in a natural way equivariant sheaves on $X_2^{[l]}$, and the equivariant wallcrossing terms 
$\widetilde\delta_{\xi,t}^X(\exp(\sum_{\rho\ge 1} \alpha_\rho))$,
$\widetilde \delta_{\xi,t}^X(\exp(\alpha z+px))$  are defined by the right-hand side of 
formulas (\ref{wallct}), where now 
$\int_{X_2^{[l]}}$ stands for  equivariant pushforward to a point, and 
\begin{gather*}
\widetilde\delta_\xi^X\Big(\exp\Big(\sum_{\rho\ge 1} \alpha_\rho\Big)\Big):= 
\Big[\widetilde\delta_{\xi,t}^X\Big(\exp\Big(\sum_{\rho\ge 1} \alpha_\rho\Big)\Big)\Big]_{t^{-1}}
,\\
\widetilde \delta_\xi^X(\exp(\alpha z+px)):=
[\widetilde \delta_{\xi,t}^X(\exp(\alpha z+px))]_{t^{-1}}.
\end{gather*}

By (\ref{cct}) we see that
$\widetilde\delta_{\xi,t}^X(\exp(\sum_{\rho\ge 1} \alpha_\rho))\in 
\Lambda^{-\xi^2-3}\Q[\ve_1,\ve_2]((t^{-1}))[[\Lambda,(\tau_k^\rho)]]$.
Thus \linebreak[1]
$\widetilde\delta_{\xi}^X(\exp(\linebreak[2]\sum_{\rho\ge 1} \alpha_\rho))\in 
\Q[\ve_1,\ve_2][[\Lambda,(\tau_k^\rho)]]$ 
and by definition 
$\widetilde\delta_{\xi,t}^X(\exp(\sum_{\rho\ge 1} \alpha_\rho))|_{\ve_1=\ve_2=0}\linebreak[1]=
\delta_{\xi,t}^X(\exp(\linebreak[2]\sum_{\rho\ge 1} \alpha_\rho))$.
\end{Definition}

\begin{Remark}
Note that the coefficient of $\Lambda^d$ in $\widetilde \delta_\xi^X(\exp(\sum_{\rho}
\alpha_\rho))$  is not a polynomial of weight $d$ in the $\tau^\rho_k$ 
(as in Remark \ref{weight}) but has contributions of different weights. 
Arguing as in Rem.~\ref{weight} one sees that 
the coefficient of $\Lambda^d$ is a sum of terms of weight $\ge d$.
Thus the variable $\Lambda$ in the definition of $\widetilde \delta_\xi^X$ 
is not redundant.
\end{Remark}

Under the assumptions of Proposition \ref{wallcr}
let $\widetilde B_+$ be a set consisting of one equivariant lift $\xi$ for  each class
of type $c_1$ with $\<\xi, H_+\> >0>\< \xi, H_-\>$. Then the same proof as before (with all sheaves and classes replaced by the equivariant versions)  shows that the statement of the proposition
holds with $\Phi^{H_+}_{c_1}$, $\Phi^{H_-}_{c_1}$,
$B_+$, $\delta_{\xi}^X$ replaced by $\widetilde \Phi^{H_+}_{c_1}$, $\widetilde\Phi^{H_-}_{c_1}$,
$\widetilde B_+$, $\widetilde \delta_{\xi}^X$ respectively, i.e. the wallcrossing of the equivariant Donaldson invariants is given by the equivariant wallcrossing terms.

In this section we want to give a formula expressing $\widetilde \delta^X_{\xi,t}$ in terms 
of the Nekrasov partition function $Z$.

\begin{Theorem}\label{main1}
$$
\widetilde \delta_{\xi,t}^X
\Big(\exp\Big(\sum_{\rho} \alpha_\rho \Big)\Big)=
\frac{1}{\Lambda}\exp\Big(\sum_{i=1}^\chi F\big(w(x_i),w(y_i),
\hbox{$\frac{t-\iota_{p_i}^*\xi}{2}$};\Lambda,((-1)^\rho\iota_{p_i}^*\alpha_\rho)_\rho\big)\Big).
$$
\end{Theorem}

Note that the left-hand side 
lies in $\Lambda^{-\xi^2-3}\Q[\ve_1,\ve_2]((t^{-1}))[[\Lambda,(\tau_k^\rho)]]$. 
In the course
of the proof we will also have to show how one can interpret the right-hand
side, so that both sides lie in the same ring. 

It is tempting  to write Theorem \ref{main1} as
$$
\widetilde \delta_{\xi,t}^X
\Big(\exp\Big(\sum_{\rho} \alpha_\rho\Big)\Big)=
\frac{1}{\Lambda}\prod_{i=1}^\chi Z\big(w(x_i),w(y_i),
\hbox{$\frac{t-\iota_{p_i}^*\xi}{2}$};\Lambda,((-1)^\rho\iota_{p_i}^*\alpha_\rho)_\rho\big),
$$
but it appears difficult to give a meaning to the  right-hand side of this equation 
(other than as an abbreviation for the right-hand side of Theorem \ref{main1}).

As a first step we will show that, up to a correction term, there is an 
expression for  $\widetilde \delta^X_{\xi,t}$ in terms of the instanton
part of the partition function. In a second step we will 
see that this correction term is accounted for by the perturbation part.

\subsection{The instanton part}
We start by reviewing  some results and definitions from \cite{EG2}.
The fixed points of the $\Gamma$-action on $X_2^{[l]}$ are the pairs 
$(Z_1,Z_2)$ of zero-dimensional subschemes with support in 
$\{p_1,\ldots,p_\chi\}$ with $\len(Z_1)+\len(Z_2)=l$ and such that each
$I_{Z_\alpha,p_i}$ is generated by monomials in $x_i,y_i$.
We associate to 
$(Z_1,Z_2)$ the $\chi$-tuple $(\vec Y^1,\ldots,\vec Y^\chi)$ with 
$\vec Y^i=(Y^i_1,Y^i_2)$, where
$$Y_\alpha^i=\big\{ (n,m)\in \Z_{>0}\times \Z_{>0} \bigm|
x_i^{n-1}y_i^{m-1}\not\in I_{Z_\alpha,p_i}\big\}.$$ 
We write $|Y^i_\alpha|$ for the number of elements of $Y^i_\alpha$ and $|\vec Y^i|:=|Y^i_1|+|Y^i_2|$.
This gives a bijection from the fixed point set $(X_2^{[l]})^\Gamma$ to the 
set of the $\chi$-tuples of pairs of Young diagrams
$(\vec Y^1,\ldots,\vec Y^\chi)$, with $\sum_{i} |Y^i|=l$.

We denote also by $x_i,y_i$ the one dimensional $\Gamma$-modules given by
$t\cdot x_i= e^{-w(x_i)} x_i$, $t\cdot y_i= e^{-w(y_i)} x_i$.
\begin{NB} Corrected 6.9.LG\end{NB}
If $L$ is an equivariant line bundle on $X$, the fibre $L(p_i)$ at a fixed point
and the cohomology groups $H^i(X,L)$ are  in a natural way  $\Gamma$-modules.

The following follows easily from 
the definition of
$N_{\alpha,\beta}^{\vec Y^i}(\ve_1,\ve_2,a)$ in  \cite{NY1} and 
\cite[Lemma 3.2]{EG2}. In fact it is basically a reformulation of
\cite[Lemma  3.2]{EG2} and a straightforward generalization of \cite[Thm.~3.4]{NY1}.
In order to get the correct result one has to take into account the following:
\begin{enumerate}
\item 
The formulas in \cite{EG2} are for $V=L\oplus (L\otimes K_X)$ instead of $L$. 
But the proof only uses that $H^0(V)=H^2(V)=0$. 
\item  Our convention for the $\Gamma$-action on $X^{[l]}_{2}$ differs from that in 
\cite{EG2}, which is  
$t\cdot  
(\I_{Y_1},\I_{Y_2})=(F_t^*\I_{Y_1},F_t^*\I_{Y_2})$). This changes the $\Gamma$-modules $x_i,y_i$ to $x_i^{-1},y_i^{-1}$.
\item 
In \cite[Thm.~3.4]{NY1} the case of $\widehat \C^2$ was studied and the argument
shows that
$\Ext^1(\I_{Z_\alpha},\linebreak[2]\I_{Z_\beta}\otimes L)= H^1(X,L)+ \bigoplus_i 
\Ext^1(\I_{Z_\alpha,p_i},\I_{Z_\beta,p_i}(-\ell_\infty))\otimes L(p_i)$.
\end{enumerate}

\begin{Lemma}
Let $(Z_1,Z_2)\in (X_2^{[l]})^\Gamma$ correspond to $(\vec Y^1,\ldots,\vec Y^\chi)$ 
under the  above bijection. Let $L$ be an $\Gamma$-equivariant 
line bundle on $X$, such that $c_1(L)$ is good. 
We have 
in the Grothendieck group of $\Gamma$-modules
\begin{align} 
T_{(Z_1,Z_2)} X_2^{[l]}&=
\sum_{i=1}^\chi\sum_{\gamma=1}^2 N_{\gamma,\gamma}^{\vec Y^i}(x_i,y_i,L(p_i)^{-\frac{1}{2}}),\\
\Ext^1(\I_{Z_2},\I_{Z_1}\otimes L)&=
H^1(X,L)+\sum_{i=1}^\chi N_{2,1}^{\vec Y^i}(x_i,y_i,L(p_i)^{-\frac{1}{2}}),\\
\Ext^1(\I_{Z_1},\I_{Z_2}\otimes L^{-1})&=
H^1(X,L^{-1})+\sum_{i=1}^\chi N_{1,2}^{\vec Y^i}(x_i,y_i,L(p_i)^{-\frac{1}{2}}).
\end{align}

\end{Lemma}
\begin{NB}\begin{proof}
Let $\vec Y=(Y_1,Y_2)$, be a pair of Young diagrams. 
Let $a'_i$ be the length of the $i$-th row of $Y_1$,
and $b'_i$ the length of the $i$-th row of $Y_2$, and let 
$a_i$ (resp. $b_i$ ) be the length  of the $i$-th column of $Y_1$ (resp. $Y_2$). 
Then $a'=(a'_i)_i$ and $b'=(b'_i)_i$ are partitions and 
$a=(a_j)_j$ and $b=(b_j)_j$ the corresponding  dual partitions.
Let 
$$I_{W_1}=(x^{a'_1},yx^{a'_2},y^2x^{a'_3},\ldots), \quad 
I_{W_2}=(x^{b'_1},yx^{b'_2},y^2x^{b'_3},\ldots).$$
These are the ideals of zero dimensional subschemes of $\A^2$.
Following \cite{EG2} we put 
$$E_{a',b'}(x,y):=
\sum_{1\le i\le j}\Big(\sum_{p=a'_{j+1}+1}^{a'_j}y^{j+1-i}x^{p-b'_i}+ 
\sum_{p=b'_{j+1}+1}^{b'_j}y^{i-j}x^{a'_i-p+1}\Big), $$
and put $E_{Y_1,Y_2}(x,y):=E_{a',b'}(x,y)$.
This coincides with the definition in \cite{EG2}, Lemma 3.2 
except for replacing
$x,y$ by $y^{-1},x^{-1}$ and replacing $a'_{i-1}$ by $a'_i$ and $b'_{i-1}$
 by $b'_i$.
The map $(p,i)\mapsto (i,a_p,p)$ is a bijection
$$Y_1\to \big\{ (i,j,p)\bigm| 1\le i\le j,\ a'_{j+1}<p\le a'_j\big\}$$
with inverse 
$(i,j,p)\mapsto (p,i)$.
Similarly $(p,i)\mapsto (i,b_p,p)$ is a bijection 
$$Y_2\to \big\{ (i,j,p)\bigm| 1\le i\le j,\ b_{j+1}<p\le b_j\big\}.$$
Thus we see  for $\alpha,\beta\in \{1,2\}$ that 
\begin{align*}
E_{Y_\alpha,Y_\beta}(y^{-1},x^{-1})
&=\sum_{s=(p,j)\in Y_\alpha} y^{a_p+1-j}x^{p-b_j'}+
\sum_{s=(p,j)\in Y_\beta}y^{j-b_p}x^{a_j'-p+1}\\
&=\sum_{s\in Y_\alpha} x^{-l_{Y_\beta}(s)}y^{a_{Y_\alpha}(s)+1}+
\sum_{s\in Y_\beta} x^{l_{Y_\alpha}(s)+1}y^{-a_{Y_\beta}(s)}=N_{\alpha,\beta}(x,y,0).
\end{align*}
By \cite{EG2} Lemma 3.2 we have in the Grothendieck group of $\Gamma$-modules 
\begin{align*}
T_{(Z_1,Z_2)}X_2^{[l]}&=
\sum_{i=1}^\chi \sum_{\gamma=1}^2 E_{Y_\gamma^i,Y_\gamma^i}(y_i^{-1},x_i^{-1})=
\sum_{i=1}^\chi \sum_{\gamma=1}^2 N_{\gamma,\gamma }^{\vec Y_i}(x_i,y_i,L(p_i)^{-
\frac12})\\
\Ext^1(\I_{Z_2},\I_{Z_1}\otimes L)&= 
H^1(X,L)+
\sum_{i=1}^\chi L(p_i) E_{Y_\alpha^i,Y_\beta^i}(y_i^{-1},x_i^{-1}),\\
&=H^1(X,L)+
\sum_{i=1}^\chi N_{\alpha,\beta}^{\vec Y_i}(x_i,y_i,L(p_i)^{-\frac{1}{2}}).\\
\Ext^1(\I_{Z_1},\I_{Z_2}\otimes L^{-1})&= 
H^1(X,L^{-1})+
\sum_{i=1}^\chi L^{-1}(p_i) E_{Y_1^i,Y_2^i}(y_i^{-1},x_i^{-1}),\\
&=
H^1(X,L^{-1})+\sum_{i=1}^\chi N_{1,2}^{\vec Y_i}(x_i,y_i,L(p_i)^{-\frac{1}{2}}).
\end{align*}
The formulas in \cite{EG2} are for $V=L\oplus (L\otimes K_X)$ instead of $L$. 
But the proof only uses that $H^0(V)=H^2(V)=0$. In \cite{EG2} the formulas contain 
$E_{Y_\alpha^i,Y_\beta^i}(y_i^{-1},x^{-1}_i)$ instead of 
$E_{Y_\alpha^i,Y_\beta^i}(x_i,y_i)$. 
Note however that now the action of $\Gamma$ is by 
$(F^{-1}_{t_1,t_2})^*$ instead of $(F_{t_1,t_2})^*$ in \cite{EG2}, which 
changes $x_i,y_i$ to $x_i^{-1},y_i^{-1}$, and that with respect to \cite{EG2},
we have exchanged the role 
of $x$ and $y$ in the definition of $E_{Y_j^i,Y_k^i}(x,y)$.
\end{proof}
\end{NB}

Let $F=\sum_{i=1}^r F_i$ be a decomposition of a $\Gamma$-module into
$1$-dimensional modules in the Grothendieck group of $\Gamma$-modules,
and let $w(F_i)$ be the weight of $F_i$. 
Then in the equivariant cohomology we get
$c^{t}(F)=\prod_{i=1}^r (w(F_i)+t)$. Thus we have the following corollary.

\begin{Corollary}\label{tann}
Let $(Z_1,Z_2)\in (X_2^{[l]})^{\Gamma}$ correspond to $(\vec Y^1,\ldots, \vec Y^\chi)$.
Write $L$ for the equivariant line bundle on $X$ whose first Chern class is 
(our chosen lifting of) $\xi$. Then in $\Q[\ve_1,\ve_2,t]$,
 we have the identity
\begin{align*}
e(T_{(Z_1,Z_2)}X_2^{[l]})&c^{-t}(\Ext^1(\I_{Z_2},\I_{Z_1}\otimes L))
c^{t}(\Ext^1(\I_{Z_1},\I_{Z_2}\otimes L^\vee))\\&=
c^{-t}(H^1(X,L))c^t(H^1(X,L^\vee))\prod_{i=1}^\chi \prod_{\alpha,\beta=1}^2n^{\vec Y^i}_{\alpha,\beta}
\big(w(x_i),w(y_i),
\hbox{$\frac{t-\iota_{p_i}^*\xi}{2}$}\big).
\end{align*}
\end{Corollary}

\begin{Lemma}\label{instan}In 
$\Lambda^{-\xi^2-3}\Q(\ve_1,\ve_2)((t^{-1}))[[\Lambda,(\tau_k^\rho)]]$ we have
\begin{align*}
\widetilde \delta_{\xi,t}^X
\Big(\exp\Big(\sum_{\rho} \alpha_\rho\Big)\Big)&=
 \frac{\prod_{i=1}^\chi \Zin\big(w(x_i),w(y_i),
\frac{t-\iota_{p_i}^*\xi}{2};\Lambda,((-1)^\rho \iota_{p_i}^*\alpha_\rho)_\rho\big)
} {\Lambda^{\xi^2+3}c^{-t}(H^1(X,L))c^{t}(H^1(X,L^\vee))}.
\end{align*}
\end{Lemma}
\begin{Remark} By Definition \ref{eqwall} the left-hand side is an element of 
$\Lambda^{-\xi^2-3}\Q[\ve_1,\ve_2]((t^{-1}))\linebreak[2][[\Lambda,(\tau_k^\rho)]]$.
Note that by  (\ref{Zinst}) the right-hand side
is an element of $\Lambda^{-\xi^2-3}\Q(\ve_1,\ve_2,t)\linebreak[2][[\Lambda,(\tau_k^\rho)]]$. Using 
$\frac{1}{b+t}=t^{-1}(\sum b^it^{-i})$, we can view $\frac{1}{\prod_{\alpha,\beta=1}^2 n_{\alpha,\beta}(w(x_i),w(y_i),(t-\iota_{p_i}^* \xi)/2)}$ as an element of 
$\Q(\ve_1,\ve_2)[[t^{-1}]]$. Then by (\ref{cct}) the right-hand side of the Lemma is
interpreted as an element of $\Lambda^{-\xi^2-3}\Q(\ve_1,\ve_2)((t^{-1}))\linebreak[2][[\Lambda,(\tau_k^\rho)]]$, and
we will show that the equality  holds here. The lemma 
shows that the right-hand
side lies even in $\Lambda^{-\xi^2-3}\Q[\ve_1,\ve_2]((t^{-1}))[[\Lambda,(\tau_k^\rho)]]$.
\end{Remark}
\begin{proof}
Let $(Z_1,Z_2)\in (X_2^{[l]})^\Gamma$ correspond to 
$(\vec Y^1,\ldots,\vec Y^\chi)$.
Let $\alpha\in \{1,2\},$ and let $p_i\in X^{\Gamma}$.
We claim that  
$$\iota_{(p_i,(Z_1,Z_2))}^*\ch(\I_\alpha)=1-(1-e^{-w(x_i)})(1-e^{-w(y_i)})\cdot
\sum_{s\in Y_{\alpha}^i} e^{-l'(s)w(x_i)-a'(s)w(y_i)}.$$
\begin{NB} Corrected 6.9. LG\end{NB}
Let $\oo_1$ (resp.\ $\oo_2$)
be the sheaf on $X\times X_2^{[l]}$ whose restriction
to $X\times X^{[n]}\times X^{[m]}$ is the pushforward of $\oo_{Z_n(X)}$
(resp. $\oo_{Z_m(X)}$) via the inclusion.
For $\alpha=1,2$  we have 
$\iota_{X\times (Z_1,Z_2)}^*(\oo_\alpha)=
\sum_{i=1}^\chi (\iota_{p_i})_*(\oo_{Z_\alpha,p_i}).$
By definition an equivariant basis of $\oo_{Z_\alpha,p_i}$ is 
$\big\{x_i^{n-1}y_i^{m-1}\bigm| (n,m)\in Y^i_\alpha\big\}$,
thus $\oo_{Z_\alpha,p_i}=\sum_{s\in Y^i_\alpha}x_i^{l'(s)} y_i^{a'(s)}$ as
$\Gamma$-modules.
By 
localization we get 
\begin{align*} i^*_{(p_i,(Z_1,Z_2))}\ch(\oo_\alpha)&=\iota_{p_i}^*\big(\sum_{j=1}^\chi (\iota_{p_j})_
*\ch(\oo_{Z_\alpha,p_j})\Big)
=\iota_{p_i}^*(\iota_{p_i})_*\ch(\oo_{Z_\alpha,p_i})
\\ &=(1-e^{-w(x_i)})(1-e^{-w(y_i)})\ch(\oo_{Z_\alpha,p_i}).
\end{align*}
\begin{NB} Corrected 6.9.LG\end{NB}
Using that $\ch(\I_\alpha)=
1-\ch(\oo_\alpha)$, we get the claim. 

We put $f_{1}:=\frac{\xi-t}{2}$, $f_2:=\frac{t-\xi}{2}$. Then the claim implies
\begin{align*}
\iota_{(p_i,(Z_1,Z_2))}^*&\big(\ch(\I_1)e^{\frac{\xi-t}2}+\ch(\I_2)e^{\frac{t-\xi}2}
\big)\\
&=\sum_{\alpha=1}^2
e^{\iota_{p_i}^*(f_\alpha)}\Big(1-(1-e^{-w(x_i)})(1-e^{-w(y_i)})
\sum_{s\in Y_{\alpha}^i} e^{-l'(s)w(x_i)-a'(s)w(y_i)}
\Big).
\end{align*}
\begin{NB} Corrected 6.9. LG\end{NB}
By the definition (\ref{FY}) of $E^{\vec Y}(\ve_1,\ve_2,a,\vec\tau)$
this gives
\begin{equation}\label{eq:EY}
\begin{split}
\iota_{(Z_1,Z_2)}^*\exp\Big(\sum_\rho(-1)^\rho \big[\ch(\I_1)&
e^{\frac{\xi-t}2}\oplus 
\ch(\I_2)e^{\frac{t-\xi}2}\big]_{\rho+1}/\alpha_\rho\Big)\\
&=\prod_{i=1}^\chi E^{\vec Y}(w(x_i),w(y_i),\hbox{$\frac{t-\iota_{p_i}^*\xi}2$},
((-1)^\rho \iota_{p_i}^*\alpha_\rho)_\rho).
\end{split}
\end{equation}
Write $|Y|:=|\vec Y_1|+\ldots+|\vec Y_\chi|$, and 
write $(Z^Y_1,Z^Y_2)$ for the point of  $X_2^{[|Y|]}$ 
determined by an $\chi$-tuple $Y= (\vec Y_1,\ldots,\vec Y_\chi)$ of pairs of
Young diagrams. Then we get by localization
\begin{align*}\widetilde\delta_{\xi,t}^X&
\Big(\exp\Big(\sum_\rho \alpha_\rho\Big)\Big)
\\
&=\sum_{Y=(\vec Y_1,\ldots,\vec Y_\chi)}
\frac{\Lambda^{4|Y|-\xi^3-3}\prod_{i=1}^\chi E^{\vec Y}\big(w(x_i),w(y_i),
\frac{t-i_{p_i}^*\xi}2,
((-1)^\rho \iota_{p_i}^*\alpha_\rho)_\rho\big)} 
{e(T_{(Z_1^Y,Z^Y_2)}X_2^{[|Y|]})c^{-t}(\Ext^1(\I_{Z_2^Y},\I_{Z_1^Y}
\otimes L))
c^{t}(\Ext^1(\I_{Z_1^Y},\I_{Z_2^Y}\otimes L^\vee))}\\
&=\Lambda^{-\xi^3-3}\frac{\prod_{i=1}^\chi
\Zin \big(w(x_i),w(y_i),\frac{t-\iota_{p_i}^*\xi}{2};\Lambda, 
((-1)^\rho \iota_{p_i}^*\alpha_\rho)_\rho\big)}{c^{-t}(H^1(X,L)) c^{t}
(H^1(X,L^\vee))},
\end{align*}
where the last step is by Cor. \ref{tann}.
\end{proof}

\subsection{The perturbation part}
Now we want to identify the contribution of the perturbation part.
We first need to review the perturbation part of the $K$-theoretic 
Nekrasov partition function from \cite[section 4.2]{NY3}.
We set
\begin{equation}\label{gammati}\begin{split}
\gamma_{\ve_1,\ve_2}(x|\beta;\Lambda)&:=
\begin{aligned}[t]
\frac{1}{2\ve_1\ve_2}\left(
-\frac{\beta}{6}\left(x+\frac{1}{2}(\ve_1+\ve_2)\right)^3 
+x^2\log(\beta\Lambda)\right)\qquad\\
+\sum_{n \geq 1}\frac{1}{n}
\frac{e^{-\beta nx}}{(e^{\beta n\ve_1}-1)(e^{\beta n\ve_2}-1)},
\end{aligned}
\\
\widetilde{\gamma}_{\ve_1,\ve_2}(x|\beta;\Lambda)
&:=\gamma_{\ve_1,\ve_2}(x|\beta;\Lambda)+
\frac{1}{\ve_1 \ve_2} \left(\frac{\pi^2 x}{6 \beta}-\frac{\zeta(3)}{\beta^2} \right)\\
&\qquad \qquad +\frac{\ve_1+\ve_2}{2\ve_1 \ve_2} 
\left( x \log (\beta \Lambda)+\frac{\pi^2}{6\beta} \right)+
\frac{\ve_1^2+\ve_2^2+3\ve_1 \ve_2}{12 \ve_1 \ve_2} \log(\beta\Lambda)
\end{split}
\end{equation}
for $(x,\beta,\Lambda)$ in a neighbourhood  of $\sqrt{-1}\R_{>0}\times\sqrt{-1}\R_{<0}
\times\sqrt{-1}\R_{>0}$.
We formally expand $\ve_1 \ve_2\widetilde{\gamma}_{\ve_1,\ve_2}(x|\beta;\Lambda)$
as a power series of $\ve_1,\ve_2$ (around $\ve_1=\ve_2=0$).
By the expansion (\ref{cm})
we obtain
$$
\sum_{n \geq 1}\frac{1}{n}
\frac{e^{-\beta nx}}{(e^{\beta n\ve_1}-1)(e^{\beta n\ve_2}-1)}
=\sum_{m \geq 0} \frac{c_m}{m!}\beta^{m-2} \mathrm{Li}_{3-m}(e^{-\beta x}),
$$
where $\Li_{3-m}$ is the polylogarithm (see \cite[Appendix B]{NY3}  for details).
Here we choose the branch of $\log$ by $\log(r\cdot e^{i\phi})=\log(r)+i\phi$
with $\log(r)\in \R$
for $\phi\in (-\pi/2,3\pi/2)$ and $r\in \R$.
We define  $\gamma_{\ve_1,\ve_2}(-x|\beta;\Lambda)$  by analytic continuation along circles in a 
counter-clockwise way. 
Finally we define 
$$\Fper_K(\ve_1,\ve_2,x|\beta;\Lambda):=-\gamma_{\ve_1,\ve_2}(2x|\beta;\Lambda)
-\gamma_{\ve_1,\ve_2}(-2x|\beta;\Lambda).$$
Then $\Fper_K(\ve_1,\ve_2,x|\beta;\Lambda)$ is  a formal power series in $\ve_1,\ve_2$ whose coefficients are 
holomorphic functions in $\Lambda\in \C\setminus \sqrt{-1}\R_{\le 0}$, $x\in \C\setminus \sqrt{-1}\R_{\le 0}$, $\beta\in \C$ with $|\beta|<\frac{\pi}{|x|}$.
In \cite[section 4.2]{NY3} it is shown that $\Fper_K(\ve_1,\ve_2,x|\beta;\Lambda)$ 
converges to $\Fper(\ve_1,\ve_2,x;\Lambda)$ when $\beta$ goes to $0$.

We will use the following obvious consequence of the localization formula on $X$.
\begin{Remark}\label{cancel}
For any class $\gamma\in H^j_\Gamma(X)$ we have 
$$\sum_{i=1}^\chi \frac{\iota_{p_i}^*\gamma}{w(x_i)w(y_i)}=\int_X \gamma\in H^{j-4}_\Gamma(pt).$$
In particular if $\gamma=1$ or $\gamma\in H^2_\Gamma(X)$, then 
$
\sum_{i=1}^\chi \frac{\iota_{p_i}^*\gamma}{w(x_i)w(y_i)}=0.
$
\end{Remark}

\begin{Lemma}\label{pertur}
\begin{multline*}
\sum_{i=1}^\chi
\Fper\big(w(x_i),w(y_i),\hbox{$\frac{t-\iota_{p_i}^*\xi}{2}$};\Lambda\big)
\\ =
(-\xi^2-2) \log\Lambda -\log(c^{-t}(H^1(X,L)))-\log(c^{t}(H^1(X,L^\vee)))
\end{multline*}
holds in $\oo[[\ve_1,\ve_2]]$, where $\oo$ denotes the
holomorphic functions in $(t,\Lambda)$ on $(\C\setminus \sqrt{-1} \R_{\le 0})^2$. 
\end{Lemma}

A priori, the left-hand side lives in $\frac{1}{\prod_{i=1}^\chi w(x_i)w(y_i)}\oo[[\ve_1,\ve_2]]$,
but in the course of the proof we show that it is, in fact, in $\oo[[\ve_1,\ve_2]]$,
and the equality holds in $\oo[[\ve_1,\ve_2]]$.
In $\oo[[\ve_1,\ve_2]]$ we can take the exponential of both sides of the equation. Note that the exponential of the right-hand side also lives in $\Lambda^{-\xi^2-2}\Q[\ve_1,\ve_2]((t^{-1}))[[\Lambda,(\tau^\rho_k)]]$. With this remark 
Lemma~\ref{pertur} and Lemma~\ref{instan} together imply  Theorem~\ref{main1}.

\begin{proof}
Let $L$ be an equivariant line bundle on $X$ whose equivariant first Chern class is $\xi$.
In particular $H^i(X,L)=0$ and $H^i(X,L^\vee)=0$  for $i\ne 1$.
Let $\ell=h^1(X,L)$, $\ell'=h^1(X,L^\vee)$, and let $\alpha_1,\ldots,\alpha_\ell$
(resp. $\alpha'_1,\ldots,\alpha'_{\ell'}$)
 be the weights of $\Gamma$ on $H^1(X,L)$ (resp. $H^1(X,L^\vee)$). 
Then in 
$\Gamma$-equivariant cohomology we get 
$$c^{-t}(H^1(X,L))=\prod_{j=1}^\ell (\alpha_j-t), \qquad 
c^{t}(H^1(X,L^\vee))=\prod_{k=1}^{\ell'} (\alpha'_k+t).$$ 
Write $p:X\to pt$ for the map to a point. Then 
the Riemann-Roch theorem gives 
$$\sum_{k=1}^{\ell'} e^{\alpha'_k+t}=-\ch(p_!(L^\vee))e^{t}=
-\sum_{k=1}^\chi \frac{e^{-i_{p_i}^*\xi+t}}{(1-e^{-w(x_i)})(1-e^{-w(y_i)})}.$$
 Thus we get 
\begin{equation}\label{eq:RR}
\begin{split}
\sum_{i=1}^\chi 
\sum_{n>0} \frac{e^{(-i_{p_i}^*\xi+t)n(-\beta)}}
{n(1-e^{-w(x_i)n(-\beta)})(1-e^{-w(y_i)n(-\beta)})}
&=-\sum_{n>0} \sum_{k=1}^{\ell'} \frac{e^{(\alpha'_k+t)n(-\beta)}}{n}\\
&= \sum_{k=1}^{\ell'} \log(1-e^{-(\alpha_k'+t)\beta})
\end{split}
\end{equation}
in ${\widetilde \oo}[[\ve_1,\ve_2]][\prod_i (w(x_i)w(y_i))^{-1}]$,
where ${\widetilde \oo}$ is the ring of holomorphic functions in  $(x,\beta,t)$ 
in a neighborhood of $\sqrt{-1}\R_{>0}\times \sqrt{-1}\R_{< 0}\times \sqrt{-1}\R_{> 0}$.

Now we apply the localization formula on $X$. Using (\ref{gammati}) and Remark
\ref{cancel}  we obtain
\begin{multline}\label{eq:RR2}
\sum_{i=1}^\chi \widetilde{\gamma}_{w(x_i),w(y_i)}
(-i_{p_i}^*\xi+t|\beta;\Lambda)\\=\sum_{k=1}^{\ell'} \log(1-e^{-(\alpha'_k+t)\beta})
+\chi(L^\vee)
\log(\beta \Lambda)
-\frac{\beta}{12}\Big(-\xi+t-\frac{K_X}{2}\Big)^3,
\end{multline}
in $\widetilde \oo[[\ve_1,\ve_2]][\prod_i (w(x_i)w(y_i))^{-1}]$.
Here we have used 
$$
\int_X\left(\frac{(-\xi+t-K_X) \cdot (-\xi+t)}{2}+\Todd_2(X)\right)
=\chi(L^\vee),
$$
which follows from Remark \ref{cancel} and the Riemann-Roch theorem. 
\begin{NB}
The same argument shows that 
\begin{equation}\label{eq:RR3}
\sum_{i=1}^\chi \widetilde{\gamma}_{w(x_i),w(y_i)}
(i_{p_i}^*\xi+t|\beta;\Lambda)=\sum_{j=1}^\ell \log(1-e^{-(\alpha_j+t)\beta})+\chi(L) \log(\beta \Lambda)
-\frac{\beta}{12}\int_X\Big(\xi+t-\frac{1}{2}K_X\Big)^3
\end{equation}
in ${\widetilde \oo}[[\ve_1,\ve_2]][\prod_i (w(x_i)w(y_i))^{-1}]$.
\end{NB}

Since $\widetilde{\gamma}_{\ve_1,\ve_2}(-x|\beta;\Lambda)$, 
is defined by an analytic continuation, we derive 
from \eqref{eq:RR2}:  
\begin{multline*}
\sum_i \widetilde{\gamma}_{w(x_i),w(y_i)}
(i_{p_i}^*\xi-t)|\beta;\Lambda)\\=\sum_{k=1}^\ell \log(1-e^{-(\alpha_k-t)\beta})+\chi(L) 
\log(\beta \Lambda)
-\frac{\beta}{12}\int_X(\xi-t-\frac{K_X}{2})^3.
\end{multline*}
Thus we get in $\widetilde \oo[[\ve_1,\ve_2]][\prod_i (w(x_i)w(y_i))^{-1}]$ that 
\begin{equation}
\label{eq:FKpert}
\begin{split}
\sum_{i=1}^\chi& \Fper_K(w(x_i),w(y_i),\hbox{$\frac{t-i_{p_i}^*\xi}{2}$}|\beta;\Lambda)\\
&=\Big(-\sum_i \widetilde{\gamma}_{w(x_i),w(y_i)}
(i_{p_i}^*\xi-t|\beta;\Lambda)-
\sum_i \widetilde{\gamma}_{w(x_i),w(y_i)}
(-i_{p_i}^*\xi+t)|\beta;\Lambda)\Big)\\
&=\Big(-(\chi(L)+\chi(L^\vee))(\log(\Lambda)+\log(\beta))
-\beta\Big(\frac{\<K_X^3\>}{48}+\frac{\<K_X\xi^2\>}{4}-\frac{t\<K_X\xi\>}{2}\Big)\\
&\quad +
\sum_{j=1}^\ell\log\Big(\frac{1}{1-e^{-(\alpha_j-t)\beta}}\Big)
+\sum_{k=1}^{\ell'}\log\Big(\frac{1}{1-e^{-(\alpha'_k+t)\beta}}\Big)\Big).
\end{split}
\end{equation}
As both sides are defined around $\beta=0$, the equality holds there.
Thus we can take $\beta=0$. Using that 
$\lim_{\beta\to 0}\log\big(\frac{\beta}{1-e^{-(\alpha_j-t)\beta}}\big)=\log(\alpha_j-t)$, 
$\lim_{\beta\to 0}\log\big(\frac{\beta}{1-e^{-(\alpha'_k+t)\beta}}\big)=\log(\alpha'_k+t)$, and that 
$\ell=-\chi(L)$, 
$\ell'=-\chi(L^\vee)$,
we obtain
$$ \Fper(w(x_i),w(y_i),\hbox{$\frac{t-i_{p_i}^*\xi}{2}$};\Lambda)
=(-\xi^2-2)\log(\Lambda)-\sum_{j=1}^\ell\log(\alpha_j-t)-\sum_{k=1}^{\ell'}\log(\alpha'_k+t).$$
Note that the right-hand side of this equation is in $\oo[[\ve_1,\ve_2]]$.
Thus, while the individual summands of the left-hand side only lie in 
$\oo[[\ve_1,\ve_2]][\prod_i (w(x_i)w(y_i))^{-1}]$, their sum lies in 
$\oo[[\ve_1,\ve_2]]$.
This shows Lemma  \ref{pertur}.
\end{proof}

Now we want to express the wallcrossing for the Donaldson invariants in
terms of the Nekrasov partition function $Z(\ve_1,\ve_2,a;\Lambda)$. 
This will be necessary because the Nekrasov conjecture determines the lowest order terms in $\ve_1,\ve_2$ of 
$F(\ve_1,\ve_2,a;\Lambda)$, but not of\linebreak[2] $F(\ve_1,\ve_2,\linebreak[1]a;\Lambda,\vec \tau)$.

\begin{Corollary}\label{deltaq}
\begin{multline*}\tag{1}
\widetilde\delta_{\xi,t}^X(\exp(\alpha z+p x))\\=
\frac{1}{\Lambda}\exp\Big(\frac{1}{2}\big\<\Todd_2(X)(\alpha z+ p x)\big\> 
\Big(\sum_{i=1}^\chi
F\big(w(x_i),w(y_i),\hbox{$\frac{t-\iota_{p_i}^*\xi}{2}$};\Lambda e^{\iota_{p_i}^*(\alpha z+px)/4}\big)\Big)\Big),
\end{multline*}
\begin{equation*}\tag{2}
\delta_{\xi,t}^X(\exp(\alpha z+p x))=
\frac{1}{\Lambda}\exp\Big(\sum_{i=1}^\chi
F\big(w(x_i),w(y_i),\hbox{$\frac{t-\iota_{p_i}^*\xi}{2}$};\Lambda e^{\iota_{p_i}^*(\alpha z+px)/4}\big)\Big)|_{\ve_1,\ve_2=0}.
\end{equation*}
\end{Corollary}
\begin{proof}
Let $\vec\tau_1:=(\tau_1,0,0,\ldots)$ be a vector with only the first entry nonzero.
Then in \cite[section 4.5]{NY2} it is shown that 
\begin{equation}\label{tauvsq}
F(\ve_1,\ve_2,a;\Lambda,\vec \tau_1)=
-\frac{\tau_1(\ve_1^2+\ve_2^2+3\ve_1\ve_2)}{24\ve_1\ve_2}+
F(\ve_1,\ve_2,a;\Lambda e^{-\tau_1/4}).
\end{equation}
Thus we get 
\begin{align*}
\widetilde\delta_{\xi,t}^X(\exp(\alpha z+p x))&=
\frac{1}{\Lambda}\exp\Bigg(\sum_{i=1}^\chi\frac{\iota_{p_i}^*(\alpha z+p x)
(w(x_i)^2+w(y_i)^2+3w(x_i)w(y_i))}{24w(x_i)w(y_i)}\Bigg)\\&\qquad \qquad 
\times \exp\Big( \sum_{i=1}^\chi
F\big(w(x_i),w(y_i),\hbox{$\frac{t-\iota_{p_i}^*\xi}{2}$};\Lambda e^{\iota_{p_i}^*(\alpha z+px)/4}
\big)\Big).
\end{align*}
By localization we get 
$$\sum_{i=1}^\chi\frac{\iota_{p_i}^*(\alpha z+p x)
(w(x_i)^2+w(y_i)^2+3w(x_i)w(y_i))}{24w(x_i)w(y_i)}=\frac{1}{2}\big\<
(\alpha z+px) \Todd_2(X)\big\>.$$
This shows (1). (2) follows immediately, because 
$\<(\alpha z+px)\Todd_2(X)\>=0$ in nonequivariant cohomology.
\end{proof}

\section{Explicit formulas in terms of modular forms}
We have expressed the wallcrossing $\delta_\xi^X$ in terms of 
the Nekrasov partition function. Now we want to use the Nekrasov conjecture
to give an explicit formula in terms of the $q$-development of modular forms. 

Let  $q:=e^{2\pi \sqrt{-1} \tau}$ for $\tau\in {\mathcal H}:=\big\{\tau\in \C\bigm| \Im(\tau)>0
\big\}$.
Recall the theta functions
$$\theta_{00}(\tau):=\sum_{n\in \Z} q^{n^2/2},\quad   \theta_{01}(\tau)
:=\sum_{n\in \Z} (-1)^n q^{n^2/2}, \quad  \theta_{10}(\tau):=\sum_{n\in \Z} 
q^{\frac{1}{2}(n+1/2)^2}.$$
Write $E_2(\tau):=1-{24}\sum_{n}\sigma_1(n) q^n$ for the normalized Eisenstein
series of weight $2$.
Denote 
\begin{equation}\label{ududaa}
u:=-\frac{\theta_{00}^4+\theta_{10}^4}{\theta_{00}^2\theta_{10}^2}\Lambda^2, \quad
\frac{du}{da}:=\frac{2 \sqrt{-1}}{\theta_{00}\theta_{10}}\Lambda, 
\quad
a:=\sqrt{-1}\frac{2E_2+\theta_{00}^4+\theta_{10}^4}{3\theta_{00}\theta_{10}}\Lambda.
\end{equation}
Finally put
$$T:=\frac{1}{24}\Big(\frac{du}{da}\Big)^2E_2 -\frac{u}{6}.$$
We can now state the formula for $\delta_\xi^X$ in terms of the $q$-development of
these functions.

\begin{Theorem}\label{main2}
Let $\xi$ be a good class. 
Then 
\begin{multline*}
\delta^X_\xi(\exp(\alpha z +px))\\=\sqrt{-1}^{\<\xi K_X\>-1}\Bigg[
q^{-\frac{1}{2}(\frac{\xi}{2})^2} \exp\Big(\frac{du}{da} \<\alpha,\xi/2\>z+T \<\alpha^2\> z^2- ux\Big)
\Big(\frac{\sqrt{-1}}{\Lambda}\frac{du}{da}\Big)^3\theta_{01}^{K_X^2}\Bigg]_{q^0}.
\end{multline*}
\end{Theorem}
We briefly review the Nekrasov conjecture.
For this we define $u$, $a$ in a different way.
Consider the family of elliptic curves 
$C_u: y^2=(z^2-u)^2-4\Lambda^4$, parametrized by $u\in \C$, which we call the
$u$-plane.  The Seiberg-Witten differential 
$dS:=-\frac{1}{2\pi} \frac{zP'(z)dz}{y}$ is a meromorphic differential
form on $C_u$. For suitable cycles $A,B$ on $C_u$ (for the definition
see \cite[section 2.1]{NY2}, here they are called $A_2,B_2$) put
$a:=\int_{A}dS$, $a^D:=2\pi \sqrt{-1} \int_{B} dS$. These are functions on the $u$-plane
($|u|\gg|\Lambda|$). By definition $a$ and $a^D$ are functions of $u$, but conversely
we will consider $u$ and $a^D$ as functions of $a$ and $\Lambda$.
The {\em period} of $C_u$ is $\tau:=\frac{1}{2\pi \sqrt{-1}}\frac{\partial a^D}{\partial a}$.
The {\em Seiberg-Witten prepotential} ${\cal F}_0$  is the (suitably normalized) 
locally defined function on the $u$-plane with $a^D=-\frac{\partial \F_0}{\partial a}$.
We choose the branch of the logarithm as $\log(re^{i\theta})=\log(r) +i\theta$
for $r\in \R^+$ and $\theta\in (-\pi,\pi)$,
with $\log(r)\in \R$. 
By \cite[sections 2.1, 2.3]{NY2} for $(a,\Lambda)$ in a neighborhood
$U\subset \C\times \C$  of  the set of $(a,\Lambda)\in \sqrt{-1} \R^+\times \sqrt{-1}\R^+$, with $|a|\gg|\Lambda|$,   ${\cal F_0}$ is a holomorphic function of $a$ and $\Lambda$, which we write as 
${\cal F_0}(a;\Lambda)$. 
\begin{NB}
I change the domain of the definition $\Lambda$ to a nbd of
$\sqrt{-1}\R^+$ to match with the convention in the previous
section. This does not match with \cite{NY2}. But the convention was
used only to fix the branch of the logarithm. Apr. 8, HN
\end{NB}
By definition we have 
$\tau=-\frac{1}{2\pi \sqrt{-1}} \frac{\partial^2 \F_0}{(\partial
a)^2}$ and $q=\exp(-\frac{\partial^2 \F_0}{(\partial a)^2})$.
Then with this definition of $\tau$  the formulas  (\ref{ududaa}) hold 
\cite[equation (1.3)]{NY2}.

The Nekrasov conjecture \cite{Nek} 
(proved in \cite{NY1},\cite{NY2},\cite{NO},\cite{BE}) says that
\begin{enumerate}
\item $\ve_1\ve_2 F(\ve_1,\ve_2,a;\Lambda)$ is regular at $\ve_1,\ve_2=0$,
\item $(\ve_1\ve_2F(\ve_1,\ve_2,a;\Lambda))|_{\ve_1=\ve_2=0}={\cal F_0}(a;\Lambda)$.
\end{enumerate}
Here we understand the equation (2) as follows: It is an abbreviation
of two equations, one for the perturbation part and the other for the
instanton part. The former is an equality for holomorphic functions
in $(a,\Lambda)\in U$, and the latter is for formal power series
in $\C[[1/a,\Lambda]]$. Equations appearing below should be understood
in the same way, until the ambiguity of the branch of the logarithm
in the perturbation
part will disappear in the expression.

In \cite{NY2} also the next higher order terms of $F(\ve_1,\ve_2,a,\Lambda)$ in $\ve_1,
\ve_2$ are determined:
We write 
\begin{align}\label{Fdev}
\ve_1\ve_2F(\ve_1,\ve_2,a;\Lambda)&=
\F_0(a;\Lambda)+(\ve_1+\ve_2)H(a;\Lambda)+\ve_1\ve_2A(a;\Lambda)+\frac{\ve_1^2+\ve_2^2}3
B(a;\Lambda)+O
\end{align}
where $O$ stands for terms of degree at least $3$ in $\ve_1$ and $\ve_2$.
It is also proved that $H$, $A$ and $B$ are holomorphic functions on
$U$.
In \cite[section 5.3]{NY2} it is shown that 
$H(a,\Lambda)=\pi \sqrt{-1} a$. By \cite[section 7.1]{NY2} we have 
\begin{equation}\label{AB}\exp(A)=\Big(\frac{2}{\theta_{00}\theta_{10}}\Big)^{1/2}
=\Big(\frac{\sqrt{-1}}{\Lambda} \frac{du}{da}\Big)^{1/2},
\qquad \exp(B -A )=\theta_{01}.
\end{equation}

\begin{Remark}
The sign of $H$ in \cite[section~5.3]{NY2} was wrong and we have
$H(\vec{a},\q,\vec{\tau}) = -\pi\sqrt{-1}\langle\vec{a},\rho\rangle$
in the first displayed formula in [loc.\ cit., p.~66]. Therefore
we have $H(a;\Lambda) = \pi\sqrt{-1}a$ in our case. The mistake
occurred when we take the sum of [loc.\ cit., (E.5)] over
$\alpha<\beta$. Accordingly blowup formulas in [loc.\ cit., section~6]
must be corrected.
\end{Remark}

\begin{pf*}{Proof of Theorem \ref{main2}} 
We  apply the localization formula to $X$. 
Note that $w(x_i),w(y_i)$ are homogeneous of degree $1$ in $\ve_1,\ve_2$.
Furthermore if $\beta\in H^{2i}_\Gamma(X)$, then $\iota_{p_i}^*(\beta)$ is homogeneous
of degree $i$ in $\ve_1,\ve_2$.
Therefore we get the  following expansion (where on the right-hand 
side
we take the values of $\F_0$ and its derivatives at
$(\frac{t}{2},\Lambda)$):
\begin{equation*}
\begin{split}
& \F_0\big(\hbox{$\frac{t-\iota_{p_i}^*(\xi)}{2}$};
\Lambda e^{\iota_{p_i}^*(\alpha z+px)/4}\big)=
\F_0-\frac{1}{2}\frac{\partial\F_0}{\partial a} \iota_{p_i}^*(\xi)+
\frac{\partial F_0}{\partial\log \Lambda} \iota_{p_i}^*(\alpha z+px)/4
\\
  &\qquad\qquad -\frac18\frac{\partial^2 F_0}{\partial a\partial\log\Lambda}
   \iota_{p_i}^*(\xi)\iota_{p_i}^*(\alpha) z
  + \frac1{2}\frac{\partial^2 F_0}{(\partial\log\Lambda)^2}
   \iota_{p_i}^*(\alpha/4)^2 z^2
   + \frac18 \frac{\partial^2 F_0}{\partial a^2}
   \iota_{p_i}^*(\xi)^2
   + O,
  \end{split}
\end{equation*}
where $O$ stands for terms of degree larger than $2$ in 
$\ve_1,\ve_2$.
By Remark \ref{cancel} we see that 
\begin{align}\label{F0}
\sum_{i=1}^\chi \frac{\F_0\big(\hbox{$\frac{t-\iota_{p_i}^*(\xi)}{2}$};
\Lambda e^{\iota_{p_i}^*(\alpha z+px)/4}\big)}{w(x_i)w(y_i)}
=
\begin{aligned}[t]
& \frac{1}{4}\frac{\partial \F_0}{\partial\log \Lambda} x
   - \frac18\frac{\partial^2 \F_0}{\partial a\partial\log\Lambda}
  \<\xi,\alpha\> z
  \\ & \qquad
  + \frac1{32}\frac{\partial^2 \F_0}{(\partial\log\Lambda)^2}
   \<\alpha^2 \>z^2
   + \frac18 \frac{\partial^2 \F_0}{\partial a^2}
   \xi^2
   + O,
\end{aligned}
\end{align}
where again $O$ stands for terms of degree larger than $2$ in $\ve_1,\ve_2$. 
Similarly we get 
\begin{equation}\label{H}\begin{split}
\sum_{i=1}^\chi \frac{w(x_i)+w(y_i)}{w(x_i)w(y_i)}
H(\hbox{$\frac{t-\iota_{p_i}^*\xi}{2}$};\Lambda e^{\iota_{p_i}^*(\alpha z+px)/4})
&=\sum_{i=1}^\chi \pi \sqrt{-1} \frac{w(x_i)+w(y_i) }
{w(x_i)w(y_i)}\frac{t-\iota_{p_i}^*\xi}{2}\\
&= \pi \sqrt{-1} \<\xi/2, K_X\>.
\end{split}\end{equation}
Finally 
\begin{equation}\label{AB1}
\sum_{i=1}^\chi \Big(A(\hbox{$\frac{t-\iota_{p_i}^*\xi}{2}; \Lambda e^{\iota_{p_i}^*(\alpha z+px)}$})
+\frac{w(x_i)^2+w(y_i)^2}{3w(x_i)w(y_i)}B(\hbox{$\frac{t-\iota_{p_i}^*\xi}{2}$};\Lambda  e^{\iota_{p_i}^*(\alpha z+px)/4})
\Big)=\chi A+ \sigma B + O,
\end{equation}
where $\sigma=\frac{1}{3}(c_1(X)^2-2\chi)$ is the signature of $X$, and the
argument of $A,B$ is $(t/2,\Lambda)$.
By the formulas (\ref{Fdev})--(\ref{AB1}) we get that 
$\sum_{i=1}^\chi \frac{1}{w(x_i)w(y_i)}\linebreak[2]
F\big(w(x_i),w(y_i),\hbox{$\frac{t-\iota_{p_i}^*\xi}{2}$},\Lambda 
e^{\iota_{p_i}^*(\alpha z+px)/4})$ is regular in $\ve_1,\ve_2$. Therefore we can take its exponential,
which is still regular in $\ve_1,\ve_2$.

Put $a:=t/2$. Then  we get by 
\cite[Prop.~2.3]{NY2}  and \cite[formula (2.12)]{NY2}  that
\begin{equation}\label{uT}
\frac{\partial \F_0}{\partial \log \Lambda}= -4u, \quad 
\frac1{32}\frac{\partial^2 \F_0}{(\partial \log \Lambda)^2}=\frac{1}{24}
E_2 \Big(\frac{du}{da}\Big)^2-\frac{u}{6}=T.
\end{equation}
As $X$ is a rational surface, we have $\chi=-\sigma+4$, thus we get by (\ref{AB})
\begin{equation}\label{AB2}
\exp(\chi A+ \sigma B)
=\Big(\frac{\sqrt{-1}}{\Lambda}\frac{du}{da}\Big)^2\theta_{01}^\sigma.
\end{equation}
These relations and \eqref{ududaa} hold in a neighborhood of
$(a,\Lambda)=(\infty,0)$.
Thus they are equalities in ${\Bbb C}[[\Lambda/a]][\Lambda,a/\Lambda]$.
Note that the $q$-development of
\begin{equation}
  \label{eq:q1/8}
\frac{\Lambda}{a}=\left(
\sqrt{-1}\frac{2 E_2+\theta_{00}^4+\theta_{10}^4}{3\theta_{00} \theta_{10}}
\right)^{-1}
\end{equation}
starts with $\frac{2}{\sqrt{-1}}q^{\frac{1}{8}}$. Thus $\C[[q^{1/8}]]
\cong \C[[\Lambda/a]]$. 
Putting (\ref{Fdev})--(\ref{AB2}) 
into  Corollary \ref{deltaq} we get
\begin{NB}
for $(\frac{t}{2},\Lambda)\in U$ : commented out by HN (Apr. 8)
\end{NB}
\begin{equation}\label{modelt}
\begin{split}
\delta_{\xi,t}^X&(\exp(\alpha z+px))=\frac{1}{\Lambda}
\lim_{\ve_1,\ve_2\to 0}\exp\Big(\sum_{i=1}^\chi \frac{1}{w(x_i)w(y_i)}
F\big(w(x_i),w(y_i),\hbox{$\frac{t-\iota_{p_i}^*\xi}{2}$},\Lambda e^{\iota_{p_i}^*(\alpha z+px)/4}\big)\Big)\\
&=\lim_{\ve_1,\ve_2\to 0}\exp\Big(\sum_{i=1}^\chi
\Big(
\frac{\F_0(\hbox{$\frac{t-\iota_{p_i}^*\xi}{2}$};\Lambda e^{\iota_{p_i}^*(\alpha z+px)/4})}{w(x_i)w(y_i)}
+\frac{w(x_i)+w(y_i)}{w(x_i)w(y_i)}H(\hbox{$\frac{t-\iota_{p_i}^*\xi}{2}$};\Lambda e^{\iota_{p_i}^*(\alpha z+px)/4})\\&\qquad
+A(\hbox{$\frac{t-\iota_{p_i}^*\xi}{2}$};\Lambda e^{\iota_{p_i}^*(\alpha z+px)/4})+\frac{w(x_i)^2+w(y_i)^2}{3w(x_i)w(y_i)}
B(\hbox{$\frac{t-\iota_{p_i}^*\xi}{2}$};\Lambda e^{\iota_{p_i}^*(\alpha z+px)/4})\Big)\Big)\\
&=\sqrt{-1}^{\<\xi,K_X\>}\Big(\frac{1}{\Lambda}
q^{-\frac{1}{2}(\frac{\xi}{2})^2} \exp\Big(\frac{du}{da} \<\alpha,\xi/2\> z+T \<\alpha^2\> z^2 - ux\Big)
\Big(\frac{\sqrt{-1}}{\Lambda}\frac{du}{da}\Big)^2\theta_{01}^\sigma\Big).
\end{split}
\end{equation}
The final equality, i.e.\ the first term $=$ the last term,   
holds in  $\C((\Lambda))((\Lambda/a))[[z,x]]=
\C((\Lambda))((1/t))[[z,x]]$.  Indeed for the left-hand side 
 the
coefficient of $z^n x^m$ is in $\Lambda^{-\xi^2-3}\linebreak[2]\times t^{\xi^2+2+n+2m}
\C[[\Lambda,1/t]]$.  Thus we get
\begin{align*}
&\delta_{\xi}^X(\exp(\alpha z+px))\\
&\ =-\sqrt{-1}^{\<\xi,K_X\>}
\res_{t=\infty}\Big(\frac{1}{\Lambda}
q^{-\frac{1}{2}(\frac{\xi}{2})^2} \exp\Big(\frac{du}{da} \<\alpha,\xi/2\> z+T \<\alpha^2\> z^2 - ux\Big)
\Big(\frac{\sqrt{-1}}{\Lambda}\frac{du}{da}\Big)^2\theta_{01}^\sigma dt\Big).
\end{align*}

Finally we want to express this result in terms of the $q$ development 
of the modular forms involved.
That is, we change the variable from $t$ to $q$.
First we determine $\frac{da}{d\tau}$.
Combining formulas (V.4.1), (V.5.2) and (V.5.6) of \cite{Ch}
(note that in the notation of \cite{Ch}  $\theta_{00}(\tau)=\theta_3(0,\tau)$, 
$\theta_{01}(\tau)=\theta_{2}(0,\tau)$,
 $\theta_{10}(\tau)=\theta_1(0,\tau)$), 
we get 
$$\frac{d\log(\theta_{00})}{d \tau}-
\frac{d\log(\theta_{10})}{d\tau}
=\frac{\pi}{4\sqrt{-1}}\theta_{01}^4.$$
By \cite[(VII.3.10)]{Ch} we have 
$\theta_{00}^4-\theta_{10}^4=\theta_{01}^4$, and thus 
\begin{align*}
\frac{d u}{d \tau}&=-\frac{d}{d\tau}\Big(\frac{\theta_{00}^2}{\theta_{10}^2}+
\frac{\theta_{10}^2}{\theta_{00}^2}\Big)\Lambda^2=-2\Lambda^2\frac{\theta_{00}^4-\theta_{10}^4}{\theta_{10}^2\theta_{00}^2}
\Big(\frac{d \log(\theta_{00})}
{d\tau}-\frac{d \log(\theta_{10})}{d\tau}\Big)
\\&
=-\frac{\Lambda^2\pi}{2 \sqrt{-1}} \frac{\theta_{01}^8}{\theta_{10}^2
\theta_{00}^2}=-\frac{\pi\sqrt{-1}}{8}\Big(\frac{du}{da}\Big)^2\theta_{01}^8.
\end{align*}
\begin{NB} Should this have a number, so that we can cite it in the $K$-theory paper 6.9. LG
\end{NB}
Thus we get
$$\frac{da}{d\tau}=\frac{da}{du}\frac{du}{d\tau}
=-\frac{\pi\sqrt{-1}}{8}\frac{du}{da}\theta_{01}^8.$$
By $a=t/2$ and $q=e^{2\pi \sqrt{-1}\tau}$, we have 
$$dt=2da=2\frac{da}{d\tau}d\tau=\frac{1}{\pi \sqrt{-1}}\frac{da}{d\tau}\frac{dq}{q}=
-\frac{1}{8}\frac{du}{da}\theta_{01}^8\frac{dq}{q}.$$
\begin{NB}
Added by K.Y. Apr. 7:
Although we do not need to mention,
I would like remark that the equality holds at $a=\infty$,
since both are defined at $a=\infty$.
\end{NB}
By \eqref{eq:q1/8}
the residue at 
$a=\infty$ is $8$ times the residue at $q=0$. 
Therefore we get 
\begin{equation*}
\begin{split}
\delta^X_{\xi}(\exp&(\alpha z +px))=\\
& \sqrt{-1}^{\<\xi,K_X\>-1}\res_{q=0}\Bigg[
q^{-\frac{1}{2}(\frac{\xi}{2})^2} \exp\Big(\frac{du}{da} \<\alpha,\xi/2\> z+T \<\alpha^2\>z^2 - ux\Big)
\Big(\frac{\sqrt{-1}}{\Lambda}\frac{du}{da}\Big)^3\theta_{01}^{\sigma+8}\frac{dq}{q}\Bigg],\end{split}
\end{equation*}
and Theorem \ref{main2} follows by $\sigma+8=3\sigma+2\chi=K_X^2$.
\end{pf*}


\begin{Remark}
(1) Denote by $u_{MW}$ and $h$ the functions denoted by $u$,  $h$ in \cite{MW}.
Note that in the notation of \cite{MW}, $\lambda_0=c_1/2$ and
$\lambda=\xi/2$. We are computing the wallcrossing for $\Phi^X_{c_1}$, whereas
in \cite{MW} the wallcrossing for $D^X_{c_1}$ is computed. Thus we have to 
multiply their formula (4.6) by $(-1)^{-(c_1^2+c_1K_X)/2}$ to compare
it with ours.
Write $\delta^X_{\xi,MW}$ for the wallcrossing formula obtained this way.
Using the fact that $u=-2u_{MW}$, $\frac{du}{da}=\frac{\sqrt{-1}\Lambda}{h}$,
we see that $\delta^X_{\xi,MW}=-\frac{1}{2}\delta^X_{-\xi}$. By definition $\delta^X_{-\xi}=-\delta^X_{\xi}$. Thus $\delta_{\xi,MW}^X=\frac{1}{2}\delta^X_\xi$.
It was observed in \cite{MW} that the formula in \cite{G} gives 
${2}\delta^X_{\xi,MW}$ for the wallcrossing of $D_{c_1}^X$. Thus our formula agrees with the  results of \cite{G}.

(2) 
Denote by $U$, $f$, $R$ the functions denoted by the same letters in 
\cite{GZ}.
Then it is easy to check that 
$$U(\tau)=-\frac{1}{\Lambda^2}u(\tau+1), \quad 
\frac{1}{f(\tau)}=\frac{1}{\Lambda}\frac{du}{da}(\tau+1),\quad
R(\tau)=- \Big(\frac{1}{\Lambda}\frac{du}{da}(\tau+1)\Big)^4
\theta_{01}(\tau+1)^8.$$
Using these formulas it  is also easy to see directly that Theorem \ref{main2} gives the
same wallcrossing formula as \cite{G},\cite{GZ}, after correcting for the different sign conventions.
\end{Remark}
\begin{NB} Detailed check of the signs.
Note that 
the formula of \cite{MW} has as sign
$-\frac{\sqrt{-1}}{2} (-1)^{(\xi-c_1)K_X/2+c_1^2/2}$, thus we have to look at
$$-\frac{\sqrt{-1}}{2} (-1)^{(\xi-c_1)K_X/2+c_1^2/2-(c_1^2+c_1K_X)/2}=
-\frac{\sqrt{-1}}{2}(-1)^{\xi K_X/2}(-1)^{-c_1K_X}
$$
By $$(-1)^{-\xi K_X/2}=(-1)^{-\xi K_X}(-1)^{\xi K_X/2}=(-1)^{c_1K_X}(-1)^{\xi K_X/2},$$
this is $-\frac{\sqrt{-1}}{2}(-1)^{\xi K_X/2}$, which is precisely the negative of the sign that our formula gives for $-\xi$. The rest of the formula of \cite{MW} coincides with ours
for $\xi$, except that we replace $-\frac{du}{da}$ by $\frac{du}{da}$, which gives precisely the correct thing if we replace $\xi$ by $-\xi$. (L.G. 3.3)

\end{NB}

\begin{NB}
I insert Lothar's message. Apr.30, HN

I wanted to explain the sign error in our paper that I found and 
also how (I think) it can be fixed.

First I explain the sign error: 
Look at $\P^2$. Then on $M_{\P^2}(H,d)$, the class $2\mu(H)$ is nef and 
big (in fact it is ample on the Uhlenbeck compactification). Thus we get 
 $\Phi_H^{\P^2}(H^{4n})>0$ for all $n\in \Z_{\ge 0}$.

Similarly on $M_{\P^2}(0,d)$, the class $\mu(H)$ is nef and big, thus 
 $\Phi_0^{\P^2}(H^{4n+1})>0$ for all $n\in \Z_{>0}$ (this is not true for 
$n=0$, because then we are not in the stable range, we cannot avoid the strictly semistable sheaves).

Now we can compute these using our wallcrossing formula:

{\bf Case $c_1=H$:} Let $\widehat \P^2$ be the blowup of $\P^2$ in a point, 
$E$ the exceptional divisor. We know that the moduli space 
$M_{\widehat P^2}^L(c_1,d)$ is empty for $L$ near to $F:=H-E$, and by the 
blowup formula we know that  $\Phi_H^{\P^2}(H^{4n})=
\Phi_H^{\widehat \P^2,H}(H^{4n})$. Thus we have to sum over all walls $\xi$
of type $H$ with $\xi H>0>\xi (H-E)$. 
Thus the walls are the $\xi=(2n-1)H-2a E$ with $a\ge n >0$.
The sign in Thm 4.2 is $\sqrt{-1}^{-\xi K_X -1}$.
By $-K=3H-E$ this is $\sqrt{-1}^{6n-3-2a-1}=(-1)^{n+a}$.
Putting this together we get
$$\Phi_{H}(\exp(Hz))=\sum_{a\ge n>0}(-1)^{n+a} \Big[q^{(a^2-(n-1/2)^2)/2} 
\exp((n-1/2) \frac{du}{da} z +Tz^2)(\sqrt{-1} \frac{du}{da})^3 \theta_{01}^8.
$$
One checks on the computer that this gives precisely the negative of the 
correct value.  

{\bf Case $c_1=0$:} In this case the blowup formula gives that 
$\Phi_0^{\P^2}(H^{4n+1})=-
\Phi_E^{\widehat \P^2,H}(EH^{4n+1})$. 
(The blowup formula for the Donaldson invariants would give the same
formula with a positive sign. But to get the formula for $\Phi$ one has to
multiply both sides with $(-1)^{(c_1^2+c_1K)/2}$, which gives a sign difference
on by $(-1)$. The corresponding computation to the case $c_1=H$ gives that
one gets precisely $(-1)$ times the correct result for  
$\Phi_0^{\P^2}(H^{4n+1})$.

In fact it is much easier to see that the sign is wrong (and this I use to find
where the error is). Lets look at the simplest case where a wallcrossing can
occur. We take $\widehat \P^2$, $H$ as first Chern class and $\xi=H-2E$.
We look at the wallcrossing for $d=0$.
Then we see that for $L\xi>0$, $M_{\widehat \P^2}^L(H,0)=pt$ and 
for $L\xi <0$, $M_{\widehat \P^2}^L(H,0)=\emptyset$.
Thus we get $\delta_{\xi}^{\widehat \P2}(1)=1$. 

Theorem 4.2 again gives  $\delta_{\xi}^{\widehat \P2}(1)=-1$
i.e. 
$$\Big[q^{3/8}\Big(\frac{\sqrt{-1}}{\Lambda} \frac{du}{da}\Big)^3 \theta_{01}^8\Big]_{q^0}=
(-1)^3=(-1).$$

The advantage is that for all formulas in the paper we can easily check,
whether they give the correct result in this example. 
Note that in this case $l=0$, $\AA_+=\oo_{pt}$, $\AA_-=0$.

We get the following:
\begin{enumerate}
\item Definition 2.5 correctly gives in this case that the coefficient of $\Lambda^0$ of 
$\delta_{\xi,t}(1)$ is $\frac{1}{t}$, and thus $\delta_{\xi}(1)=1$.
\item Also Lemma 3.9 correctly gives that the coefficient of $\Lambda^0$ of 
$\delta_{\xi,t}(1)$ is $1$.
\item On the other hand already formula (4.11) gives the wrong sign:
Formula (4.11) gives 
$$\delta_{\xi,t}=\sqrt{-1}^{-\<\xi K_X\>} \frac{q^{3/8}}{\Lambda}
\Big(\frac{\sqrt{-1}}{\Lambda}\frac{du}{da}\Big)^2.$$ As $-\xi K=1$,
This is congruent modulo $q^{1/2}$ to 
$\sqrt{-1} \frac{q^{3/8}}{\Lambda} (-q^{-1/8})^2=\frac{\sqrt{-1}q^{1/8}}{\Lambda}$. 
On the other hand $t=2a\equiv \sqrt{-1} q^{-1/8}\Lambda$. Thus we get 
 $\delta_{\xi,t}(1)=-\frac{1}{t}+O(t^{-2})$. 
\end{enumerate}

Thus the sign error must occur between these two places.
Not very much is happening in between, basically we only put in the 
perturbation part. I think the error lies there, that we did not handle
the perturbation part correctly:
In fact I think that $H(a,\Lambda)$ is wrong.

I think that the formula (on the bottom of page 27) after formula (4.3)
(Fdev)  is incorrect. The formula there is $H(a,\Lambda)=-\pi \sqrt{-1} a$.
We refer this to [NY2] section 5.3; there the statement is $H(a,\Lambda)=-\log(-1)a$. The first thing is that this sign is not evident from the formula (E.5)
in [NY2]. It will depend on the choice of the branch of $\log$.

On the  other hand in the K-theoretic perturbation part you have 
$H(a,\Lambda)=\pi \sqrt{-1} a$ (Section 4.3 of [NY3]).
Note that the K-theory perturbation part should tend to the topological 
perturbation part as $\beta\to 0$ (and $H$ does not depend on $\beta$). Thus there is an incompatibility. I expect that this is due to the difference in convention for the branch of $\log$ (although I have not understood it). In any case we have computed the contribution of 
the perturbation part via the $K$-theoretic version, thus we have to use the 
$$H(a,\Lambda)=\pi \sqrt{-1} a$$.

Thus in formula 
[H] (3.6) er get $\pi \sqrt{-1} \<xi/2,K_X\>$.

This changes the sign at RHS of formula [nodelt] (4.11) to 
$\sqrt{-1}^{\<\xi,K_X\>}$, and in all the other formulas 
$\sqrt{-1}^{-\<\xi,K_X\>}$ is replaced by $\sqrt{-1}^{\<\xi,K_X\>}$.
In particular the sign on the RHS of Theorem 4.2 is 
$\sqrt{-1}^{\<\xi,K_X\>-1}$.
Thus the sign has changed by $(-1)^{\<xi,K_X\>}= (-1)^{\xi^2}$. 
This will give the correct results.

Surprisingly remark 4.12 (comparing with Moore-Witten) is still true after making this correction, because
I also found a sign error in the proof of Remark 4.12. 
This is just an instance of the fact that it is easy to get whatever sign one
wants:

When comparing with the sign in Moore-Witten's paper (this is formula (4.6))
I made the following mistake. I used only the sign that is written in front
of the formula to compare, but there is also a hidden sign: 
Moore-Witten have $\frac{1}{h^3}$ instead of our 
$\big(\frac{\sqrt{-1}}{\Lambda}
\frac{du}{da}\big)^3$. However by definition 
$\frac{du}{da}=\frac{\sqrt{-1}\Lambda}{h}$, thus 
$\frac{1}{h^3}=(-1)^3\big(\frac{\sqrt{-1}}{\Lambda}
\frac{du}{da}\big)^3$. This gives an extra sign, and my statement 
is still correct.
\end{NB}




\section{Generalization to non-toric surfaces}\label{cobor}

In this section we show that the wallcrossing term is given by the
same formula as in \thmref{main2} for a good wall of
an {\it arbitrary\/} simply connected projective surface $X$. The
proof is based on \cite{EGL} for Chern numbers of Hilbert schemes of
points. 

We consider the Grothendieck group $K(Y)$ of locally free sheaves on a
smooth projective variety $Y$.
It is isomorphic to that of coherent sheaves.
It has a ring structure from the tensor product. We denote it by
$\otimes$.
For a morphism $f\colon Y_1\to Y_2$ we have a pushforward homomorphism
$f_!\colon K(Y_1)\to K(Y_2)$, and the pullback homomorphism $f^!\colon
K(Y_2)\to K(Y_1)$. We also have the involution $\vee$ on $K(Y)$ given
by the dual vector bundle for a vector bundle.

Let $X$ be a projective surface and $X^{[n]}$ denote the Hilbert
scheme of $n$ points on $X$.
As before let $X_2 = X \sqcup X$ be the disjoint union of two copies of $X$. Let
\(
  X_2^{[l]}
\)
be the Hilbert scheme of $l$ points on $X_2$, i.e.\ 
\(
  X_2^{[l]} = \bigsqcup_{m+n=l} X^{[m]}\times X^{[n]},
\)
and let $\mathcal I_1$ (resp.\ $\mathcal I_2$) be the sheaf on
$X\times X_2^{[l]}$ whose restriction to $X\times X^{[m]}\times
X^{[n]}$ is $p_{12}^*(\mathcal I_{Z_{m}(X)})$ (resp.\ 
$p_{13}^*(\mathcal I_{Z_{n}(X)})$).
Let us define $p$, $q$ by
\begin{equation*}
   p\colon X_2^{[l]}\times X \to X_2^{[l]}, \qquad
   q\colon X_2^{[l]}\times X \to X.
\end{equation*}
These maps depend on $l$, but we suppress the dependence from the
notation hoping that they do not lead to confusion, though we will vary
$l$ later.

In this section we prove the following:
\begin{NB}
\begin{Theorem}
There exist universal power series
$A_i\in \Q((t^{-1}))[[\Lambda,(\tau_k^\rho)]]$, $i=1,..,7$, depending
on *** such that
\begin{equation*}
\begin{split}
   & 
   (-1)^{\chi(\mathcal O_X)+\xi(\xi-K_X)/2} t^{-\xi^2-2\chi(\mathcal O_X)}
   \Lambda^{\xi^2+3\chi(\mathcal O_X)}
   \delta^X_{\xi,t}
   \left(\exp\left(\sum_\rho \alpha_\rho\right)\right)
\\
   =\; & \exp( \xi^2 A_1 + \xi\cdot c_1(X) A_2
   + c_1(X)^2 A_3 + c_2(X) A_4 + \alpha\cdot \xi A_5
   + \alpha\cdot c_1(X) A_6 + \alpha^2 A_7).
\end{split}
\end{equation*}
\end{Theorem}
\end{NB}
\begin{NB}
The above is not quite correct. We need to include various
intersection numbers of $b_k^\rho$'s, etc. 

What I meant by putting this into NB is that I wanted to not keep it. 
It is very complicated to formulate, the proof becomes more complicated and we
are not using it. 11.4. LG

Well, I would like to keep it in \verb+NB+ for a possible use in a future.  
Apr. 17, HN
\end{NB}

\begin{Theorem}\label{univ}
There exist universal power series
$A_i\in \Q((t^{-1}))[[\Lambda]]$, $i=1,..,8$, such that for all projective surfaces $X$ and all $\xi\in \Pic(X)$
\begin{equation*}
\begin{split}
   & 
   (-1)^{\chi(\mathcal O_X)+\xi(\xi-K_X)/2} t^{-\xi^2-2\chi(\mathcal O_X)}
   \Lambda^{\xi^2+3\chi(\mathcal O_X)}
   \delta^X_{\xi,t}
   \left(\exp\left(\alpha z+px\right)\right)
\\
   =\; & 
   \begin{aligned}[t]
   \exp( \xi^2 A_1 + \xi\cdot c_1(X) A_2
   + c_1(X)^2 A_3 + c_2(X) A_4 + \alpha\cdot \xi A_5 z \qquad\qquad \\
   + \alpha\cdot c_1(X) A_6 z + \alpha^2 A_7 z^2+x A_8).
   \end{aligned}
\end{split}
\end{equation*}
\end{Theorem}

Here 
$\delta^X_{\xi,t}$ is defined for arbitrary projective surface by the
same formula \eqref{wallct} except that we change $\Lambda^{4l-\xi^2-3}$ into
$\Lambda^{4l-\xi^2-3\chi(\mathcal O_X)}$ and also $\mathcal
A_{\xi,-}$, $\mathcal A_{\xi,+}$ into
\begin{equation*}
   -p_! (\mathcal I_2^\vee\otimes \mathcal I_1\otimes q^!\xi),
\qquad
   -p_! (\mathcal I_1^\vee\otimes \mathcal I_2\otimes q^!\xi^\vee)
   \in K(X_2^{[l]})
\end{equation*}
respectively.

\begin{NB}
We can consider $\xi$ as an element of either $\operatorname{Pic}(X)$
or $\operatorname{NS}(X)$. Here is a record of discussion:

In order to consider $\xi\in H^2(X,\Z)$ as a class in $K$-theory, we
need to use the topological $K$-theory. But for our purpose, it is
probably more natural to assume $\xi\in\operatorname{NS}(X)$
(N\'eron-Severi group). Apr. 9 HN.

Yes: $\xi$ should not just be any class in $H^2(X,\Z)$, but e.g. an element of 
$NS(X)$. I actually think it is best to just take $\xi$ an element of $Pic(X)$, because 
this is the way we use it. Because of this I have changed it to $Pic(X)$ above. 
11.4 LG.

I am not so sure any more that $Pic(X)$ is better. If you prefer it
can be changed back to $NS(X)$ 21.4 LG

It seems that a map $\Pic(X)\to K(X)$ does not descend to $NS(X)\to
K(X)$. If this is true, we need to use $\Pic(X)$. Apr. 24, HN.
\end{NB}

When $\xi$ is good, both
$\Ext_p^0(\mathcal I_2,\mathcal I_1(\xi))$,
$\Ext_p^2(\mathcal I_2,\mathcal I_1(\xi))$ vanish
\cite[Lemma 4.3]{EG1}. Therefore we have
\begin{equation*}
 \Ext_p^1(\mathcal I_2,\mathcal I_1(\xi)) 
 = -p_! (\mathcal I_2^\vee\otimes \mathcal I_1\otimes q^!\xi) 
\end{equation*}
and the same for $\Ext_p^1(\mathcal I_1,\mathcal I_2(-\xi))$.

\begin{NB}
In fact, $\Ext_p^\bullet$ are the higher derived functors of the
composite functor $p_*\circ{\mathcal H}om$. The induced homomorphism
in the Grothendieck group
\[
   \sum_{i=0}^2 (-1)^i \Ext_p^i 
   = \sum_{i,j} (-1)^i R^i p_* (-1)^j {\mathcal E}xt^j
\]
is the composite of the corresponding homomorphisms $p_!$ and $\sum
(-1)^j{\mathcal E}xt^j(\bullet,\bullet) = \bullet^\vee \otimes
\bullet$.
\end{NB}

\begin{NB}
By Riemann-Roch we have
\begin{equation*}
\begin{split}
   & \chi(L) =
   \int_X e^{\xi}\left(1 + \frac12 c_1(X) +
   \operatorname{Todd}_2(X)\right)
  = \chi(\mathcal O_X) + \frac{\xi(\xi-K_X)}2,
\\
 & \chi(L^\vee) =
   \int_X e^{-\xi}\left(1 + \frac12 c_1(X) +
   \operatorname{Todd}_2(X)\right)
  = \chi(\mathcal O_X) + \frac{\xi(\xi+K_X)}2.
\end{split}
\end{equation*}
\end{NB}

The proof is a straightforward modification of that of
\cite[Th.~4.2]{EGL}, so we only give a sketch of the proof. The
essential point is to use the incidence variety to compute the
intersection products on Hilbert schemes recursively. A slight
difference is that we need to introduce {\it two\/} incidence
varieties because we study Hilbert schemes of a nonconnected surface
$X_2$.

For $\alpha=1,2$ let $X^{[l,l+1]}_{2,\alpha}$ be the variety of
pairs $Z$, $Z'$ in $X_2^{[l]}\times X_2^{[l+1]}$ satisfying $Z\subset
Z'$ and $Z'\setminus Z$ is a point in the
$\alpha^{\mathrm{th}}$-factor of $X_2$.  This is an obvious
generalization of the incidence variety $X^{[l,l+1]}$, studied by
various people and used in \cite{EGL}.
Let $\phi_\alpha$ and $\psi_\alpha$ be the projections from
$X_{2,\alpha}^{[l,l+1]}$ to $X_2^{[l]}$ and $X_2^{[l+1]}$ respectively.
Let $\rho_\alpha$ be the map $X_{2,\alpha}^{[l,l+1]}\to X$ defined by letting
$\rho(Z,Z')$ be the unique point in $Z'\setminus Z$. Let $\mathcal L$
be the line bundle whose fiber at $(Z,Z')$ is the kernel of the
homomorphism  $H^0(\mathcal O_{Z'}) \to H^0(\mathcal O_Z)$.
We have
\begin{equation*}
\begin{CD}
   X @<{\rho_\alpha}<< X_{2,\alpha}^{[l,l+1]} @>\psi_\alpha>> X_2^{[l+1]}
\\
   @. @VV{\phi_\alpha}V @.
\\
   @. X_2^{[l]} @.
\end{CD}
\end{equation*}
We also define 
\(
   j_\alpha = \rho_\alpha\times\operatorname{id}\colon
   X_{2,\alpha}^{[l,l+1]}\to X\times X_{2,\alpha}^{[l,l+1]}
\)
and
\(
  \sigma_\alpha = \rho_\alpha\times\phi_\alpha \colon 
X_{2,\alpha}^{[l,l+1]}\to X\times X_2^{[l]}.
\)

We first have the following analog of [loc.\ cit., (5)]
\begin{equation}\label{eq:ideal}
   \psi_{\alpha X}^! \mathcal I_\beta
   = \phi_{\alpha X}^! \mathcal I_\beta 
   - \delta_{\alpha\beta} j_{\alpha !} \mathcal L
   = \phi_{\alpha X}^! \mathcal I_\beta 
   - \delta_{\alpha\beta}p^!\mathcal L\otimes
   \rho_{\alpha X}^!\mathcal O_\Delta,
   \qquad \text{for $\alpha,\beta=1,2$},
\end{equation}
where $p\colon X\times X_{2,\alpha}^{[l,l+1]}\to
X_{2,\alpha}^{[l,l+1]}$ and $f_X = f\times \id_X$ for $f=\phi_\alpha$,
$\psi_\alpha$.

Next we have an analog of [loc.\ cit., (8)]
\begin{equation}
  \psi_\alpha^* \operatorname{ch}(\mathcal I_\beta)/c
  = \phi_\alpha^* \operatorname{ch}(\mathcal I_\beta)/c
  - \delta_{\alpha\beta} \operatorname{ch}(\mathcal L)
  \cdot \rho_\alpha^* c
\end{equation}
for $c\in H_*(X)$.

We also get an analog of \cite[Prop.\ 2.3]{EGL} using \eqref{eq:ideal}
\begin{equation}
\begin{split}
   & \psi_\alpha^!p_!(\mathcal I_2^\vee\otimes \mathcal I_1\otimes q^!\xi)
\\
   =\; &
   \phi_\alpha^!p_!(\mathcal I_2^\vee\otimes \mathcal I_1\otimes q^!\xi)
   -\delta_{\alpha 1}\sigma_\alpha^!\mathcal I_2^\vee
   \otimes \rho_\alpha^!\xi \otimes \mathcal L
   -\delta_{\alpha 2}\sigma_\alpha^!\mathcal I_1
   \otimes \rho_\alpha^!(\xi
   \otimes\omega_X^\vee)\otimes \mathcal L^\vee.
\end{split}
\end{equation}
More precisely, we do not get a term corresponding to the third term in [loc.\ cit., (10)] coming
from the product of two copies of the diagonal, because 
$\delta_{\alpha 1}\delta_{\alpha 2}$ is always $0$.

\begin{NB}
\begin{proof}
\begin{equation*}
\begin{split}
   & \psi_\alpha^!p_!(\mathcal I_2^\vee\otimes \mathcal I_1\otimes
   q^!\xi)
   = p_! \psi_{\alpha X}^! (\mathcal I_2^\vee\otimes  \mathcal I_1\otimes
   q^!\xi)
\\
   =\; &
   p_! \left((\phi_{\alpha X}^! \mathcal I_2 - 
     \delta_{\alpha 2} p^!\mathcal L\otimes
   \rho_{\alpha X}^!\mathcal O_{\Delta})^\vee\otimes
   (\phi_{\alpha X}^! \mathcal I_1 - 
     \delta_{\alpha 1} p^!\mathcal L\otimes
   \rho_{\alpha X}^!\mathcal O_{\Delta}) \otimes q^!\xi\right)
\\
   =\; &
   \phi_\alpha^! p_!(\mathcal I_2^\vee\otimes \mathcal I_1\otimes
   q^!\xi)
   - \delta_{\alpha 1} p_! (\rho_{\alpha X}^!\mathcal O_{\Delta}^\vee
   \otimes \phi_{\alpha X}^!\mathcal I_1 \otimes q^! \xi)\otimes 
   \mathcal L^\vee
   - \delta_{\alpha 2} p_! (\rho_{\alpha X}^!\mathcal O_{\Delta}
   \otimes \phi_{\alpha X}^!\mathcal I_2^\vee \otimes q^! \xi)\otimes 
   \mathcal L.
\end{split}
\end{equation*}
And
\begin{equation*}
   \rho_{\alpha X}^!\mathcal O_\Delta \otimes 
   \phi_{\alpha X}^! \mathcal I_2^\vee
   = j_{\alpha !}j_\alpha^!\phi_{\alpha X}^! \mathcal I_2^\vee
   = j_{\alpha !} \sigma_\alpha^! \mathcal I_2^\vee,
\end{equation*}
where $j_\alpha = \rho_\alpha\times \id\colon X_{2,\alpha}^{[n,n+1]}
\to X\times X_{2,\alpha}^{[n,n+1]}$.
\end{proof}
\end{NB}



Using these results,  the same argument as in \cite[Prop.~3.1, Thm~4.1]{EGL} 
shows the following.

\begin{Lemma}\label{pol} Fix $l\ge 0$. 
Let $P$ be any polynomial in the $c_{i_1}(\AA_+)$, $c_{i_2}(\AA_-)$, 
$\ch_{i_3}(\I_1) \xi^{i_4}/(\alpha z +p x)$, $\ch_{i_5}(\I_2) \xi^{i_6}/(\alpha z +p x)$ for $i_1,\ldots,i_6\in 
\Z_{\ge 0}$,
then there exists a universal polynomial $Q$ 
(depending only on $P$,)
in $\xi^2$, $\xi c_1(X)$, $c_1(X)^2$, 
$c_2(X)$, $\xi \alpha z$, $\alpha c_1(X) z$, $\alpha^2 z^2, x$, such that 
$\int_{X_{2}^{[l]}}P=Q$.
\end{Lemma}

We denote the left-hand-side of Theorem \ref{univ} by 
$\overline \delta_{\xi,t}^X$. By definition we have 
$$(-1)^{rk(\AA_-)} t^{rk(\AA_+)+rk(\AA_-)}\frac{1}{c^t(\AA_+)c^{-t}(\AA_-)}=
\sum_{i,j}s_i(\AA_-)s_{j}(\AA_+)(-1)^i t^{-i-j},
$$ where
$rk(\AA_-)=l-\chi(\oo_X)-\frac{\xi(\xi-K_X)}2$, 
 $rk(\AA_+)=l-\chi(\oo_X)-\frac{\xi(\xi+K_X)}2$.
Therefore by Lemma \ref{pol}
we can write $\overline \delta_{\xi,t}^X(\alpha z +p x)=
\sum_{l\ge 0} \sum_{i\in \Z} \Lambda^{4l} P_{l,i} t^i$, where 
$P_{l,i}$ is a universal polynomial in  $\xi^2$, $\xi c_1(X)$, $c_1(X)^2$, $c_2(X)$,  $\xi \alpha z$, $\alpha c_1(X) z$, $\alpha^2 z^2$, $x$, depending only 
on $l$ and $i$.  It is easy to see from the definition that the coefficient of $\Lambda^0$ of $\overline \delta_{\xi,t}^X$ as a power series in $\Lambda$ is $1$. Thus there is a universal power series 
$G_{t,\Lambda}\in \Q((t^{-1}))[x_1,\ldots,x_8][[\Lambda]]$, such that
$\overline \delta_{\xi,t}^X(\alpha z +px)=\exp(G_{t,\Lambda}(\xi^2,\xi c_1(X),c_1(X)^2,c_2(X),\xi \alpha z, \alpha c_1(X) z, \alpha^2 z^2, x))$.

Now assume that $X=Y\sqcup Z$ for $Y,Z$ not necessarily connected projective surfaces, and $\xi\in \Pic(X)$, $\beta\in H_2(X)z\oplus H_0(X)x$ 
satisfy $\xi|_Y=\xi_1$, $\xi|_Z=\xi_2$,
$\beta|_Y=\beta_1$, $\beta|_Z=\beta_2$. Then 
$X^{[l]}_2=\coprod_{n+m=l} Y^{[n]}_2\times Z^{[m]}_2$, and denoting 
$\AA_{-,X}$, $\AA_{-,Y}$, $\AA_{-,Z}$ respectively, the bundles $\AA_-$ on $X^{[l]}_2$, $Y^{[n]}_2$ and $Z^{[m]}_2$, it is obvious that $\AA_{-,X}|_{Y^{[n]}_2\times Z^{[m]}_2}=\AA_{-,Y}\boxplus\AA_{-,Z}$ and similarly for $\AA_+$, $\I_1$, $\I_2$.
Thus it follows from the definitions that 
\begin{equation}\label{prod}
\overline\delta^X_{\xi,t}(\beta)=\overline\delta^Y_{\xi_1,t}(\beta_1)\overline\delta^Z_{\xi_2,t}(\beta_2).
\end{equation}
To a triple $(X,\xi,\beta)$ of a projective surface $X$, a class $\xi\in \Pic(X)$ 
and  $\beta\in H_2(X)z\oplus H_0(X)x$ we associate the vector
$v(X,\xi,\beta) :=(\xi^2,\xi c_1(X),c_1(X)^2,c_2(X),\xi \beta, \beta c_1(X) , \beta^2, \beta)\linebreak[2]\in \Z^8$, where we suppress $\int_X$ in the notation. Then we know 
$\overline\delta^X_{\xi,t}(\beta)=\exp(G_{t,\Lambda}(v(X,\xi,\beta))$.
Choose
triples $(X_i,\xi_i,\beta_i),\  i=1,\ldots, 8$ as above such that 
the $w_i:=v(X_i,\xi_i,\beta_i)$ form a basis of $\Q^8$.
Let $(a_{i,j})_{i,j=1}^8$ be the matrix such that $\sum_{j} a_{i,j}w_j=e_i$ for all $i$,
where $e_i$ is the vector with $i$-th entry $1$ and all others zero.
For all $i$ put $A_i:=\sum_j a_{i,j}G_{t,\Lambda}(w_j)$. Let $(X,\xi,\beta)$ be a triple, such
that $(v^1,\ldots,v^8):=v(X,\xi,\beta)=\sum_{i=1}^8 n_i w_i$, with $n_i\in \Z_{\ge 0}$. 

Then by \eqref{prod} we get
$\overline\delta_{\xi,t}^X(\beta)=\exp(\sum_j n_j G_{t,\Lambda}(w_j))$. Thus by  
$\sum_i v^i a_{i,j}=n_j$ for all $j$, 
we get
$\overline\delta_{\xi,t}^X(\beta)=\exp\big(\sum_{i} v^i A_i\big)$.
\begin{NB} corrected 6.9. LG\end{NB} 
Note that the $v^i$ are just the intersection numbers $\xi^2,\ldots,\beta$.
As the set of all vectors $\sum_{i=1}^8 n_i w_i$  with all $n_i\in \Z_{\ge 0}$  is Zariski dense in 
$\Q^8$, the last equality holds for all triples $(X,\xi,\beta)$
of a projective surface $X$, a class $\xi\in \operatorname{Pic}(X)$
and  $\beta\in H_2(X)z\oplus H_0(X)x$. This proves the  Theorem.

\begin{Corollary}\label{cor:main3} 
\textup{(1)} Theorem \ref{main2} holds for any simply connected smooth  projective surface with 
$p_g=0$
and any good class $\xi$.

\noindent \textup{(2)} More generally for any smooth projective
surface $X$ and any $\xi\in \operatorname{Pic}(X)$, 
we have
\begin{equation*}
\delta_{\xi,t}^X(\exp(\alpha z+px))=
\begin{aligned}[t]
\sqrt{-1}^{\<\xi,K_X\>}\Big(\frac{q^{-\frac{1}{2}(\frac{\xi}{2})^2}}{\Lambda^{\chi(\oo_X)}}
\exp\Big(\frac{du}{da} \<\alpha,\xi/2\> z+T \<\alpha^2\> z^2 - ux\Big)
\qquad\qquad\\
\times
\Big(\frac{\sqrt{-1}}{\Lambda}
\frac{du}{da}\Big)^{2\chi(\oo_X)}\theta_{01}^\sigma\Big).
\end{aligned}
\end{equation*}
\end{Corollary}
\begin{proof}
In the notations of section 4, putting $t=2a$, we can rewrite Theorem \ref{univ} in terms of $q$. For $f,g\in \C((q^{1/8},\Lambda))$, we write $f\equiv g$ if $f/g=\exp(h)$ with 
$h\in q^{1/8}\C[[q^{1/2}]]$. Note that  $\frac{du}{da}\equiv \sqrt{-1}\Lambda q^{-1/8}$, 
$t\equiv\sqrt{-1}{\Lambda q^{-1/8}}$. Thus 
$$ (-1)^{\chi(\mathcal O_X)+\xi(\xi-K_X)/2} t^{-\xi^2-2\chi(\mathcal O_X)}
   \Lambda^{\xi^2+3\chi(\mathcal O_X)}\equiv
   \sqrt{-1}^{-\<\xi,K_X\>} \Lambda^{\chi(\oo_X)}q^{\frac{\xi^2}{8}}\Big(\frac{\sqrt{-1}}{\Lambda}\frac{du}{da}\Big)^{-2\chi(\oo_X)}.$$
 Thus for any triple $(X,\xi,\beta)$ with $v(X,\xi,\beta)=(v^1,\dots,v^8)$ we get 
 $$\delta_{\xi,t}^X(\beta)=\frac{\sqrt{-1}^{\<\xi,K_X\>}}{ \Lambda^{\chi(\oo_X)}}
 q^{-\frac{\xi^2}{8}}\Big(\frac{\sqrt{-1}}{\Lambda}\frac{du}{da}\Big)^{2\chi(\oo_X)}
 \exp\Big(\sum_{i=1}^8 v^i B_i \Big),$$ 
 for some universal power series $B_i\in \C((q^{-\frac{1}{8}}))[[\Lambda]]$.
 As  the $v(X,\xi,\beta)$ with $X$ a toric surface and $\xi$ a good class generate $\Q^8$ as a vector space, the
 $B_i$ are determined by their values for toric surfaces and good classes, i.e.\ they are given by 
 \eqref{modelt}. \begin{NB} changed 7.9. LG\end{NB}
 Note that the proof of \eqref{modelt} still works without any  changes also if $\xi$ is not good (replacing $\mathcal A_{\xi,-}$, $\mathcal A_{\xi,+}$ by
  $ -p_! (\mathcal I_2^\vee\otimes \mathcal I_1\otimes q^!\xi),$   
  $-p_! (\mathcal I_1^\vee\otimes \mathcal I_2\otimes q^!\xi^\vee)$.)
%
\end{proof}

\begin{Remark}\label{rem:c_1}

(1) 
Using \cite[Thm1.12]{Moc}, we get that  Thm \ref{main1} and part (1) of Cor. \ref{cor:main3} hold also if $\xi$ is not good. 
\begin{NB} added 6.9. LG\end{NB}

(2) As we mentioned in the introduction, the assertion that $\xi c_1(X)$
appears only as a sign in $\delta^X_{\xi,t}$ is one of statements of
the Kotschick-Morgan conjecture. This comes from $H(a,\Lambda) =
\pi\sqrt{-1}a$, as $\ve_1+\ve_2$ is the equivariant first Chern class
of $\mathbb A^2$. The latter statement, proved in
\cite[section~5.3]{NY2}, is a consequence of the blowup equation
\cite[(6.14)]{NY1}. This is by no means simple to check directly from
the definition of Nekrasov partition function.
\end{Remark}

\section{Equivariant Donaldson invariants for $\CP^2$}\label{P2}

Let us consider the complex projective plane $X = \CP^2$ and let $H$
be the hyperplane bundle. Let $M_H(n)$ be the moduli space of
$H$-semistable sheaves $E$ on $X$ with $\rank E = 2$, $c_1(E) = H$,
$c_2(E) - \frac{1}{4}c_1(E)^2 \equiv \Delta(E) = n$. As
$\operatorname{GCD}(2,c_1(E)) = 1$,
$M_H(n)$ is nonsingular of
dimension $4n-3$. Let $\mathcal E$ be the universal bundle.
Our method works also for $\rank E = 2$, $c_1(E) = 0$, $c_2(E) = n
\equiv 1 \mod 2$. But the moduli space becomes singular when
$c_2(E)$ is even, so our localization technique fails.

Let us consider the Donaldson invariants
\[
   \Phi_H^H(\alpha z +p x) = 
   \sum_{n\ge 0} \Lambda^{4n-3}
   \int_{M_H(n)} \exp\left(-\overline{\ch}_2(\mathcal
   E)/(\alpha z +p x)\right).
\]
Hereafter we denote this just by $\Phi(\alpha z +px)$ for brevity as we
will not vary $H$ in this section.

Let $\TT$ be the $2$-dimensional torus acting on $X$ by
$[x:y:z]\mapsto [t_1 x: t_2 y:z]$. We have three fixed points
$[1:0:0]$, $[0:1:0]$, $[0:0:1]$, and their characters of tangent
spaces are
\(
   1/t_1 + t_2/t_1
\),
\(
   t_1/t_2 + 1/t_2
\)
and
\(
   t_1 + t_2
\)
respectively.
We set $p_x = [1:0:0]$, $p_y = [0:1:0]$, $p_z = [0:0:1]$. We take the
coordinates around each $p_i$ and define their weights as 
\begin{equation*}
  (w(x_i),w(y_i)) = (-\ve_1,\ve_2-\ve_1),
\qquad
  (\ve_1-\ve_2,-\ve_2),
\qquad
  (\ve_1,\ve_2)
\end{equation*}
for $i=x,y,z$ respectively.
\begin{NB}
We could switch $w(x_i)$ and $w(y_i)$. But the partition function is
invariant under $\ve_1\leftrightarrow\ve_2$. So either choice is OK.
\end{NB}
We consider the induced $\TT$-action on $M_H(n)$. It also lifts to the
universal bundle $\mathcal E$, so we can define the equivariant
Donaldson invariants $\widetilde\Phi(\alpha z+px)$, where $\alpha$, $p$
are equivariant cohomology classes. In the nonequivariant limit
$\ve_1,\ve_2\to 0$, they go
to the ordinary Donaldson invariants
$\Phi(\alpha z+p x)$, where $\alpha$, $p$ are replaced by their
nonequivariant limit. For example, there are three equivariant lifts
$[p_x]$, $[p_y]$, $[p_z]$ of the point class $p$ given by the
three fixed points. Then $\widetilde\Phi(p_i x)$ depends on $i=x,y,z$,
but its nonequivariant limit $\lim_{\ve_1,\ve_2\to 0}
\widetilde\Phi(p_i x)$ is equal to $\Phi(p x)$.

\begin{Proposition}[\cite{Klyachko}]
\textup{(1)} A sheaf $E\in M_H(n)$ is fixed by the $\TT$-action if
  and only if there exists an $\TT$-equivariant structure on $E$.
  
\textup{(2)} A sheaf $E\in M_H(n)$ is fixed by the $\TT$-action if
and only if both its reflexive hull $E^{\vee\vee}$ and the quotient
$E^{\vee\vee}/E$ have $\TT$-equivariant structures.
\end{Proposition}
\begin{NB} 
  Although it is not explained carefully in \cite{Klyachko}, Kota
  pointed out me that the result is true when the moduli space is
  fine, but he does not know how to prove, say when $c_1 = 0$, $c_2$
  even. Apr. 17 HN
\end{NB}

For a stable sheaf $E$, its $\TT$-equivariant structure is unique up
to a twist by a character. We normalized it so that $\det
E^{\vee\vee}$ is trivial. This may not be possible in general, but it
is possible if we formally tensor by 
a square root of a line bundle. In
particular, the actions on the fibers over fixed points are
well-defined if we lift the action to a double covering
$\widetilde\TT\to \TT$. We consider the $\widetilde\TT$-structure as
if it is a $\TT$-structure hereafter.

Let $\shfO(x)$, $\shfO(y)$, $\shfO(z)$ be the $\TT$-equivariant line
bundles, where the $\TT$-structures are given so that the homomorphism
$x\colon \shfO\to\shfO(x)$ is equivariant, etc. The characters of the
fiber of $\shfO(x)$ at the fixed points $[1:0:0]$, $[0:1:0]$,
$[0:0:1]$ are given by
\(
   1, 
   t_2/t_1, 
   1/t_1
\)
respectively.
\begin{NB}
Let us check the equivariant Euler characteristic by Bott's formula:
\begin{equation*}
   \chi_T(\CP^2,\shfO(x)) = 
   \frac{1}{(1 - 1/t_1)(1 - t_2/t_1)}
   +
   \frac{t_2/t_1}{(1 - t_1/t_2)(1 - 1/t_2)}
   +
   \frac{1/t_1}{(1 - t_1)(1 - t_2)}
  = 1 + 1/t_1 + t_2/t_1.
\end{equation*}
\end{NB}

By a result of \cite{be,Klyachko}, we have
\begin{Proposition}
\textup{(1)} A $\TT$-equivariant rank $2$ vector bundle $E$ with $c_1(E)
= 1$ is classified by a triple $(p,q,r)\in \Z^3_{>0}$ with
$p+q+r\equiv 1\mod 2$.

\textup{(2)} The above $E$ is stable if and only if $p$, $q$, $r$
satisfy the strict triangle inequality, i.e.\  
$p+q < r$, $q+r < p$, $r+p < q$.
\end{Proposition}

In fact, the vector bundle $E$ is given as a cokernel of
\begin{equation*}
  {\cal O} \to {\cal O}(px) \oplus {\cal O}(qy) \oplus {\cal O}(rz)
\end{equation*}
for some $(p,q,r)\in\Z^3_{> 0}$ after a twist by a line bundle. Let
$E^{(p,q,r)}$ be the corresponding $\TT$-equivariant vector bundle. We
have
\begin{equation*}
   \Delta(E^{(p,q,r)}) 
    = (pq + qr + rp)/2 - p^2/4 - q^2/4 - r^2/4.
\end{equation*}
Let us denote this by $\Delta(p,q,r)$.
Note that $E^{(p,q,r)}$ is an isolated fixed point in
$M(H,\Delta(p,q,r))$. This assertion fails for higher ranks or toric
surfaces other than $\proj^2$. 

We have the decomposition of the fixed point set:
\begin{equation*}
  M_H(n)^{\TT} = \bigsqcup M_{p,q,r}(n - \Delta(p,q,r)),
\end{equation*}
where $M_{p,q,r}(m)$ denote the set of $\TT$-equivariant sheaves $E$
with $E^{\vee\vee} = E^{(p,q,r)}$ and
$\operatorname{length}(E^{\vee\vee}/E) = m$. The quotient sheaf
$E^{\vee\vee}/E$ is supported at $\{ p_x, p_y, p_z\}$. Accordingly we
have a factorization
\begin{equation*}
  M_{p,q,r}(m) 
   = \bigsqcup_{m_x+m_y+m_z = m}
      M_{p,q,r}^x(m_x)\times M_{p,q,r}^y(m_y) \times
    M_{p,q,r}^z(m_z),
\end{equation*}
where $M_{p,q,r}^x(m_x)$ denotes the set of $\TT$-equivariant sheaves
supported at $p_x$, etc.

The character of the fiber of $E^{(p,q,r)}$ at the fixed
point $p_z = [0:0:1]$ is given by
\begin{equation}\label{eq:weight_z}
\ch E^{(p,q,r)}_{p_z}
   = t_1^{p/2} t_2^{q/2} \left[t_1^{-p} + t_2^{-q}\right]
   = t_1^{-p/2} t_2^{q/2} + t_1^{p/2} t_2^{-q/2}
\end{equation}
Similarly the characters of the fibers at $p_y = [0:1:0]$, $p_x =
[1:0:0]$ are given by
\begin{align}
   & \ch E^{(p,q,r)}_{p_y}
   = (t_1/t_2)^{p/2} (1/t_2)^{r/2} \left[(t_1/t_2)^{-p} + (1/t_2)^{-r}\right]
   = t_1^{-p/2} t_2^{(p-r)/2} + t_1^{p/2} t_2^{(r-p)/2},
   \label{eq:weight_y}
\\
   & 
   \ch E^{(p,q,r)}_{p_x}
   = (1/t_1)^{r/2} (t_2/t_1)^{q/2} \left[(1/t_1)^{-r} + (t_2/t_1)^{-q}\right]
   = t_1^{(r-q)/2}t_2^{q/2} + t_1^{(q-r)/2} t_2^{-q/2},
   \label{eq:weight_x}
\end{align}
respectively.

Let us study a $\TT$-equivariant sheaf $E\in M_{p,q,r}^z(m_z)$,
i.e.\ $E^{\vee\vee} = E^{(p,q,r)}$ and $\Supp(E^{\vee\vee}/E)\linebreak[2] =
\{p_z\}$.
Using the 
coordinate system $(x/z, y/z)$ around $p_z$, we can identify
$E^{\vee\vee}/E$ with a $\TT$-equivariant quotient sheaf $Q =
\shfO^{\oplus 2}/F$, where $\TT$ acts on the trivial bundle
$\shfO^{\oplus 2}$ so that the character of the fiber at the origin is
\eqref{eq:weight_z}.

Let $M(n)$ be the framed moduli space of rank $2$ torsion-free sheaves
on $\P^2$ as in \subsecref{subsec:Nek}. This is the Gieseker
partial compactification of framed moduli spaces of instantons on
$\R^4$. Let $M_0(n)$ be the Uhlenbeck
partial compactification of framed moduli spaces of instantons on
$\R^4$, and let $\pi\colon M(n)\to M_0(n)$ be the natural projective
morphism. (See \cite[\S2]{NY1}, \cite[\S3]{NY2}.)
We have an action of $\hT = (\C^*)^2 \times \C^*$ on
$M(n)$, $M_0(n)$ such that $\pi$ is equivariant.
According to \eqref{eq:weight_z}, we define $\rho^z_{p,q,r}\colon
\TT\to \hT$ by
\begin{equation*}
  \rho^z_{p,q,r}(t_1,t_2) = (t_1,t_2 
  , t_1^{p/2}t_2^{-q/2}    
  ).
\end{equation*}
Note that there is no reason to prefer $t_1^{p/2} t_2^{-q/2}$ instead
of $t_1^{-p/2}t_2^{q/2}$. Either choice will work in the following
argument.

The following result follows from \cite{Furuta-Hashimoto}, but we give
a direct proof:
\begin{Lemma}\label{lem:fixed}
\textup{(1)} The origin $(\shfO^{\oplus 2}, m[0])$ is the only
$\TT$-fixed point in $M_0(m)$.

\textup{(2)} A point $(F,\varphi)\in M(m)$ is fixed by the
$\TT$-action if and only if $F^{\vee\vee} = \shfO^{\oplus 2}$ and
$F^{\vee\vee}/F$ is a $\TT$-equivariant sheaf.
\end{Lemma}

\begin{proof}
(2) follows from (1). Let us prove (1). Let us use the ADHM
description $(B_1,B_2,i,j)$ for $M_0(m)$. The coordinate ring of
$M_0(m)$ is generated by the following two types of functions
\begin{aenume}
\item $\tr(B_{\alpha_1}B_{\alpha_2}\dots B_{\alpha_N})$,
\item $\langle \chi, jB_{\alpha_1}B_{\alpha_2}\dots B_{\alpha_N} i\rangle$
\end{aenume}
where $\alpha_i = 1$ or $2$ and $\chi$ is a linear form on
$\End(W)$. Let us take a $\TT$-equivariant form $\chi$. Then its weight
is either $1$, $t_1^{-p} t_2^q$ or $t_1^p t_2^{-q}$. Under the
$\TT$-action, $B_{\alpha_1}B_{\alpha_2}\dots B_{\alpha_N}$ is
multiplied by
\(
   t_{\alpha_1} t_{\alpha_2} \dots t_{\alpha_N}.
\)
Therefore the first type of functions are never preserved by the
$\TT$-action. Similarly the second type of functions 
are multiplied by
\(
  t_{\alpha_1} t_{\alpha_2} \dots t_{\alpha_N} t_1 t_2,
\)
\(
  t_{\alpha_1} t_{\alpha_2} \dots t_{\alpha_N} t_1^{1-p} t_2^{q+1}
\)
or
\(
  t_{\alpha_1} t_{\alpha_2} \dots t_{\alpha_N} t_1^{p+1} t_2^{1-q}.
\)
These are never $1$ as $p$, $q > 0$.
\end{proof}

Thanks to this lemma we have
\begin{Corollary}
\begin{equation*}
   M^z_{p,q,r}(m) \cong M(m)^{\rho^z_{p,q,r}(\TT)}.
\end{equation*}
\end{Corollary}

We define $\rho^x_{p,q,r}$, $\rho^y_{p,q,r}\colon
\TT\to \hT$ by 
\begin{equation*}
  \begin{split}
   & \rho^x_{p,q,r}(t_1,t_2) = (1/t_1, t_2/t_1,
   t_1^{(q-r)/2} t_2^{-q/2}),
\\
   & \rho^y_{p,q,r}(t_1,t_2) = (t_1/t_2, 1/t_2,
   t_1^{p/2} t_2^{(r-p)/2}).
  \end{split}
\end{equation*}
(See \eqref{eq:weight_x} (resp.\ \eqref{eq:weight_y}.)
The above lemma and corollary hold also for these homomorphisms.

Let $N_{p,q,r;m}$ be the normal bundle of $M_{p,q,r}(m)$ in
$M(H,m+\Delta(p,q,r))$. Its fiber at $E$ is the sum of nonzero
weight spaces in $\Ext^1(E,E)$. We decompose $E^{\vee\vee}/E$ to
$Q^x$, $Q^y$, $Q^z$ according to the support $p_x$, $p_y$,
$p_z$. By the above corollary, we identify them with
$(F^x,\varphi)$, $(F^y, \varphi)$, $(F^z, \varphi)$ as elements of
$M(m_x)$, $M(m_y)$, $M(m_z)$. We have
\begin{equation*}
   \Ext^1(E,E) = \Ext^1(E^{(p,q,r)}, E^{(p,q,r)})
   + \sum_{i=x,y,z} \Ext^1(F^i, F^i(-\linf))
\end{equation*}
in the Grothendieck group of $\TT$-equivariant vector bundles on 
\(
   M_{p,q,r}(m)
   = \bigsqcup_{m=m_x+m_y+m_z}\linebreak[4]
   \prod_{i=x,y,z} M(m_i)^{\rho^i_{p,q,r}(\TT)}.
\)
The first factor of the right hand side is the tangent space of
$M_H(n)$ at $E^{(p,q,r)}$. Let us denote it by $T_{p,q,r}$. Then the
equivariant Euler class of $N_{p,q,r;m}$ is given by
\begin{equation*}
   e(T_{(p,q,r)})
   \prod_{i=x,y,z} e(N_{(F^i,\varphi)}),
\end{equation*}
where $N_{(F^i,\varphi)}$ denotes the fiber of the normal bundle of
the fixed point component containing $(F^i,\varphi)$ in $M(m_i)$.
We also have
\begin{equation*}
\begin{split}
   & - \overline{\ch}_2(\mathcal E)
   = - \overline{\ch}_2(E^{p,q,r}) + \ch_2(Q^x) + \ch_2(Q^y) +
   \ch_2(Q^z)
\\
   =\; & - \overline{\ch}_2(E^{p,q,r})
   + m_x [p_x] + m_y [p_y] + m_z [p_z],
\end{split}
\end{equation*}
where we have identified homology classes $[p_x]$, $[p_y]$, $[p_z]$
with their Poincar\'e dual.
We get
\begin{equation}\label{eq:a}
\begin{split}
   & \widetilde\Phi(\alpha z+px) =
   \sum_{p,q,r} \Lambda^{4\Delta(p,q,r)-3}
   \sum_{m} \Lambda^{4m} \int_{M_{p,q,r}(m)} 
   \frac{\exp\left(-\overline{\ch}_2(\mathcal
   E)/(\alpha z+px)\right)}{e(N_{p,q,r}(m))}
\\
   =\; & \sum_{p,q,r} \Lambda^{4\Delta(p,q,r) - 3}
   \frac{\exp\left(-\overline{\ch}_2(E^{(p,q,r)})/(\alpha z+px)\right)}
    {e(T_{p,q,r})}
    \\
   & \qquad\qquad\times
   \prod_{i=x,y,z}
   \sum_{m_i}
   \Lambda^{4m_i}
   \exp\left(m_i \iota_{p_i}^* (\alpha z+px)\right)
   \int_{M(m_i)^{\rho^i_{p,q,r}(\TT)}}
   \frac1{e(N_{(F^i,\varphi)})},
\end{split}
\end{equation}
where $\iota_{p_i}$ denotes the inclusion map $\{ p_i\} \to X$.

We study the first term and the second term separately.

\subsection{Quotient sheaf part}
Recall that the instanton part of Nekrasov's partition function is
written by
\(
  (\iota_{0*})^{-1}\pi_*[M(m)] = (\iota_{0*})^{-1} [M_0(m)],
\)
where $\iota_0$ is the inclusion
of the $\hT$-fixed point in $M_0(m)$. Here the equivariant
homology groups are taken with respect to the $\hT$-action. 
By \lemref{lem:fixed}(1) we replace them by those with respect to the
$\TT$-action and get an element
\(
  (\iota_{0*})^{-1} [M_0(m)]
\)
in the quotient field of $H^*_{\TT}(\mathrm{pt})$.
In order to distinguish this from the above element, we denote them by
\(
   (\iota_{0*})^{-1} [M_0(m)]_{\hT}
\)
and
\(
   (\iota_{0*})^{-1} [M_0(m)]_{\TT}.
\)

We set $S(\TT) = H^*_{\TT}(\mathrm{pt})$, $S(\hT) =
H^*_{\hT}(\mathrm{pt})$. We denote their quotient fields by $\mathcal
S(\TT)$ and $\mathcal S(\hT)$ respectively.
Let 
\(
  d\rho^i_{p,q,r}\colon
  \operatorname{Lie}(\TT)\to \operatorname{Lie}(\hT)
\)
be the differential of the homomorphism $\rho^i_{p,q,r}$.  It
induces the restriction homomorphism
\(
   \left(d\rho^i_{p,q,r}\right)^* \colon S(\hT) \to S(\TT).
\)

\begin{Lemma}\label{lem:specialization}
The rational function
\(
   (\iota_{0*})^{-1} [M_0(m)]_{\hT}\in   \mathcal S(\hT)
\)
can be restricted under the homomorphism $(d\rho^i_{p,q,r})^*$,
and is mapped to
\(
   (\iota_{0*})^{-1} [M_0(m)]_{\TT}\in
   \mathcal S(\TT).
\)
\end{Lemma}

\begin{proof}
{}From the proof of the localization theorem (see e.g., \cite{AB}), 
\(
   (\iota_{0*})^{-1} [M_0(m)]_{\hT}
\)
can be defined in a localized module $S(\hT)_f$ with a polynomial $f$
which vanishes on all Lie subalgebras of stabilizer subgroups $\neq \hT$.
Under the homomorphism
\(
   \rho^i_{p,q,r}\colon \TT\to \hT,
\)
stabilizer subgroups in $\hT$ are mapped to stabilizer subgroups in $\TT$.
By \lemref{lem:fixed}(1), if a stabilizer subgroup is not $\hT$, then
it is mapped to a subgroup $\neq\TT$. Therefore $f$ is restricted to a
nonzero polynomial under
\(
  d\rho^i_{p,q,r}\colon
  \operatorname{Lie}(\TT)\to \operatorname{Lie}(\hT)
\)
and we have an induced homomorphism
\begin{equation*}
   \left(d\rho^i_{p,q,r}\right)^* \colon
   S(\hT)_f \to S(\TT)_f.
\end{equation*}
{}From the definition we clearly have the assertion.
\end{proof}

By the localization theorem, 
\(
   (\iota_{0*})^{-1}[M_0(m_i)]_{\TT}
   = (\iota_{0*})^{-1}\pi_*[M(m_i)]_{\TT}
\)
is equal to
\begin{equation*}
   \int_{M(m_i)^{\rho^i_{p,q,r}(\TT)}}
   \frac1{e(N_{(F^\alpha,\varphi)})}.
\end{equation*}
Therefore we get
\begin{equation}\label{eq:quot}
\begin{split}
    & \sum_{m_i} \q^{m_i}
    \exp\left(m_i \iota_{p_i}^* (\alpha+pX)\right)
    \int_{M(m_i)^{\rho^i_{p,q,r}(\TT)}}
    \frac1{e(N_{(F^i,\varphi)})}
\\
    =\; &
    \Zin(w(x_i),w(y_i),-\textstyle{\frac{\xi^i_{p,q,r}}2};
    \q e^{\iota_{p_i}^*(\alpha+pX)}),
\end{split}
\end{equation}
where
\begin{equation}\label{eq:xi}
   \xi_{p,q,r}^z = -p\ve_1 + q\ve_2, \qquad
   \xi_{p,q,r}^y = -p\ve_1 + (p-r)\ve_2, \qquad
   \xi_{p,q,r}^x = (r-q)\ve_1 + q\ve_2.
\end{equation}

\subsection{Vector bundle part}

Let us calculate
\begin{equation*}
   \ch T^{(p,q,r)} = \ch \Ext^1(E^{(p,q,r)}, E^{(p,q,r)})
   = \sum_{p=0}^2 (-1)^{p+1} \ch H^p(\proj^2,
     \mathcal{E}{nd}_0(E^{(p,q,r)}))
\end{equation*}
where $\mathcal{E}{nd}_0$ means the trace-free part. We calculate this
by the localization theorem, i.e.
\begin{equation*}
\begin{split}
  & \ch \Ext^1(E^{(p,q,r)},E^{(p,q,r)})
\\
  =\; & - \frac{\ch\mathcal{E}{nd}_0(E^{(p,q,r)}))_{[0:0:1]}}
  {(1-t_1)(1-t_2)}
  - \frac{\ch\mathcal{E}{nd}_0(E^{(p,q,r)}))_{[0:1:0]}}
  {(1-t_1/t_2)(1-1/t_2)}
  - \frac{\ch\mathcal{E}{nd}_0(E^{(p,q,r)}))_{[1:0:0]}}
  {(1-1/t_1)(1-t_2/t_1)}.
\end{split}
\end{equation*}
We have
\begin{equation*}
\begin{split}
 & \ch \mathcal{E}{nd}_0(E^{(p,q,r)}))_{[0:0:1]}
   = 1 + t_1^{-p}t_2^q + t_1^p t_2^{-q}, \\
 & \ch \mathcal{E}{nd}_0(E^{(p,q,r)}))_{[0:1:0]}
   = 1 + t_1^{-p}t_2^{p-r} + t_1^{p} t_2^{r-p},\\
 & \ch \mathcal{E}{nd}_0(E^{(p,q,r)}))_{[1:0:0]}
   = 1 + t_1^{r-q}t_2^q + t_1^{q-r} t_2^{-q}.
\end{split}
\end{equation*}

A calculation shows the following:
\begin{Lemma}
Let us define the convex region $D^{(p,q,r)}$ as follows:

\textup{(1)} Case $p=q=r=1$: $D^{(p,q,r)} = \{ (0,0) \}$.

\textup{(2)} Case $p=1$, $q=r\neq 1$: 
\(
   D^{(p,q,r)} = \operatorname{Conv}((0,q-1), (-1,q-1), (-1,-q+2),
   (0,-q+1)).
\)

\textup{(3)} Case $q=1$, $p=r\neq 1$:
\(
   D^{(p,q,r)} = \operatorname{Conv}((p-1,0), (-p+1,0), (-p+2,-1),
   (p-1,-1)).
\)

\textup{(4)} Case $r=1$, $p=q\neq 1$:
\(
   D^{(p,q,r)} = \operatorname{Conv}((p-1,1-p), (p-1,2-p), (2-p,p-1),
   (1-p,p-1)).
\)

\textup{(5)} Case $p+q=r+1$, not above:
\(
   D^{(p,q,r)} = \operatorname{Conv}((p-1,q-1), (-p,q-1), (-p,-q+2),
   (-p+2,-q), (p-1,-q)).
\)

\textup{(6)} Case $r+p=q+1$, not above:
\(
   D^{(p,q,r)} = \operatorname{Conv}((p-1,1-q), (p-1,r-p+1), (2-p,q-1),
   (-p,q-1), (-p,p-r+1)).
\)

\textup{(7)} Case $q+r=p+1$, not above:
\(
   D^{(p,q,r)} = \operatorname{Conv}((p-1,-q), (p-1,r-p+1), (r-q+1,q-1),
   (1-p,q-1), (q-r+1,-q)).
\)

\textup{(8)} Otherwise:
\(
   D^{(p,q,r)} = \operatorname{Conv}((p-1,-q), (p-1,r-p+1),
   (r-q+1,q-1), (-p,q-1), (-p,p-r+1), (q-r+1,-q))
\)

Here $\operatorname{Conv}$ denotes the convex hull. 
Then $\ch \Ext^1(E^{(p,q,r)},E^{(p,q,r)})$ is the sum of monomials $t_1^m
t_2^n$ where $(m,n)\in\Z^2$ runs over $D^{(p,q,r)}\setminus\{(0,0)\}$.
\end{Lemma}

Note that the origin $(0,0)$ is in $D^{(p,q,r)}$ in all
cases. We thus have
\begin{equation}\label{eq:vect}
    e(T_{p,q,r}) = \prod_{(m,n)\in D^{(p,q,r)}\cap\Z^2\setminus\{(0,0)\}}
    (m\ve_1 + n\ve_2).
\end{equation}

We can also express $\ch_2(\mathcal E)/(\alpha z +px )$ by the
localization formula:
\begin{equation}\label{eq:mu_vect}
   -\ch_2(\mathcal E)/(\alpha z +p x)
   = -\frac14 \sum_{i=x,y,z} 
   \frac{(\xi_{p,q,r}^i)^2 \iota_{p_i}^*(\alpha z +p x)}{w(x_i)w(y_i)},
\end{equation}
where $\xi_{p,q,r}^i$ is as in \eqref{eq:xi} and $w(x_i)w(y_i)$
appears as the Euler class $e(T_{p_i}\proj^2)$ of the tangent space at
$p_i$.

Substituting (\ref{eq:quot}, \ref{eq:vect}, \ref{eq:mu_vect}) into
\eqref{eq:a}, we get the following:
\begin{Theorem}\label{thm:P2}
The equivariant Donaldson invariants of $\proj^2$ are given by
\begin{equation*}
  \begin{split}
  \widetilde\Phi(\alpha z+px) &=
  \sum_{p,q,r} \Lambda^{4\Delta(p,q,r)-3}
  \exp\left(-\frac14 \sum_{i=x,y,z} 
   \frac{(\xi_{p,q,r}^i)^2 \iota_{p_i}^*(\alpha z+px)}
   {w(x_i)w(y_i)}\right)
\\
  & \qquad\times
    {\prod_{(m,n)\in \substack{D^{(p,q,r)}\cap\Z^2\\
       \ \ \setminus\{(0,0)\}}}
      \frac1{m\ve_1 + n\ve_2}}
   \prod_{i=x,y,z}
    \Zin(w(x_i),w(y_i),-\textstyle{\frac{\xi^i_{p,q,r}}2};
    \q e^{\iota_{p_i}^*(\alpha z+px)}),
  \end{split}
\end{equation*}
where $p,q,r$ runs over $\Z^3_{>0}$ satisfying
$p+q+r\equiv 1\mod 2$ and the strict triangle inequality.
\end{Theorem}

Ordinary Donaldson invariants $\Phi(\alpha z+px)$ are given by
$\lim_{\ve_1,\ve_2\to 0} \widetilde\Phi(\alpha z+px)$. But the
solution of Nekrasov's conjecture does not say anything about this
limit, so we do not know how to get an explicit formula from the above.
Note that the summation over $p$,$q$,$r$ is related to 
Hurwitz class numbers (according to \cite{Klyachko1}), 
which appeared in \cite[\S9]{MW} in the formula for the Donaldson
invariants of $\proj^2$.
 
On the other hand we have
\begin{Theorem}
Let $P\colon Y\to \P^2$ be the blowup of the fixed point $p_z$
\textup(different from $p=p_x$\textup). Then
\begin{equation*}
  \widetilde\Phi(H z + p x) = \Bigg[
   \begin{aligned}[t]
   & \frac{1}{\Lambda}
   \exp\Big(\frac{1}{2}\big\langle\Todd_2(Y)(P^* H z+ P^*p
   x)\big\rangle 
\\
   & \quad
   \times\sum_{\substack{\xi = (2n-1)P^* H - 2a E\\ a\ge n\in\Z_{>0}}}
   \exp\left(
     \sum_{i=x,y,z_1,z_2}
   F\big(w(x_i),w(y_i),\hbox{$\frac{t-\iota_{p_i}^*\xi}{2}$};\Lambda
   e^{\iota_{p_i}^*(P^*H z+P^*px)/4}\big)
   \right)
   \Bigg]_{t^{-1}},
   \end{aligned}
\end{equation*}
where $p_{z_1}$, $p_{z_2}$ are the fixed
points in the exceptional set of $Y$.
\end{Theorem}

The formula, when compared with the one in \thmref{thm:P2}, probably
gives us a nontrivial identity on the partition function.

The idea of the proof is the same as in \cite[Th.\ 3.5]{G}, but we put
a little more care as we consider the equivariant Donaldson invariants.

\begin{proof}
Let $P\colon Y\to \P^2$ be the blowup of the fixed point $p_z$.
We first assume that the line $H$ is $H_{xy}$. In particular,
$H$ does not pass through the point $p_z$ which we blowup.
Let $\widehat M_{H}(n)$ be
the moduli space of $(P^*H-\ve E)$-stable rank $2$ sheaves on $Y$ with
$c_1 = P^* H$, $\Delta = n$. By \cite[App.~F]{NY2} there exists a
projective morphism $\widehat\pi\colon \widehat M_{H}(n)\to N_{H}(n)$,
where $N_H^{\P^2}(n)$ is the Uhlenbeck compactification of the moduli
space of locally free sheaves on $\P^2$ with $c_1 = H$, $\Delta =
n$. 

\begin{NB}
H.N. Apr. 17

Kota give two proof for the assertion that we can put a $T$-action on
$N_H(n)$ so that $\widehat\pi$ is $T$-equivariant by the definition of
$\widehat\pi$ in \cite{NY2}.

{\bf 1st proof}:

For the $T$-action on $N_H(n)$, the invariance of $W^a_n$
follows from its definition.
$W^a_n$ is generated by sections corresponding Cartier divisors
$$
  D_L:=\{E \in M_H(n)|H^1(E \otimes L) \ne 0 \},
$$
where $L$ is a line bundle of a suitable degree (we choose $L$ so that 
$\chi(E\otimes L)=0$)
on a smooth curve $\not \ni p_z$ (the center of the blow-up).

This condition is clearly invariant under the action of $T$:
for $t \in T$,
$$
t^*(D_L):=\{E \in M_H(n)|H^1(t_*(E) \otimes L) \ne 0 \}
         =\{E \in M_H(n)|H^1(E \otimes t^*L) \ne 0 \}=D_{t^*(L)}.
$$
So we have an action of $T$ on $N_H(n)$ which is compatible
with the action on $M_H(n)$. This is enough for our purpose:

(Since $M_H(n) \to N_H(n)$ is surjective, the compatible action on $N_H(n)$ is
set-theoretically unique. Since we regard $N_H(n)$ as a reduced scheme,
the compatible action is really unique.
Thus the action coincides with the naive action.)

For the case of $\widehat{\pi}:\bM_H(n) \to N_H(n)$,
it is sufficient to see that
the map is equivariant for all $t \in T$ on an open subset (depending on t)
of $\bM_H(n)$, but this is obvious.

{\bf 2nd proof}:

Notation:
Let $G$ be a vector bundle on $X$.
A coherent sheaf $E$ is $G$-twisted semi-stable, iff
$\chi(G^{\vee} \otimes F(nH))/\rank F \leq
\chi(G^{\vee} \otimes E(nH))/\rank E$, $n \gg 0$ for all subsheaf
$F$ of $E$.
We denote the moduli space by $M_H^G([E])$.

We set $E_0:=\Omega_{{\Bbb P}^2}(2)$.
Then $E_0$ is a stable vector bundle of rank 2 such that
$c_1(E_0)=H$ and $\chi(E_0,E_0)=1$.
In my paper (math.AG/0106042), I gave a moduli-theoretic
construction for the contraction
map $\pi$ in the case of the moduli of stable sheaves of rank 2 
and $c_1=H$ on 
${\Bbb P}^2$ (see Example 2.1 in the paper).
Unfortunately the construction depends on the property
$\chi(E_0,E_0)=1$, I can only treat this spacial case.

Roughly speaking, the contraction map is given by
sending $E$ to the $S$-equivalence class $gr(\tilde{E})$ of 
the universal extension
$$
0 \to E \to \tilde{E} \to E_0 \otimes \Ext^1(E_0,E) \to 0.
$$
(Note that $\tilde{E}$ is $E_0$-twisted semi-stable).
Therefore the Uhlenbeck compactification is
the image of $M_H([E])=M_H^{E_0}([E]) \to M_H^{E_0}([\tilde{E}])$.
The same construction works for $M_{P^*(H)}^{P^*(E_0)}(P^*[E])=
M_{P^*(H)-\ve E}^{P^*(E_0)}(P^*[E])$.

In the following, I will construct $\widehat{\pi}$
by constructing $\widehat{\pi}:M_{P^*(H)-\ve C}^{P^*(E_0)}(P^*[\tilde{E}])\to
 M_H^{E_0}([\tilde{E}])$:
 \begin{equation}
 \begin{CD}
 M_{P^*(H)}^{P^*(E_0)}(P^*[E]) @>{\pi}>> 
 M_{P^*(H)-\ve C}^{P^*(E_0)}(P^*[\tilde{E}])\\
@. @VV{\widehat{\pi}}V \\
M_{H}^{E_0}([E]) @>{\pi}>> M_{H}^{E_0}([\tilde{E}])
\end{CD}
\end{equation}

Since the construction is moduli-theoretic, it is easy to see that
the Uhlenbeck compactification equips $T$-action and 
$\pi,\widehat{\pi}$ are $T$-equivariant.

Let $F$ be a $P^*(E_0)$-twisted semi-stable sheaf on $\widehat{\Bbb P}^2$
with respect to $P^* H-\ve C$ such that
$c_1(P^*(E_0)^{\vee} \otimes F)=0$ and
$\chi(P^*(E_0),F)=0$, where $C$ is the exceptional divisor. 
We first note that $P_*(F)^{\vee \vee}$ is 
$\mu$-semi-stable with respect to $H$.
Hence 
\begin{equation}\label{eq:ext2}
\Ext^2(P^*(E_0),F(kC))=\Hom(F,P^*(E_0)(K_{\widehat{\Bbb P}^2}-kC))^{\vee} 
\subset \Hom(P_*(F)^{\vee \vee},E_0(K_{{\Bbb P}^2}))^{\vee}=0
\end{equation} 
for all $k$.
\begin{Lemma}\label{lem:vanish}
\textup{(1)} $\Hom(P^*(E_0),F)=\Ext^1(P^*(E_0),F)=0$.

\textup{(2)} $H^1(F \otimes {\cal O}_C)=0$.
\end{Lemma}

\begin{proof}
(1) Since $\chi(P^*(E_0),F)=0$, it is sufficient to prove 
$\Hom(P^*(E_0),F)=0$.
Since $P^*(E_0)$ is a $P^*(E_0)$-twisted stable sheaf
with $\chi(P^*(E_0),P^*(E_0))=1$,
there is no homomorphism $P^*(E_0) \to F$.

(2) It is sufficient to prove 
$\Ext^1(P^*(E_0),F \otimes {\cal O}_C)=0$.
This follows from (1) and \eqref{eq:ext2}.
\end{proof}

\begin{Corollary}
$R^1 P_*(F)=0$ and $P_*(F)$ is an $E_0$-twisted semi-stable
sheaf with $\chi(E_0,P_*(F))=0$.
\end{Corollary}

\begin{proof}
We have an exact sequence
\begin{equation}
0 \to {\cal O}_C (n) \to {\cal O}_{(n+1)C} \to {\cal O}_{nC}
\to 0.
\end{equation}
Hence Lemma \ref{lem:vanish} (2) implies that 
$H^1(F \otimes {\cal O}_{nC})=0$ for all $n>0$.
Hence $\varprojlim H^1(F \otimes {\cal O}_{nC})=0$.
By the theorem of formal function, 
$\varprojlim R^1 P_*(F) \otimes {\cal O}_{{\Bbb P}^2}/{\frak m}_{p_z}^n
=\varprojlim H^1(F \otimes {\cal O}_{nC})=0$.
Thus $R^1 P_*(F)=0$.
We prove the second claim.
The claim $\chi(E_0,P_*(F))=0$ is obvious.
Assume that $P_*(F)$ is not $E_0$-twisted semi-stable.
Then there is a $E_0$-twisted semi-stable subsheaf $G$ such that
$c_1(G \otimes E_0^{\vee})=0$ and 
$\chi(E_0,G)>0$.
Since $\Ext^2(E_0,G)=\Hom(G,E_0(K_{{\Bbb P}^2}))^{\vee}=0$,
there is a non-trivial homomorphism 
$E_0 \to G \to P_*(F)$.
On the other hand, 
$0=\Hom(P^*(E_0),F)=\Hom(E_0,P_*(F))$, this is a contradiction.
\end{proof}

Therefore we have a contraction map
\begin{equation}
\begin{matrix}
\widehat{\pi}:& M_{P^*(H)-\ve C}^{P^*(E_0)}([F])&  \to & M_H^{E_0}([P_*(F)])\\
& F & \mapsto & P_*(F).
\end{matrix}
\end{equation}
\end{NB}

By the definition of $\widehat\pi$ the class $\mu(P^*H)$ on
$\widehat M_{H}(n)$ is the pullback of the class $\mu(H)$ on $N_H(n)$
by $\widehat\pi$.
In fact, by
\begin{equation*}
   H_{xz} = H_{xy} - \ve_2, \qquad
   H_{yz} = H_{xy} - \ve_1,
\end{equation*}
and
\begin{equation*}
  c_2(\mathcal E) - \frac14 c_1(\mathcal E)^2/[\proj^2]
  \in H^0_{\TT}(N_H(n)) \cong \C,
\end{equation*}
$\mu(H_{xz})$, $\mu(H_{yz})$ are equal to $\mu(H) = \mu(H_{xy})$ modulo
classes from $H^*_{\TT}(\mathrm{pt})$. Therefore this assertion is true
for any $H$.

By \cite[Th.~6.9]{FM}, $\mu(p)$ extends to a class on the Donaldson
compactification $N_H(n)$. The extension can be made so that the class
is equivariant with respect to the compact form of $\TT$, and it is
enough for our purpose.
\begin{NB}
  The proof in \cite{FM} is topological. I do not know how to prove
  the same assertion for the (equivariant) Chow group.
\end{NB}
Then we have $\mu(P^*p) = \widehat\pi^*\mu(p)$ as we blowup at a point
different from $p$.
Therefore
\begin{equation*}
\begin{split}
   & \exp(\mu(P^*H) z + \mu(P^*p) x)\cap [\widehat M_H(n)]
   = \widehat\pi^*\left(\exp(\mu(H)z + \mu(p)x)\right)
   \cap [\widehat M_H(n)]
\\
   =\; & \exp(\mu(H)z + \mu(p)x)\cap \widehat\pi_*[\widehat M_H(n)]
   = \exp(\mu(H)z + \mu(p)x)\cap [N_H(n)].
\end{split}
\end{equation*}
There is an alternative way to prove this formula. Restrict the maps
$\pi$, $\widehat\pi$ to the fixed point set: $\pi\colon
M_H(n)^{\TT}\to N_H(n)^{\TT}$, $\widehat\pi\colon \widehat
M_H(n)^{\TT}\to N_H(n)^{\TT}$. We have
\begin{equation*}
   N_H(n)^{\TT} = 
 \bigsqcup_{\substack{p,q,r,m_x,m_y,m_z \\ m_x+m_y+m_z+\Delta(p,q,r) = n}}
   \{ (E^{(p,q,r)}, m_x [p_x] + m_y [p_y] + m_z [p_z]) \}
\end{equation*}
by the same argument as above. In particular, $N_H(n)^{\TT}$ consists
of finitely many points. We have direct sum decompositions
of $H^4_{\TT}(M_H(n))$ and $H^4_{\TT}(\widehat M_H(n))$ correspondingly.
{}From the expression \eqref{FY} we see that
$\mu(p)\in H^4_{\TT}(M_H(n)^{\TT})$ and
$\mu(P^* p)\in H^4_{\TT}(\widehat M_H(n)^{\TT})$ are
pullbacks of the same class in
\( 
   H^4_{\TT}(N_H(n)^{\TT}) 
  = \bigoplus_{p,q,r,m_x,m_y,m_z}H^4_{\TT}(\mathrm{pt}).
\)
This assertion is enough for the above calculation.

Therefore
\[
   \widetilde\Phi^{\P^2,H}_H(\exp(H z + p x))
  = \widetilde\Phi^{Y,H-\ve E}_H(\exp(P^*H z + P^*p x)).
\]

On the other hand we have
\begin{equation*}
   \widetilde\Phi^{Y,F+\ve E}_H(\exp(P^*H z + P^*p x)) = 0
\end{equation*}
for $F = P^*H - E$ is the fiber class and $\ve$ is a sufficiently
small number by \cite{Q}. Therefore by the proof of \cite[Th.~3.5]{G}
we have
\begin{equation*}
   \widetilde\Phi^{Y,H-\ve E}_H(\exp(P^*H z + P^*p x)) 
   = \sum_{\substack{\xi = (2n-1)P^* H - 2a E\\ a\ge n\in\Z_{>0}}}
   \widetilde\delta^{Y}_\xi(\exp(P^*H z+ P^*p x)).
\end{equation*}

Let $p_x$, $p_y$ denote the inverse image of $p_x$, $p_y$ under
$P$. Let $p_{z_1}$, $p_{z_2}$ be the two fixed points in the
exceptional set $E$. By \corref{deltaq} we have
\begin{equation*}
\begin{split}
   & \widetilde\delta^{Y}_{\xi,t}(\exp(P^*H z+ P^*p x))
\\
   = \; &
   \begin{aligned}[t]
   \frac{1}{\Lambda}
   \exp\Big(&\frac{1}{2}\big\langle\Todd_2(Y)(P^* H z+ P^*p x)\big\rangle 
   \\
   &\qquad\qquad\times\sum_{i=x,y,z_1,z_2}
   F\big(w(x_i),w(y_i),\hbox{$\frac{t-\iota_{p_i}^*\xi}{2}$};\Lambda
   e^{\iota_{p_i}^*(P^*H z+P^*px)/4}\big)\Big).
   \end{aligned}
\end{split}
\end{equation*}
Now the assertion follows.
\end{proof}

\begin{NB}
Note that $\iota_{p_i}^* P^*H = 0 = \iota_{p_i}^* P^* p$ for $i=z_1$,
$z_2$. Therefore we have
\begin{equation*}
\begin{split}
  & \sum_{i=z_1,z_2}
   F\big(w(x_i),w(y_i),\hbox{$\frac{t-\iota_{p_i}^*\xi}{2}$};\Lambda
   e^{\iota_{p_i}^*(P^*H z+P^*px)/4}\big)
\\
  =\; &
   F(\ve_1,\ve_2-\ve_1,
   \textstyle{\frac{t}2}+a\ve_1;\Lambda)
   + 
   F(\ve_1-\ve_2,\ve_2, \textstyle{\frac{t}2}+a\ve_2;\Lambda).
\end{split}
\end{equation*}
On the other hand we have $\iota_{p_i}^*\xi = 0$ for $i=x,y$. Then
\begin{equation*}
\begin{split}
  &\sum_{i=x,y}
  F\big(w(x_i),w(y_i),\hbox{$\frac{t-\iota_{p_i}^*\xi}{2}$};\Lambda
   e^{\iota_{p_i}^*(P^*H z+P^*px)/4}\big)
\\
   = \; &
   F\left(-\ve_1,\ve_2-\ve_1,
   \textstyle{\frac{t}2+(2n-1)\ve_1};\Lambda
   e^{(\ve_1 z + \ve_1(\ve_1-\ve_2)x)/4}\right)
   + F\left(\ve_1-\ve_2,-\ve_2,
   \textstyle{\frac{t}2+(2n-1)\ve_2};\Lambda e^{\ve_2 z/4}\right).
\end{split}
\end{equation*}
\end{NB}

\appendix
\section{Generic smoothness after blowup}\label{sec:app}

For a $\mu$-semistable rank $2$ sheaf $F$ on $X$, there exists a
constant $\beta_\infty$ depending only on $X$ (and the rank of $F$)
(see \cite[4.5.7]{HL}) such that
\[
  \dim \Ext^2(F,F)_0 \leq \beta_\infty.
\]
Therefore
\begin{equation}
\label{eq:beta}
  \dim M_H^X(c_1,n) \leq \exp\dim M_H^X(c_1,n) +\beta_\infty,
\end{equation}
where $\exp\dim M_H^X(c_1,n)$ is the expected dimension of $M_H^X(c_1,n)$.

By the result of Donaldson, Zuo, Gieseker-Li, O'Grady (see
\cite[\S9]{HL}) there exists a constant $m_0$ depending only on $X$,
$H$ (and rank) such that $M_H^X(c_1,m)$ is irreducible and of expected
dimension for $m\ge m_0$.

Let $P\colon \widehat X\to X$ be blowup at points $p_1$, \dots, $p_N$
as before. We take a polarization $H$ on $X$ and consider the
polarization $P^*H$ on $\widehat X$ as above. For simplicity we assume
$(c_1,H)$ is odd. By \cite[App.~F]{NY2} we have a projective morphism
\(
  \widehat\pi\colon M^{\widehat X}_{P^*H}(P^* c_1,m) \to N^X_H(c_1,m),
\)
where $N^X_H(c_1,m)$ is the Uhlenbeck compactification, which is
set-theoretically equal to
\(
  N^X_H(c_1,m) = \bigsqcup_k M_H^X(c_1,m-k)_{\mathrm{lf}} \times S^k X,
\)
where $M_H^X(c_1,m-k)_{\mathrm{lf}}$ is the open subscheme
of $M_H^X(c_1,m-k)$ consisting of stable vector bundles.

A point in $S^k X$ can be written as $[Z] = \sum m_i [p_i] + \sum
\lambda_p [x_p]$ where $p_i$, $x_p$ are disjoint and $\lambda_p\ge 1$.
Then we have a stratification of $S^k X$ parametrized by $(m_i)_i\in
\Z_{\ge 0}^N$ and the partition $\lambda = \{ \lambda_p \}_p$ of $k -
\sum m_i$.
By \cite[App.~F]{NY2} the fiber of $\widehat\pi$ over $(E, [Z])\in
M_H^X(c_1,m-k)_{\mathrm{lf}} \times S^k X$ depends only on $m-k$ and
the stratum containing $[Z]$. And it is also equal to the fiber of the
morphism defined for the framed moduli spaces on $\bp$ and
$\proj^2$. The homology of central fibers (i.e.\ $\lambda = \emptyset$,
$m_1 = n$, $m_i = 0$ ($i\ge 2$)) was calculated in
\cite[Th.~3.8$\sim$10]{NY2}. We find its dimension is given by
\begin{equation*}
   2n + \max_{l\in \Z : l^2\le n} l \le 3n.
\end{equation*}
Therefore
\begin{equation*}
   \dim \widehat\pi^{-1}(E,[Z])
   \le 3 \sum_{i=1}^N m_i + \sum_p (2\lambda_p - 1).
\end{equation*}
Therefore we have
\begin{equation*}
\begin{split}
  & \dim \widehat\pi^{-1}(M_H^X(c_1,m-k)_{\mathrm{lf}}\times S^k X)
\\
   \le \; & 
   \dim M_H^X(c_1,m-k)_{\mathrm{lf}} + 
   \max_{\sum m_i + |\lambda_p| = k}\left\{ \sum_i 3 m_i
   + \sum_p (2 \lambda_p + 1)\right\}
\\
   \le\; & \dim M_H^X(c_1,m-k)_{\mathrm{lf}} + 3k.
\end{split}
\end{equation*}
Let us take $m \ge m_0 + \beta_\infty$. For $k > \beta_\infty$ we
have
\begin{equation*}
   \dim M_H^X(c_1,m-k)_{\mathrm{lf}} + 3k
   \le \exp\dim M_H^X(c_1,m) - k + \beta_\infty
   < \exp\dim M_H^X(c_1,m)
\end{equation*}
by \eqref{eq:beta}.
For $k \le \beta_\infty$, we have $m - k \ge m_0$. Therefore
$M_H^X(c_1,m-k)$ is of expected dimension. Therefore
\begin{equation*}
   \dim M_H^X(c_1,m-k)_{\mathrm{lf}} + 3k
   = \exp\dim M_H^X(c_1,m) - k < \dim M_H^X(c_1,m)
\end{equation*}
unless $k= 0$.
The open locus $\widehat\pi^{-1}(M_H^X(c_1,m)_{\mathrm{lf}})$ consists
of pullbacks $P^*E$ of $E\in M_H^X(c_1,m)$ and $\widehat\pi$ is an
isomorphism there. Therefore $M^{P^*H}_{\widehat X}(P^*c_1,m)$ is
of expected dimension (and irreducible).

\end{document}